\documentclass[11pt,reqno]{amsart}
\usepackage{amsmath}
\numberwithin{equation}{section}
\usepackage{hyperref}
\hypersetup{nesting=true,debug=true,naturalnames=true}
\usepackage{graphicx,amssymb,upref}
\usepackage{enumerate}
\usepackage[all]{xy}
\usepackage{color}
\usepackage{mathrsfs}
\usepackage{amsfonts}
\usepackage{caption}
\usepackage{paralist}
\usepackage{color}
\usepackage{float}
\usepackage{tikz,pgf}
\usetikzlibrary{calc}
\usetikzlibrary{intersections,decorations.markings}
\usepackage{tikz-cd}

\theoremstyle{definition}
\newtheorem{theorem}{\bf Theorem}[section]

\newtheorem{lemma}[theorem]{\bf Lemma}
\newtheorem{corollary}[theorem]{\bf Corollary}

\theoremstyle{definition}

\newtheorem{definition}[theorem]{\bf Definition}
\newtheorem{remark}[theorem]{\bf Remark}

\newtheorem{proposition}[theorem]{\bf Proposition}
\newcommand{\mm}[1]{\mathrm{#1}}
\newcommand{\mb}[1]{\mathbb{#1}}
\newcommand{\CP}{\mathbb {CP}}
\newcommand{\Sq}{\mathrm {Sq}}
\newcommand{\K}{\mathrm {K}}
\newcommand{\Z}{\mathbb {Z}}
\newcommand{\SU}{\mathrm {SU}}
\newcommand{\Sp}{\mathrm {Sp}}
\newcommand{\MO}{\mathrm {MO}}
\newcommand{\diag}{\mathrm {diag}}

\newcommand{\q}{\bold q}
\newcommand{\BO}{\mathrm{BO}}
\newcommand{\Tor}{\operatorname{Tor}}
\newcommand{\Ext}{\operatorname{Ext}}

\newcommand{\bord}{\Omega^{\mathrm{O}\langle8\rangle}}
\newcommand{\M}{\mathrm{M}}

\makeatletter
\@namedef{subjclassname@2020}{\textup{2020} Mathematics Subject Classification}
\makeatother

\begin{document}

\title[Topological classification of Bazaikin spaces]
{Topological classification of Bazaikin spaces}

\author[F. Fang]{Fuquan Fang}
\thanks{The first author was supported by
NSF Research Group Grant 11821101 of China and National Center of Applied Mathematics, Shenzhen}
\address{Department of Mathematics, Capital Normal University,
Beijing, P.R.China}
  \email{ffang@nankai.edu.cn}
\author[Wen Shen]{Wen Shen}
\address{Department of Mathematics, Wenzhou University,
Wenzhou, P.R.China}
\email{shenwen121212@163.com}

\begin{abstract}
 In this paper, we establish a topological classification of Bazaikin spaces whose sixth cohomology groups have orders prime to \(5\). Moreover, we obtain an estimate for the number of diffeomorphism types among these spaces.
\end{abstract}

\subjclass[2020]{Primary 57R65, 57R90, 55T10, 55T25}

\maketitle
\tableofcontents

\section{Background and Statements of the Main Theorems}\label{introduction}

 Manifolds with positive sectional curvature have been of interest since the beginning of Riemannian geometry. Except for the rank one symmetric spaces, there are very few examples of compact manifolds with positive sectional curvature, which exist only in dimensions $6$, $7$, $12$, $13$, and $24$ as identified by Wallach \cite{Wallach} and Berger \cite{Be}. Additionally, there are two infinite families: one in dimension $7$ (Aloff-Wallach spaces and Eschenburg spaces; see \cite{Aloff} \cite{E1} \cite{E2}) and the other in dimension $13$ (Bazaikin spaces; see \cite{Ba}).
The first family in dimension 7 consists of the Eschenburg biquotients
$$\mathcal{E}_{k,l}=diag(z^{k_1},z^{k_2},z^{k_3})\backslash \mm{SU}(3)/diag(z^{l_1},z^{l_2},z^{l_3})^{-1},$$
with $\Sigma k_i=\Sigma l_i$ which includes the infinite subfamily of homogeneous Aloff-Wallach spaces
$$ \mathcal{W}_l= \mm{SU}(3)/diag(z^{l_1},z^{l_2},z^{l_3})\quad \Sigma l_i=0.$$

Based on the modified surgery theory developed by Kreck \cite{K}, Kreck and Stolz \cite{KS1,KS} obtained a complete topological classification of $7$-manifolds modeled on the Aloff-Wallach spaces. As an important by-product,  it was shown that the space of positively curved metrics on a certain Aloff-Wallach space is not path-connected. Further classification results on the Eschenburg spaces $\mathcal{E}$ were obtained in \cite{kb1}, \cite{kb2}.

The second family exists in dimension $13$ and consists
of the Bazaikin biquotients (as described in \cite{E3})
\begin{equation}
	\mathcal{B}_\bold{q}=diag(z^{q_0},z^{q_1},z^{q_2},z^{q_3},z^{q_4},z^{q_5})\backslash \mm{SU}(6)/\mm{Sp}(3) \label{Bazaibiquotient}
\end{equation}
where $\bold{q}= (q_0, \cdots , q_5)$ is a $6$-tuple of odd integers satisfying $\Sigma q_i = 0$. This biquotient is a manifold if and only if
\begin{equation}
	\text{each $q_i$ is odd with $\Sigma q_i = 0$ and $gcd(q_{\tau(1)}+q_{\tau(2)},q_{\tau(3)}+q_{\tau(4)})=2$} \label{freeactionBq}
\end{equation}
for all permutations $\tau \in S_6$ \cite{FL}. We call $\mathcal{B}_\bold{q}$ a Bazaikin space when the biquotient is a manifold. Meanwhile, the left $S^1$-action in \eqref{Bazaibiquotient} on $\SU(6)/\Sp(3)$ has ineffective kernel $\Z_2$.

Moreover, by  \cite{FL}  a Bazaikin space $\mathcal{B}_\bold{q}$ with a certain deformed biquotient metric has positive sectional curvature  if there exists $0 \le i_0 \le 5$ such that
\begin{equation}
	\text{  $q_i + q_j > 0$ (or $< 0$), for all $i_0 \ne i < j \ne i_0$ \label{positivicur}}
\end{equation}
Since we can always reorder the integers and change the sign of $\bold q$ if necessary so that $q_i\le q_j$ and $|q_2|\le q_3$,  the  condition (\ref{positivicur}) is equivalent to $q_1+q_2>0$. This provides infinitely many Bazaikin spaces $\mathcal{B}_\bold{q}$ with metrics of positive sectional curvature, e.g.,  ${\bold{q}}=(-303-4n,-3,5,77,83 , 141+4n)$	where $n=11\times 37\times 410 \times m$ for any $m\in \mb{Z}$.

In \cite{FL} Florit and Ziller investigated the topological invariants of the Bazaikin spaces. Among others,  the cohomology rings, the linking forms, and the Pontryagin classes were calculated. Based on these basic calculations, the following finiteness assertions hold:
\begin{enumerate}
\item For a given cohomology ring, there are only finitely many positively curved Bazaikin spaces.
\item For a given homeomorphism type, there are only finitely many Bazaikin spaces.
\end{enumerate}

On the other hand, Kerin  \cite{KerM} presented infinitely many pairwise non-homeomorphic Bazaikin spaces admitting quasi-positive curvature sharing the same cohomology ring.

The main purpose of this paper is to address the topological classification of the Bazaikin spaces $\mathcal{B}_{\bold{q}}$ where $\bold{q} = (q_0, \ldots, q_5)$. This can be regarded as the $13$-dimensional counterpart of Kreck-Stolz's work \cite{KS1,KS}. But it is much more difficult in the calculations and characterizations of certain bordism groups. On the other hand, in contrast, there is no additional surgery obstruction in dimension $13$
which is different from the $7$-dimensional case where the so-called $s$-invariant occurs.

Throughout the rest of the paper let $H^i(-)$ (or $H_i(-)$) denote the integral cohomology (or homology) group $H^i(-;\mathbb{Z})$ (or $H_i(-;\mathbb{Z})$) unless other statements are given.

To state our main result, let us recall some basic facts about the topology of the Bazaikin space.
By \cite{Ba}, $\mathcal{B}_\bold{q}$ is simply connected and the non-vanishing cohomology groups are given by
$$H^0(\mathcal{B}_\bold{q})\cong   H^2(\mathcal{B}_\bold{q}) \cong   H^4(\mathcal{B}_\bold{q}) \cong   H^9(\mathcal{B}_\bold{q}) \cong   H^{11}(\mathcal{B}_\bold{q}) \cong   H^{13}(\mathcal{B}_\bold{q}) \cong  \mathbb{ Z},$$
$$\text{and} \quad H^6(\mathcal{B}_\bold{q}) \cong   H^8(\mathcal{B}_\bold{q}) \cong  \mathbb{ Z}_s$$
where the digit $s=\frac 18 |\sigma _3(\bold{q}) |\in \Bbb Z^+$, the $3$rd elementary symmetric polynomial of $\bold{q}$ which is congruent to $\pm 1\pmod 6$ since (1.2) holds. Moreover, if $x\in H^2(\mathcal{B}_\bold{q})$ is a generator, then $x^i$ is a generator of $H^{2i}(\mathcal{B}_{\bold q})$ for $i = 2$, $3$, $4$. This completely determines the cohomology ring structure of $\mathcal{B}_\bold{q}$. By \cite{FL}, the first Pontryagin class $p_1(\mathcal{B}_\bold{q}) = kx^2$ where $k=7$ or $15$ mod $24$. Note that replacing $\bold q$
with $-\bold q$, one obtains the same Bazaikin space. Therefore we can assume that $\bold q$ satisfies $\sigma _3(\bold{q})>0$ in the rest of the paper.

Now, we are ready to present the main results of this paper.

\begin{theorem}\label{classify}
Assume $s\ne 0\mod 5$, then the
Bazaikin spaces $\mathcal{B}_{\bold q}$ and $\mathcal{B}_{\bold q^\prime}$  are homeomorphic if and only if $$\sigma_2({\bold q})=\sigma_2({\bold q}^\prime),\quad  \sigma_3({\bold q})=\sigma_3({\bold q}^\prime),$$  $$ 7\sigma_4({\bold q}) =7\sigma_4({\bold q}^\prime)\pmod s,\quad \sigma_5({\bold q})=\sigma_5({\bold q}^\prime)\pmod s. $$
\end{theorem}

Indeed, $\sigma _2({\bold q})$, $\sigma _4({\bold q})$ and $\sigma _5({\bold q})$ reflect the Pontryagin classes $p_1(\mathcal{B}_\bold q)$, $ p_2(\mathcal{B}_\bold q)$, and the linking form of  $\mathcal{B}_\bold{q}$ \cite{FL}.
All of the $\sigma _i({\bold q})$ are diffeomorphism invariants of $ \mathcal{B}_{\bold q}$, but $\sigma _4(\bold q)$ may not be a homeomorphism invariant. However, $7\sigma _4({\bold q})\mod s$ is a homeomorphism invariant \cite{Sharma}. This raises an interesting question: if $s=0\pmod 7$ but $s \ne 0\pmod 5$, and the conditions in Theorem 1.1 for $\bold q$ are satisfied, how large is the set  $\{ \sigma _4(\bold q ) \pmod s\}$?



In order to obtain a diffeomorphism classification we have to determine the normal bordism group and furthermore to characterize the bordism class of the Bazaikin space  $\mathcal{B}_{\bold q}$ in the smooth category in terms of the $6$-tuple ${\bold q}$. At present, we have not obtained a complete set of numerical invariants characterizing the diffeomorphism types of Bazaikin spaces. Instead, we establish the following upper bound.
\begin{theorem}\label{diff}
Assume $s\ne 0\mod 5$.	There are  at most $12\cdot 7^n$ diffeomorphism structures on the underlying topological manifold of a Bazaikin space $\mathcal B_\q$ where $n=1$  or $0$, depending on if  $s=0 \pmod 7$ or not.
\end{theorem}

The organization of this paper is as follows. In Section \ref{Normalbg}, we establish the normal $5$-type of the Bazaikin space in the smooth and the PL category and present the necessary results of this paper. In Section \ref{sec:spectral-prelim}, we give some conventions for spectral sequence. In Section \ref{approach}, we outline our main idea and provide a reading guide.
 In Sections \ref{cohomologybarB}, \ref{pprimary}, \ref{2primary}, \ref{koMketa}, \ref{3primary}, we mainly focus on the computation of a certain bordism group. The bordism invariants of certain direct summands of the bordism group are given in Section \ref{Zs}. In Section \ref{oriented}, we present some fundamental homotopy properties of the Bazaikin spaces. In Sections  \ref{sec:Y6-restriction}, \ref{Detect}, we deal with the bordism invariants of the quotient group of the bordism group. Section \ref{PL} discusses the bordism theory in the topological category. In Sections \ref{mainproof}, we  prove the main theorems announced above.

The first four sections should be read sequentially, while the remaining sections may be read following the proof outline given in Section \ref{approach}


\section{The normal $5$-type}\label{Normalbg}

We first recall the normal-type language \cite{K} used throughout the paper.
Let
$\pi\colon \mathscr B\longrightarrow \mm{BO}$
be a fibration, and let
$\nu\colon M\longrightarrow \mm{BO}$
denote the stable normal Gauss map of a closed smooth manifold \(M\).
A normal \(\mathscr B\)-structure on \(M\) is a homotopy lift
\[
\bar\nu\colon M\longrightarrow \mathscr B
\]
of \(\nu\) along the fibration $\pi$.
Such a lift is called a normal
\(n\)-smoothing if
it is an \((n+1)\)-equivalence: it induces an isomorphism on
\(\pi_i\) for \(i\leq n\) and an epimorphism on
\(\pi_{n+1}\).

The fiber-homotopy type of a fibration
\(\pi\colon\mathscr B\to \mm{BO}\) is the normal
\(n\)-type of \(M\) if it admits a normal \(n\)-smoothing of $M$ and the homotopy groups of its homotopy fiber vanish in dimensions $\ge n+1$. The normal
\(n\)-type of a manifold is unique up to fiber homotopy equivalence.

In this section, we construct the normal \(5\)-type of a Bazaikin space.

Let \(\mathcal B_{\bold q}\) be a Bazaikin space with
$
H^6(\mathcal B_{\bold q})\cong \mb Z_s.$
Let $x\in H^2(\mathcal B_{\bold q})$ be a generator,
and let
\[
c\colon
\mathcal B_{\bold q}\longrightarrow
\mb{CP}^{\infty}\simeq\mm K(\mb Z,2)
\]
be a classifying map satisfying
$
c^*(x)=x,
$
where, by abuse of notation, $x$ also denotes the
standard generator of $H^2(\mb{CP}^{\infty})$ such that
$$H^\ast(\mb{CP}^{\infty})=\Z[x].$$ 

Consider the map
\[
r_s\colon
\mb{CP}^{\infty}\longrightarrow \mm{K}(\mb Z,6)
\]
representing the class
$sx^3\in H^6(\mb{CP}^{\infty}).$

Define \(B\) to be the homotopy fiber of \(r_s\).  Equivalently,
there is a homotopy pullback square
\begin{equation}\begin{tikzcd}
 &    &\mm{K}(\mathbb{Z},5) \ar[rr] \ar[d,""]&  & \mm{K}(\mathbb{Z},5)  \ar[d]\\
&     & B\ar[rr] \ar[d,"p"]& &\mm{PK}(\mathbb{Z},6)\ar[d]\\
\mathcal B_{\bold q}\ar[rr,"c"]\ar[rru,dotted,"f"]&  &{\mathbb{CP}}^\infty \ar[rr,"r_s"] &  & \mm{K}(\mathbb{Z},6)
\end{tikzcd}
\label{fib}
\end{equation}
where \(\mm{PK}(\mb Z,6)\) is contractible.  Hence
\begin{equation}\label{eq:Bs-fibration}
\mm{K}(\mb Z,5)
\longrightarrow
B
\stackrel{p}{\longrightarrow}
\mb{CP}^{\infty}
\end{equation}
is a homotopy fibration.

Let $l_6$ denote the fundamental class of $\mm{K}(\mb Z,6)$.
Since
\[
(r_s\circ c)^{*}(l_6)
=
sx^3
=
0
\in H^6(\mathcal B_{\bold q})
\cong\mb Z_s,
\]
the composite
$r_s\circ c\colon
\mathcal B_{\bold q}\longrightarrow \mm{K}(\mb Z,6)$
is null-homotopic.  Therefore, \(c\) admits a lift
\begin{equation}\label{eq:fq}
f\colon
\mathcal B_{\bold q}\longrightarrow B,
\qquad
p\circ f\simeq c.
\end{equation}
Moreover, since $H^5(\mathcal B_{\bold q})=0$, such a lift is unique up to homotopy.

 The space \(B\) is precisely the fifth Postnikov stage relevant
to the Bazaikin spaces. Indeed, the Postnikov \(k\)-invariant in degree \(6\) is represented by
\(sx^3\).
The fibration
\eqref{eq:Bs-fibration} gives
\[
\pi_i(B)\cong
\begin{cases}
\mb Z,&i=2,5,\\
0,&1\leq i\leq6,\ i\neq2,5.
\end{cases}
\]
On the other hand, the low-dimensional homotopy groups of $\mathcal B_{\bold q}$ are
\[
\pi_1(\mathcal B_{\bold q})
=
\pi_3(\mathcal B_{\bold q})
=
\pi_4(\mathcal B_{\bold q})
=0,
\]
\[
\pi_2(\mathcal B_{\bold q})
\cong
\pi_5(\mathcal B_{\bold q})
\cong\mb Z,
\qquad
\pi_6(\mathcal B_{\bold q})\cong\mb Z_2;
\]
see \cite{FL}.
  Consequently, 
$f\colon
\mathcal B_{\bold q}\longrightarrow B$
is a \(6\)-equivalence. 

We next incorporate the stable normal bundle.
Write
\begin{equation}\label{eq:p1Bq}
p_1(\mathcal B_{\bold q})=kx^2
\end{equation}
for the first Pontryagin class.
From \cite{FL} one has 
\begin{equation}\label{Pontry_mod8}
	k\equiv7\pmod8.
\end{equation}

Let \(L\) denote the Hopf line bundle over
\(\mb{CP}^{\infty}\), and consider the stable virtual bundle
\[
\xi=-kL
\]
over \(\mb{CP}^{\infty}\).
Its total Chern class is
$c(\xi)=(1+x)^{-k}.$
Hence
\begin{equation}\label{eq:xi-char}
p_1(\xi)
=
c_1(\xi)^2-2c_2(\xi)
=
-kx^2.
\end{equation}
Moreover, since \(k\equiv7\pmod8\),
\[
w_2(\xi)=k\bar x=\bar x,\quad w_4(\xi)
=
\frac{k(k+1)}{2}\bar x^2
=
0
\]
where $\bar x$ denotes the mod $2$ reduction of $x$.

Pull \(\xi\) back to \(B\):
\begin{equation}\label{eq:eta-def}
\eta=p^{*}\xi.
\end{equation}
Thus
\begin{equation}\label{eq:eta-char}
p_1(\eta)
=
-kp^{*}(x^2),
\qquad
w_2(\eta)
=
p^{*}\bar x,
\qquad
w_4(\eta)=0.
\end{equation}

Now let
\[
\nu_{{\bold q}}\colon
\mathcal B_{\bold q}\longrightarrow \mm{BO}
\]
be the stable normal Gauss map and denote the corresponding stable
normal bundle by
\(\mathscr N_{{\bold q}}\).
Since \(\mathcal B_{\bold q}\) is oriented,
\[
w_2(\mathscr N_{{\bold q}})
=
w_2(T\mathcal B_{\bold q})
=
k\bar x
=
\bar x\in H^2(\mathcal B_{\bold q};\Z_2).
\]
Furthermore,
\[
p_1(\mathscr N_{{\bold q}})
=
-p_1(T\mathcal B_{\bold q})
=
-kx^2\in H^4(\mathcal B_{\bold q}).
\]
Since
$p\circ f\simeq c,$
equation \eqref{eq:eta-char} gives
\[
p_1(f^{*}\eta)
=
-kx^2,
\qquad
w_2(f^{*}\eta)
=
\bar x.
\]
Therefore, the stable difference bundle
\begin{equation}\label{eq:difference-bundle}
\delta_{\bold q}
=
\mathscr N_{\bold q}
-
f^{*}\eta
\end{equation}
over $\mathcal B_{\bold q}$ satisfies
\[
w_2(\delta_{\bold q})=0,
\qquad
p_1(\delta_{\bold q})=0.
\]
Since
$H^4(\mathcal B_{\bold q})\cong\mb Z$
is torsion-free, the spin characteristic class
$
\frac{p_1}{2}(\delta_{\bold q})$
also vanishes.  Hence the classifying map of
\(\delta_{\bold q}\) lifts through the \(7\)-connected cover
\[
\mm{BO}\langle8\rangle
=
\mm{BString}
\longrightarrow
\mm{BO}.
\]
Choose such a lift
\begin{equation}\label{eq:stringlift}
g\colon
\mathcal B_{\bold q}
\longrightarrow
\mm{BO}\langle8\rangle
\end{equation}
so that 
$
\delta_{\bold q}$ is stably isomorphic to
$g^{*}\gamma_8$
where \(\gamma_8\) denotes the universal stable bundle over
\(\mm{BO}\langle8\rangle\). Moreover, since $H^1(\mathcal B_{\bold q};\mb{Z}_2)=0$ and $H^3(\mathcal B_{\bold q})=0$, such a lift $g$ is unique up to homotopy.

We now define
$
\mathscr B_{s,k}
=
\mm{BO}\langle8\rangle\times B.$
Let
\[
\pi_{s,k}\colon
\mathscr B_{s,k}\longrightarrow \mm{BO}
\]
be the classifying map of the stable  bundle
$\gamma_8\times\eta$
over \(\mathscr B_{s,k}\). By checking the homotopy groups of $B$, $\mm{BO}\langle 8\rangle$ and $\mm{BO}$, we have that the homotopy groups of the homotopy fiber of $\pi_{s,k}$ vanish in dimensions $\ge 6$.  
Define
\begin{equation}\label{eq:normal-smoothing}
\bar\nu_{\bold q}
=
(g,f)
\colon
\mathcal B_{\bold q}
\longrightarrow
\mathscr B_{s,k}.
\end{equation}
By 
$\mathscr N_{\bold q}
\cong_s
g^{*}\gamma_8
\oplus
f^{*}\eta$, we have
$\pi_{s,k}\circ\bar\nu_{\bold q}
\simeq
\nu_{\bold q}.$
Thus \(\bar\nu_{\bold q}\) is a normal \(\mathscr B_{s,k}\)-structure on \(\mathcal B_{\bold q}\). 

Since \(\mm{BO}\langle8\rangle\) is \(7\)-connected, the homotopy homomorphisms in dimensions $\le 6$ induced by
\(\bar\nu_{\bold q}\) are precisely the homomorphisms induced by
\(f\). Hence, \(\bar\nu_{\bold q}\) is a normal $5$-smoothing of $\mathcal B_{\bold q}$, and thus the fiber homotopy type of $\pi_{s,k}$ is the normal $5$-type of $\mathcal B_{\bold q}$,  

The normal $5$-type of $\mathcal B_{\bold q}$ gives a bordism group $\Omega_{13}^{\mm{O}\langle 8\rangle}(\eta)$, and the normal $5$-smoothing \eqref{eq:normal-smoothing} represents a class $[\mathcal B_{\bold q},(g,f)]$ in $\Omega_{13}^{\mm{O}\langle 8\rangle}(\eta)$.   

\begin{proposition}\label{prop:normal5type}
Any two Bazaikin spaces $\mathcal B_\q $ and $\mathcal B_{\q'} $ with the same \(s\) and the
same first Pontryagin coefficient \(k\) admit normal
\(5\)-smoothings \eqref{eq:normal-smoothing}
 into the common normal \(5\)-type
\[
\pi_{s,k}:\mathscr B_{s,k}
=
\mm{BO}\langle8\rangle\times B\to \mm{BO}.
\]

Furthermore, if the normal $5$-smoothings are bordant, then $\mathcal B_\q $ is diffeomorphic to  $\mathcal B_{\q'}$.
\end{proposition}
\begin{proof}
	The construction of the normal $5$-type depends only on the two integers
$s=|H^6(\mathcal B_{\bold q})|$
and $k$
where
$p_1(\mathcal B_{\bold q})=kx^2.$
If one take $-x$ as the generator of $H^2(\mathcal B_{\bold q})$, the coefficient $k$ is unchanged since $(-x)^2=x^2$. Hence, the normal \(5\)-type
\(\mathscr B_{s,k}\) is independent of this choice.

The diffeomorphism statement relies on  \cite{K}. Readers may consult Appendix \ref{Surgeryobstruction} for more details.
\end{proof} 

\begin{remark}
Once the classifying map
\[
c\colon
\mathcal B_{\bold q}\longrightarrow\mb{CP}^{\infty}
\]
is fixed, the lift \(f\) in \eqref{eq:fq} is unique up to
homotopy.
The String lift in \eqref{eq:stringlift} is likewise unique up to homotopy.
  Thus the normal
\(5\)-smoothing \(\bar\nu_{\bold q}\) is canonical up to homotopy
after the class \(c\) has been fixed. Indeed, there is a canonical classifying map
$
c\colon
\mathcal B_{\bold q}\longrightarrow\mb{CP}^{\infty}$, which classifies the principle $S^1$ bundle
\begin{equation}
	S^1\to \SU(6)/\Sp(3)\to \mathcal B_\q =S^1\backslash \SU(6)/\Sp(3) \label{principleS1}
\end{equation}
coming from the construction \eqref{Bazaibiquotient}. After this, we always choose such classifying map and get the canonical normal $5$-smoothing of $\mathcal B_\q $.
\end{remark}

Let $\mm{BPL}$ be the classifying space for stable piecewise linear (PL) bundles. The natural map $\mm{BO}\to \mm{BPL}$ induces a map $T:\mm{BO}\langle 8\rangle\to \mm{BPL}\langle 8\rangle$ between their $7$-connected covers. Let 
$$\pi_{s,k}^{pl}:\mathscr B_{s,k}
^{pl}=\mm{BPL}\langle 8\rangle\times B\to \mm{BPL}$$ be the classifying map of the stable PL bundle $\gamma_8^{pl}\times \eta$ where $\gamma_8^{pl}$ is the universal PL bundle over $\mm{BPL}\langle 8\rangle$.

The fiber homotopy type of $\pi_{s,k}^{pl}:\mathscr B^{pl}_{s,k}\to \mm{BPL}$ is the PL normal $5$-type of $\mathcal B_\q $, and admits the PL normal $5$-smoothing of $\mathcal B_\q$: 
\begin{equation}\label{eq:PLnormal-smoothing}
	(T\circ g,f):\mathcal B_\q\stackrel{\Delta}{\longrightarrow }\mathcal B_\q\times \mathcal B_\q\stackrel{g\times f}{\longrightarrow}\mm{BO}\langle 8 \rangle \times B\stackrel{T\times \mm{id}}{\longrightarrow} \mm{BPL}\langle 8 \rangle \times B.
\end{equation} 

Let $\Omega_{13}^{\mm{PL}\langle 8 \rangle}(\eta)$ denote a bordism group given by the PL normal $5$-type. The PL normal $5$-smoothing of $\mathcal B_\q$ determines an element in $\Omega_{13}^{\mm{PL}\langle 8 \rangle}(\eta)$.

\begin{proposition}\label{PLhcobor}
Any two Bazaikin spaces $\mathcal B_\q $ and $\mathcal B_{\q'} $ with the same \(s\) and the
same first Pontryagin coefficient \(k\) admit PL normal
\(5\)-smoothings \eqref{eq:PLnormal-smoothing}
 into the common PL normal \(5\)-type
\[
\pi^{pl}_{s,k}:\mathscr B_{s,k}^{pl}
=
\mm{BPL}\langle8\rangle\times B\to \mm{BPL}.
\]

Furthermore, if the PL normal $5$-smoothings are bordant, then $\mathcal B_\q $ is PL homeomorphic to  $\mathcal B_{\q'}$.
\end{proposition}
\begin{proof}
	 The first claim is trivial. The second statement uses the PL counterpart of \cite{K}. 
\end{proof}

\section{Spectral sequences and filtration conventions}
\label{sec:spectral-prelim}

In this section, we recall the spectral sequences and
filtration conventions used throughout the paper.
We also recall some standard properties concerning
naturality, detection.

\subsection{Filtrations and detection}
\quad

Let
$
E^r_{p,q}\Longrightarrow G_{p+q}$
be a convergent homological spectral sequence.
It carries an increasing bidgree filtration
\[
0=F_{-1,n+1}
\subseteq F_{0,n}
\subseteq F_{1,n-1}
\subseteq\cdots\subseteq F_{n,0}=G_n
\]
such that
\[
E^\infty_{p,n-p}
\cong F_{p,n-p}/F_{p-1,n+1-p}.
\]

\begin{definition}\label{provided}
A class $x\in E^2_{p,q}$ is called a
\emph{permanent cycle} if it survives to a
nonzero class in $E^\infty_{p,q}$.

 A permanent cycle $x$ is
the image of $z\in F_{p,\ast}$ in
$
F_{p,\ast}/F_{p-1,\ast+1}$.
And we say that $z$ is \emph{detected by} $x$,
or that $z$ is \emph{provided by}
$E^2_{p,q}$.
In this case, $z$ has filtration $p$.
\end{definition}


We shall use the following elementary observation.

\begin{proposition}\label{cosetspectralsequence}
Let $x\in E^\infty_{p,q}$ be a nonzero class.
Then the elements of $G_{p+q}$ detected by $x$
form a coset of $F_{p-1,q+1}$ in
$F_{p,q}$.
\end{proposition}

\begin{proof}
This follows immediately from $E^\infty_{p,n-p}
\cong F_{p,n-p}/F_{p-1,n+1-p}$.
\end{proof}

Let $\{x\}$ denote the coset detected by $x\in E^\infty_{p,q}$.

In particular, a nonzero class on the
$E^\infty$-page determines an element of $G$ only up to terms of lower filtration
for an increasing filtration.

Similarly, let
$
E_r^{s,t}\Longrightarrow G$
be a convergent cohomological spectral sequence.
The abutment carries a decreasing bidegree filtration.
The analogues of
Definition~\ref{provided} and
Proposition~\ref{cosetspectralsequence}
hold with respect to this filtration.

\subsection{Naturality of spectral sequences}\quad

For two convergent homological spectral sequences
\[
E^r_{p,q}\Longrightarrow G_{p+q},
\qquad
\bar E^r_{p,q}\Longrightarrow \bar G_{p+q},
\]
a morphism between them consists of
homomorphisms
\[
h_r:E^r_{p,q}\longrightarrow \bar E^r_{p,q},
\qquad r\ge2,
\]
which commute with the differentials,
\[
h_r\circ d_r=\bar d_r\circ h_r,
\]
and are compatible with the passage from
$E_r$ to $E_{r+1}$.

We shall frequently use naturality to compare
differentials and detect permanent cycles.
For instance, the Atiyah--Hirzebruch spectral sequence
\[
E^2_{p,q}(E,X)
=
H_p(X;\pi_q(E))
\Longrightarrow E_{p+q}(X)
\]
is natural in both the coefficient spectrum $E$
and the spectrum $X$.
In particular, a map of coefficient spectra
\[
c:E\longrightarrow E'
\]
induces a morphism of spectral sequences
\[
h_r:E^r_{p,q}(E,X)
\longrightarrow E^r_{p,q}(E',X),
\]
whose map on the $E_2$-page is induced by
\[
c_*:\pi_q(E)\longrightarrow\pi_q(E').
\]
Finally, this morphism induces a homomorphism on the $E^\infty$-pages
$$h_\infty:E^\infty_{p,q}(E,X)
\longrightarrow E^\infty_{p,q}(E',X).$$

On the other hand, the coefficient map $c$ induces the
homomorphism
\[
c_*:E_n(X)\longrightarrow E'_n(X),
\]
which preserves the Atiyah--Hirzebruch (AH) filtrations, i.e.,
\[
c_*(F_{p,q})\subseteq F'_{p,q}.
\]
Moreover, the induced homomorphism on the $E^\infty$-pages agrees with
$h_\infty$ under the identifications
\[
E^\infty_{p,q}
\cong F_{p,q}/F_{p-1,q+1},\quad
E'^\infty_{p,q}
\cong F'_{p,q}/F'_{p-1,q+1}.
\]

For the map $t:X\to X'$, we also have the analogues as above.

 Suppose that
$z\in E_{p+q}(X)$ is detected by  $x\in E^2_{p,q}(E,X)$,
and let $[x]\in E^\infty_{p,q}(E,X)$
denote its nonzero surviving class.
If
$
h_\infty([x])\ne0,$
then $c_*(z)\ne0$ and is detected by
$h_2(x)$ in the target spectral sequence.

It is particularly important to note that if
$
h_\infty([x])=0,$
then
\[
c_*(z)\in F'_{p-1,q+1}.
\]
In particular, the vanishing of
$h_\infty([x])$ does not necessarily imply
that $c_*(z)=0$; its image may remain
nonzero in a lower filtration.

We shall also use the following comparison
principle.

\begin{lemma}\label{ahsspro}
Let
$
h_r:E^r_{p,q}\longrightarrow\bar E^r_{p,q}$
be a morphism of the homological AHSSs, with
$
E^2_{p,q}=\bar E^2_{p,q}=0$ for $p<0$ or $q<0$.  
Suppose that, for $n>0$,
\[
h_2:E^2_{p,q}\xrightarrow{\cong}\bar E^2_{p,q},
\qquad 0\le p<n.
\]
Then:

\begin{enumerate}
\item[(1)]
For every $r\ge2$, the map $h_r$ is an
isomorphism for $0\le p\le n-r$ and an
epimorphism for $n-r<p<n$.

\item[(2)]
If $x\in E^2_{n,m}$ satisfies
$h_2(x)=\bar x\ne0$, and $\bar x$ survives
to $\bar E^{n+1}_{n,m}$, then $x$ survives
to $E^{n+1}_{n,m}$.
\end{enumerate}
\end{lemma}

\begin{proof}
The first assertion follows inductively from
the naturality of the differentials and the
definition of $E^r$-pages.

For the second assertion, suppose that $x$
has survived to $E^r_{n,m}$ where $2\le r<n+1$.  Since $\bar x\in \bar E^r_{n,m}$ is a cycle, i.e., $\bar d_r(\bar x)=0$, we have
\[
h_r(d_r x)=\bar d_r(h_r x)=0.
\]
The target of $d_rx$ lies in filtration $n-r$, where $h_r$ is
injective by (1), and hence $d_rx=0$.  Moreover, if
$x=d_ry$,
then
$
\bar x=\bar d_r(h_r y),
$
a contradiction.  Thus $x$ survives to
$E^{r+1}$. We complete the proof by induction. 
\end{proof}

\section{Outline of the Proof}\label{approach}

Let $\mathcal B_\q$ and $\mathcal B_{\q'}$ be two Bazaikin spaces satisfying
\[
\sigma_i(\q)=\sigma_i(\q'),\qquad i=2,3.
\]
Our aim is to determine when they are PL homeomorphic or diffeomorphic.
By Propositions \ref{prop:normal5type} and \ref{PLhcobor}, this problem
reduces to comparing their canonical normal $5$-smoothings in the
corresponding normal bordism groups. More precisely, we need to determine
when
\[
[\mathcal B_\q,(g,f)]
=
[\mathcal B_{\q'},(g',f')]
\in \bord_{13}(\eta)
\]
in the smooth category, or
\[
[\mathcal B_\q,(g,T\circ f)]
=
[\mathcal B_{\q'},(g',T\circ f')]
\in \Omega_{13}^{\mm{PL}\langle 8\rangle}(\eta)
\]
in the PL category.

The proof proceeds in two main stages. First, we establish equality in
a suitable quotient of the smooth normal bordism group by using
connective real $K$-theory. The key step is to construct a common
finite-dimensional space on which the two Bazaikin spaces induce
isomorphic bundles, thereby relating the corresponding $ko$-classes.
We then use the Atiyah--Hirzebruch filtration to identify the remaining
obstruction and show that it vanishes after passing to the PL category.

\medskip
\noindent
\textbf{Step 1. Reduction to a quotient of $\bord_{13}(\eta)$.}

Let $\MO\langle 8\rangle$ and $\M\eta$ denote the Thom spectra of
$\gamma_8$ and $\eta$, respectively. By the Thom--Pontryagin
isomorphism \cite{Sw}, we have
\[
\Omega_{13}^{\mm{O}\langle 8\rangle}(\eta)
\cong
\pi_{13}(\mm{MO}\langle 8\rangle\wedge\mm{M}\eta).
\]
Under this isomorphism, the normal $5$-smoothing
$[\mathcal B_{\bold q},(g,f)]$ is represented by the stable homotopy
class
\begin{equation}\label{Bqhomoclass}
\{\mathcal B_{\bold q}\}:
S^{13}\stackrel{\ell}{\longrightarrow}
\mm{M}_{\mathcal B_{\bold q}}
\longrightarrow
\mm{MO}\langle 8\rangle\wedge\mm{M}\eta,
\end{equation}
where $\mm{M}_{\mathcal B_{\bold q}}$ denotes the Thom spectrum of
$\mathscr N_{\mathcal B_{\bold q}}$. Here, $\ell$ is the
Thom--Pontryagin map, and the second map is induced by
\[
(g,f):\mathcal B_{\bold q}
\xrightarrow{g\times f}
\mm{BO}\langle 8\rangle\times B.
\]
To make the dependence on $f$ more explicit, we factor this map as
\[
\mathcal B_{\bold q}
\xrightarrow{g\times 1}
\mm{BO}\langle 8\rangle\times\mathcal B_{\bold q}
\xrightarrow{1\times f}
\mm{BO}\langle 8\rangle\times B.
\]
Consequently, the class in \eqref{Bqhomoclass} admits a factorization
\[
\{\mathcal B_{\bold q}\}:
S^{13}\xrightarrow{\varsigma_\q}
\mm{MO}\langle 8\rangle\wedge\mm{M}\varphi
\xrightarrow{1\wedge f}
\mm{MO}\langle 8\rangle\wedge\mm{M}\eta,
\]
where $\M\varphi$ is the Thom spectrum of the pullback bundle
$
\varphi:=f^\ast\eta.$
We use the same notation for the map induced between the corresponding
Thom spectra. Thus, the class $\varsigma_\q$ records the normal
bordism data before applying the map induced by $f$.

Recall that the Atiyah--Hirzebruch (AH) filtration on
$\Omega_{13}^{\mm{O}\langle 8\rangle}(\eta)$ is given by
\[
\Omega_{13}^{\mm{O}\langle 8\rangle}(\eta)
=F_{13,0}\supset F_{12,1}\supset\cdots
\supset F_{1,12}\supset F_{0,13}.
\]
We first compare the two bordism classes in the quotient
\[
\Omega_{13}^{\mm{O}\langle 8\rangle}(\eta)/F_{8,5}.
\]

By the computations in Section \ref{pprimary}, the quotient
$\bord_{13}(\eta)/F_{8,5}$ is finite and admits the decomposition
\begin{equation}\label{decompo_bordism_Section4}
\bord_{13}(\eta)/F_{8,5}
\cong
\Z_s\oplus\Z_5\oplus A_{(2)}\oplus A_{(3)},
\end{equation}
where $A_{(p)}$ denotes the $p$-localization of $A$.
By Lemma \ref{130Zs}, the component of
$\{\mathcal B_\q\}$ in the summand $\Z_s$ is determined by its
linking form $\ell k(\mathcal B_\q)$, which in turn is determined by
$\sigma_5(\mathcal B_\q)\pmod s$.
Thus, once the congruence condition on $\sigma_5$ is imposed,
the remaining problem is to compare the components in the other
summands of \eqref{decompo_bordism_Section4}.

To carry out this comparison, we use connective real $K$-theory.
Let $ko$ denote the connective spectrum of real $K$-theory, and
consider the composition
\[
\alpha_8:
\mm{MO}\langle 8\rangle
\longrightarrow
\mm{MSpin}
\stackrel{\alpha}{\longrightarrow}ko,
\]
where $\alpha$ is the Atiyah--Bott--Shapiro orientation
\cite{AtiyahBottShapiro}. The first map is induced by the map
between the Thom spaces of the corresponding universal bundles.
It induces a homomorphism
\[
\alpha_{8\ast}:
\Omega_{13}^{\mm{O}\langle 8\rangle}(\eta)
\longrightarrow
ko_{13}(\mm{M}\eta)
\cong
\pi_{13}(ko\wedge\mm{M}\eta).
\]

The following result, proved in Section \ref{koMeta85_quot},
shows that $ko$-theory detects all the information needed in the
quotient under consideration.

\begin{theorem}\label{pFZ4_Sec4}
The homomorphism induced by $\alpha_8$,
\[
(\alpha_8)_*:
{\Omega_{13}^{\mm O\langle8\rangle}(\eta)}/{F_{8,5}}
\longrightarrow
{ko_{13}(\mm M\eta)}/{\bar F_{8,5}},
\]
is injective, where $\bar F_{8,5}$ is the AH filtration of
$ko_{13}(\mm M\eta)$.
\end{theorem}

In view of this injectivity, it suffices to prove that the two
corresponding $ko$-classes agree modulo $\bar F_{8,5}$.

\medskip
\noindent
\textbf{Step 2. Relating the two $ko$-classes through a common complex.}

For simplicity, let
\[
\rho_\q=\alpha_{8\ast}(\varsigma_\q)
\in \pi_{13}(ko\wedge\mm{M}\varphi).
\]
By naturality, we obtain the commutative diagram
\[
\xymatrix@C=.7cm{
\varsigma_\q\in
\pi_{13}(\mm{MO}\langle 8\rangle\wedge\mm{M}\varphi)
\ar[d]^-{\alpha_{8\ast}}
\ar[r]^-{f_\ast}
&
\pi_{13}(\mm{MO}\langle 8\rangle\wedge\mm{M}\eta)
\ar[d]^-{\alpha_{8\ast}}
\\
\rho_\q\in
\pi_{13}(ko\wedge\mm{M}\varphi)
\ar[r]^-{f_\ast}
&
\pi_{13}(ko\wedge\mm{M}\eta).
}
\]
It follows from Theorem \ref{pFZ4_Sec4} that if
\[
f_\ast(\rho_\q)=f'_\ast(\rho_{\q'})
\in
ko_{13}(\mm M\eta)/\bar F_{8,5},
\]
then
\[
\{\mathcal B_\q\}
=f_\ast(\varsigma_\q)
=f'_\ast(\varsigma_{\q'})
=\{\mathcal B_{\q'}\}
\in
\bord_{13}(\eta)/F_{8,5}.
\]
Thus, the main difficulty is to compare the maps $f$ and $f'$,
which have a common target $\M\eta$ but different source spaces.

The key idea is to construct a common finite complex over which the
two Bazaikin spaces determine isomorphic $P$-bundles. Let
$
P=\SU(6)/\Sp(3).$
We consider the Borel construction
\begin{equation}\label{universal-Pq_Section4}
P
\longrightarrow
S^\infty\times_{S^1}P
\longrightarrow
\CP^\infty,
\end{equation}
where the $S^1$-action on $S^\infty$ is the right Hopf action,
and the action on $P$ is the right free action arising from the
principal bundle \eqref{principleS1}. The total space of \eqref{universal-Pq_Section4} is homotopy
equivalent to $\mathcal B_\q$.

The idea of comparing Bazaikin spaces by restricting their associated bundles \eqref{universal-Pq_Section4} to a suitable finite-dimensional complex originates in the Shen's doctoral dissertation. In an earlier arXiv preprint, we pursued this approach by restricting the bundles \eqref{universal-Pq_Section4} to $\CP^3$. However, we did not obtain a complete topological classification through this restriction. In the present paper, we refine this approach by choosing a more carefully constructed finite complex (Y), which allows us to establish a stronger relation between the two restricted $P$-bundles and ultimately obtain the desired classification.

We now pull back \eqref{universal-Pq_Section4} to a fixed finite
complex $Y$, chosen to be the $10$-skeleton of the sphere bundle of
the complex vector bundle
\[
L\oplus L\oplus L\oplus L^{\otimes 6}
\]
over $\CP^5$ where $L$ is the Hopf line bundle over $\CP^5$. This gives the commutative diagram
\begin{equation}\label{eq:H9-circle-square_Section4}
\begin{tikzcd}[column sep=large,row sep=large]
P\arrow[r,equal]\arrow[d]
  &P\arrow[r," "]\arrow[d,"\pi"']
  &\K(\Z,5)\arrow[d,"i"]\\
X_\q\arrow[r,"j"]\arrow[d,"\pi_Y"']
  &S^\infty\times_{S^1} P\simeq\mathcal B_\q
       \arrow[r,"f"]\arrow[d]
  &B\arrow[d,"p"]\\
Y\arrow[r,"c_Y"]
  &\CP^\infty\arrow[r,equal]
  &\CP^\infty.
\end{tikzcd}
\end{equation}

The importance of this construction is that the base $Y$ is independent
of $\q$. In Proposition \ref{prop:Y-restriction-isomorphic}, we prove
that
\[
X_\q\cong X_\q'
\]
as $P$-bundles. Consequently, the two pullback constructions give
the same map to $\CP^\infty$ up to homotopy:
\[
p\circ f\circ j\simeq p\circ f'\circ j'.
\]
Since $H^5(X_\q)=0$ (cf. Proposition \ref{prop:restricted-P-ring})
and $p$ is a fibration with fiber $\K(\Z,5)$, the homotopy above
has a lift. In particular,
\begin{equation}
	f\circ j\simeq f'\circ j' \label{basehomotopy_Section4}.
\end{equation}

Let $\omega$ denote the pullback bundle
\[
\omega=j^\ast\varphi=(f\circ j)^\ast\eta,
\]
and let $\M\omega$ be its Thom spectrum. Lifting the homotopy
\eqref{basehomotopy_Section4} to the corresponding Thom spectra,
we obtain the following homotopy commutative diagram:
\begin{equation}\label{commuta_Thomspec_Section4}
\begin{tikzcd}[column sep=large,row sep=large]
\mm M\omega
\arrow[r,"j"]
\arrow[d,"j'"']
&
\mm M\varphi
\arrow[d,"f"]
\\
\mm M\varphi'
\arrow[r,"f'"']
&
\mm M\eta .
\end{tikzcd}
\end{equation}
This diagram provides the common framework needed to compare
$\rho_\q$ and $\rho_{\q'}$: both classes can now be studied through
the same group $ko_{13}(\M\omega)$.

\medskip
\noindent
\textbf{Step 3. Comparing the classes in the common $ko$-group.}

The next task is to construct an element
$\theta\in ko_{13}(\M\omega)$ whose images under the two
directions of \eqref{commuta_Thomspec_Section4} recover suitable
multiples of $\rho_\q$ and $\rho_{\q'}$. More precisely, we will
show that
\[
j_\ast(\theta)
=s\rho_\q+n\bar\rho_\q
\in ko_{13}(\M\varphi)/\mm{F}_{8,5}
\]
and
\[
j'_\ast(\theta)
=s\rho_{\q'}+n'\bar\rho_{\q'}
\in ko_{13}(\M\varphi')/\mm{F}_{8,5}' ,
\]
where $\mm{F}_{8,5}$ and $\mm{F}_{8,5}'$ denote the AH filtrations
for $ko_{13}(\M\varphi)$ and $ko_{13}(\M\varphi')$, respectively.

The role of $\theta$ is to place the two classes in a common
$ko$-theoretic setting. The remaining terms involving
$\bar\rho_\q$ and $\bar\rho_{\q'}$ are then handled by constructing
additional elements of $ko_{13}(\M\omega)$. Using these elements
together with the diagram \eqref{commuta_Thomspec_Section4}, we
will prove that
\[
f_\ast(n\bar\rho_\q)
=
f'_\ast(n'\bar\rho_{\q'})
\in
(ko_{13}(\M\eta)/\bar F_{8,5})_{(p)}
\]
for every prime $p$ satisfying $\gcd(p,s)=1$.
Subtracting these equalities from the images of $\theta$ gives
\[
f_\ast(s\rho_\q)
=
f'_\ast(s\rho_{\q'})
\in
(ko_{13}(\M\eta)/\bar F_{8,5})_{(p)}.
\]

Recall that $\gcd(s,6)=1$. If we further assume that
$\gcd(s,5)=1$, then $s$ is invertible at the primes
$p=2,3,5$. Therefore,
\[
f_\ast(\rho_\q)
=
f'_\ast(\rho_{\q'})
\in
(ko_{13}(\M\eta)/\bar F_{8,5})_{(p)},
\qquad p=2,3,5.
\]
Together with the decomposition \eqref{decompo_bordism_Section4},
this yields
$$
\{\mathcal B_\q\}
=
\{\mathcal B_{\q'}\}
\in
\bord_{13}(\eta)/F_{8,5}.
$$ In other words, 
$$
\{\mathcal B_\q\}
-
\{\mathcal B_{\q'}\}
\in
F_{8,5}\subset \bord_{13}(\eta).
$$

%

\medskip
\noindent
\textbf{Step 4. Passing from smooth bordism to PL bordism.}

It remains to understand the subgroup $F_{8,5}$.
In Section \ref{PL}, we show that
\[
F_{8,5}\cong\Z_s\oplus F_{4,9}.
\]
Under the forgetful map
\[
T_\ast:
\bord_{13}(\eta)
\longrightarrow
\Omega_{13}^{\mm{PL}\langle 8\rangle}(\eta),
\]
the subgroup $F_{4,9}$ vanishes. Moreover, if
\[
7\sigma_4(\q)\equiv 7\sigma_4(\q')\pmod s,
\]
then the remaining $\Z_s$-component of the difference $\{\mathcal B_\q\}
-
\{\mathcal B_{\q'}\}$ also maps
to zero. Consequently,
\[
0=T_\ast\big(
\{\mathcal B_\q\}-\{\mathcal B_{\q'}\}
\big)
\in\widetilde F_{8,5}
\subset
\Omega_{13}^{\mm{PL}\langle 8\rangle}(\eta).
\]
Thus, the canonical normal $5$-smoothings of
$\mathcal B_\q$ and $\mathcal B_{\q'}$ are bordant in the PL
normal bordism group. By Proposition \ref{PLhcobor}, it follows
that $\mathcal B_\q$ and $\mathcal B_{\q'}$ are PL homeomorphic.

\medskip
\noindent
\textbf{Guide to the readers.}

The proof is organized according to the four steps above.
Section \ref{pprimary} computes the relevant quotient of the
smooth normal bordism group and Section \ref{Zs} gives the geometric invariants of certain summands of
$\bord_{13}(\eta)$. Section \ref{koMeta85_quot}
establishes the injectivity of the $ko$-theoretic detection map,
which reduces the bordism comparison to a calculation in
connective real $K$-theory.  Section \ref{sec:Y6-restriction} constructs the
finite complex $Y$, establish the isomorphism of the restricted
$P$-bundles, and then Section \ref{Detect} uses the resulting common Thom spectrum $\M\omega$
to compare the two $ko$-classes. Finally, Section \ref{PL}
identifies the remaining filtration subgroup and shows that its
contribution disappears in PL bordism. 
The remaining sections make the necessary computations for the bordism group $\bord_{13}(\eta)$ and other homotopy groups of certain spectra.

\section{The cohomology ring of \(B\)}
\label{cohomologybarB}

Let \(B\) be the space defined by the homotopy fibration
\[
\mm{K}(\mathbb Z,5)\stackrel{i}{\longrightarrow}
B
\stackrel{p}{\longrightarrow}
\CP^\infty ,
\]
constructed from the map
$
r_s:\CP^\infty\longrightarrow \mm{K}(\mathbb Z,6)
$
representing the cohomology class
$
sx^3\in H^6(\CP^\infty ).$
In our application, by \cite{FL},
\[
s=\left|H^6(\mathcal B_{\bold q})\right|\equiv \pm1\pmod 6.
\]

In this section, we compute the cohomology ring of \(B\) in low
dimensions. We also determine the structure of
$
H^*(B;\mathbb Z_2)
$
as a module over the Steenrod sub-algebra $\mathscr{A}_2$ generated by
\(\Sq^1,\Sq^2,\Sq^4\).

\subsection{\(H^*(B;\mathbb Z_2)\) and its module over the Steenrod
algebra}
\label{mod2cohomologyB}\quad

We first compare the Serre spectral sequences associated with the
following two fibrations:
\begin{equation}\label{B_commutediagram}
\begin{tikzcd}
\K(\mathbb Z,5)\ar[rr]\ar[d,"i"]
&&\K(\mathbb Z,5)\ar[d]
\\
B\ar[rr]\ar[d,"p"]
&&\mathrm{PK}(\mathbb Z,6)\simeq *
\ar[d]
\\
\CP^\infty\ar[rr,"r_s"]
&&\K(\mathbb Z,6).
\end{tikzcd}
\end{equation}
The corresponding Serre spectral sequences with
\(\mathbb Z_2\)-coefficients are
\[
\{E_r^{p,q}(*),\bar d_r
\}
\Longrightarrow
H^{p+q}(*;\mathbb Z_2),
\]
\[
\{E_r^{p,q}(B), d_r
\}\Longrightarrow
H^{p+q}(B;\mathbb Z_2).
\]
In particular, $$E_2^{p,q}(*)=
H^p(\K(\mathbb Z,6);
H^q(\K(\mathbb Z,5);\mathbb Z_2)),$$ $$E_2^{p,q}(B)
=
H^p(\CP^\infty;
H^q(\K(\mathbb Z,5);\mathbb Z_2)).$$


Recall that
$H^*(\K(\mathbb Z,n);\mathbb Z_2)$
is generated by the fundamental class \(l_n\) and the classes
$
\Sq^I l_n,$
where \(I=(i_1,\dots ,i_r)\) is an admissible sequence satisfying
$
e(I)=2i_1-(i_1+\cdots+i_r)<n,$
$i_r>1 .
$ In particular, 
$$\Sq^1l_n=0.$$ Let $d(I)=\Sigma_{j=1}^ri_j$.

\begin{lemma}\label{KZ5SSS}
In the spectral sequence $\{E_r^{p,q}(*),\bar d_r
\}$, any generator $\Sq^Il_5$ of $H^\ast(\K(\mb{Z},5);\mb{Z}_2)$ survives to $E_{6+d(I)}^{0,5+d(I)}(\ast)$. Moreover, 
$$
E_6^{p,q}(*)=E_2^{p,q}(*).$$
\end{lemma}

\begin{proof}
First $l_5$ survives to $E_6^{0,5}(\ast)$. Then \cite[Corollary 6.9]{McCleary}  implies that $\Sq^Il_5$ survives to $E_{6+d(I)}^{0,5+d(I)}(\ast)$. Meanwhile, 
\begin{equation}\label{eq:d6inKZ5}
	\bar d_r(\Sq^Il_5)=0 \quad \text{for $2\le r\le 5$}.
\end{equation}

The cohomology ring $H^*(\K(\mathbb Z,5);\mathbb Z_2)$
is generated by  \(l_n\) and the classes
$
\Sq^I l_n.$
 By the derivation property of differentials, we have 
 $$\bar d_r(E_r^{p,q}(*))=0\quad \text{for $2\le r\le 5$}.$$
Hence $E_6^{p,q}(*)=E_2^{p,q}(*)$.
\end{proof}

The commutative diagram \eqref{B_commutediagram} induces a morphism between the
spectral sequences of the two fibrations. Since 
$\bar d_6(l_5)=l_6,$
we have
$$
d_6(l_5)=r_s^*(l_6)=s  x^3$$
where $x$ denotes the generator of $ H^2(\CP^\infty;\Z_2)$.
Since \(s\) is odd, this becomes
\[
d_6(l_5)= x^3.
\]
Moreover, by Lemma \ref{KZ5SSS}, all differentials \(d_r\) with \(2\le r<6\) vanish.
So,
\[
E_6^{p,q}(B)=E_2^{p,q}(B).
\]
\begin{lemma}\label{vanish0}
For the spectral sequence of \(B\),
$
E_7^{p,q}(B)=0$
whenever
$p\ge5 .$
\end{lemma}

\begin{proof}
The groups
$E_2^{p,q}(B)$
vanish for odd \(p\), hence in particular
\[
E_2^{5,q}(B)=0.
\]

For \(p\ge6\), every element of \(E_6^{p,q}(B)=E_2^{p,q}(B)\) can be written as
\[
x^3a,
\]
where
$
x^3\in E_6^{6,0}(B),$ $a\in E_6^{p-6,q}(B)$.
If
$
d_6(x^3a)=0,$
then, by the derivation property, we have
\[
0=d_6(x^3a)=d_6(x^3)a+x^3d_6(a)=x^3d_6(a).
\]
This implies $d_6(a)=0$. Thus,
\[
d_6(l_5a)
=
x^3a+l_5d_6(a)
=
x^3a.
\]
Hence every \(d_6\)-cycle with \(p\ge6\) is a boundary. 
\end{proof}

\begin{proposition}\label{convergence}
The spectral sequence converges at the \(E_7\)-page:
\[
E_\infty^{p,q}(B)=E_7^{p,q}(B).
\]
\end{proposition}

\begin{proof}
By Lemma~\ref{vanish0}, for \(r\ge7\), all targets of differentials
\(d_r\) vanish.
\end{proof}

\begin{lemma}
In the range up to dimension $15$, the mod $2$ cohomology ring of $B$ is as in the following table

\begin{table}[!h]
\begin{center}
\caption{$H^\ast(B;\mb{Z}_2)$}
\label{HBZ2}
\begin{tabular}{|c|c|c|c|c|c|c|c|c|c|c|c|c|c|c|c|}
\hline
$n$ & 1 & 2 & 3 &4& 5& 6&7 & 8& 9 & 10 \\
\hline
$H^n(B;\mathbb{Z}_2)$ &0& $x$&0& $x^2$& 0&0& $u$& $\mathrm{Sq}^1u$&$ux$, $v$&$\mathrm{Sq}^2\mathrm{Sq}^1u$, $x\mathrm{Sq}^1u$\\
\hline
\end{tabular}
\end{center}

\begin{center}
\begin{tabular}{|c|c|c|c|c|c|c|c|c|c|c|c|c|c|c|c|}
\hline
 11 & 12 & 13  \\
\hline
 $ux^2$, $vx$, $\mathrm{Sq}^4 u$& $x^2\mathrm{Sq}^1 u$, $x\mathrm{Sq}^2\mathrm{Sq}^1 u$, $\mathrm{Sq}^5u$&$vx^2$, $x\mathrm{Sq}^4 u$, $\mathrm{Sq}^6 u$\\
\hline
\end{tabular}
\end{center}
\begin{center}
\begin{tabular}{|c|c|c|c|c|c|c|c|c|c|c|c|c|c|c|c|}
\hline
14& 15 \\
\hline
 $u^2$, $x^2\mathrm{Sq}^2\mathrm{Sq}^1 u$, $x\mathrm{Sq}^5u$, $\mathrm{Sq}^6 \mathrm{Sq}^1 u$& $x^2\mathrm{Sq}^4 u$, $ x\mathrm{Sq}^6 u$, $u\mathrm{Sq}^1u$, $\mathrm{Sq}^7\mathrm{Sq}^1 u$\\
\hline
\end{tabular}
\end{center}
\end{table}
where $u \in H^7(B;\mathbb{Z}_2)$ and $v \in   H^9(B;\mathbb{Z}_2)$ satisfy
\begin{equation}
	i^\ast (u)=\mathrm{Sq}^2l_5\quad i^\ast (v)=\mathrm{Sq}^4l_5 \quad \mathrm{Sq}^1v=\mathrm{Sq}^2\mathrm{Sq}^1u.  \label{chooseuv}
\end{equation}
$$i^\ast:H^{\ast}(B;\mb{Z}_2)\to H^{\ast}(\mm{K}(\mb{Z},5)
;\mb{Z}_2)$$

\end{lemma}
\begin{proof}
By Proposition~\ref{convergence}, it is enough to compute the
\(E_7\)-page in total degrees at most \(15\).
By Lemma~\ref{vanish0}, it suffices to determine
\[
\ker\bigl(
d_6:E_6^{p,q}(B)\longrightarrow E_6^{p+6,q-5}(B)
\bigr)
\]
for \(p<5\), hence for \(p=0,2,4\). By the derivation property, the kernel is isomorphic to $E_6^{p,0}(B)\otimes G$ where $G$ is the kernel of the differential $$d^0_6: E_6^{0,q}(B)\to E_6^{6,q-5}(B).$$

Write some generators of $H^\ast(\K(\mathbb Z,5);\mathbb Z_2)$ as follows:
\[
\begin{aligned}
a_7&=\Sq^2l_5, &
a_8&=\Sq^3l_5, &
a_9&=\Sq^4l_5,\\
a_{11}&=\Sq^4\Sq^2l_5, &
a_{12}&=\Sq^5\Sq^2l_5, &
a_{13}&=\Sq^6\Sq^2l_5,\\
a_{14}&=\Sq^6\Sq^3l_5, &
a_{15}&=\Sq^7\Sq^3l_5 .
\end{aligned}
\]
Consequently, in degrees \(5\le q\le15\), one has
\begin{table}[htbp]
\centering
\small
\setlength{\tabcolsep}{5pt}
\renewcommand{\arraystretch}{1.1}
\begin{tabular}{c|cccccc}
\hline
$q$ & $5$ & $6$ & $7$ & $8$ & $9$ & $10$ \\
\hline
$E^{0,q}_6$
& $\langle l_5\rangle$
& $0$
& $\langle a_7\rangle$
& $\langle a_8\rangle$
& $\langle a_9\rangle$
& $\langle l_5^2\rangle$ \\
\hline
$q$ & $11$ & $12$ & $13$ & $14$ & $15$ & \\
\hline
$E_6^{0,q}$
& $\langle a_{11}\rangle$
& $\langle a_{12},l_5a_7\rangle$
& $\langle a_{13},l_5a_8\rangle$
& $\langle a_{14},l_5a_9,a_7^2\rangle$
& $\langle a_{15},l_5^3,a_7a_8\rangle$
& \\
\hline
\end{tabular}
\end{table}

By Lemma \ref{KZ5SSS}, 
$
\bar d_6(a_i)=0$ for $7\le i\le 15$.
By naturality, 
\[
d_6(a_i)=0 \quad \text{for $7\le i\le 15$,}
\]
\[
d_6(a_7^2)=0 \quad \text{and}\quad d_6(a_7a_8)=0.
\]
Moreover, by
$d_6(l_5)=x^3,$ we have
$$d_6(l_5a_i)=x^3a_i \quad \text{for $7\le i\le 9$},$$
and 
$$d_6(l_5^2)=2x^3=0,\quad d_6(l_5^3)=3x^3=x^3.$$
By the edge homomorphism \cite[15.28]{Sw} of the Serre spectral sequence, 
 there exist $u\in H^{7}(B;\mb{Z}_2)$ and $v\in H^{9}(B;\mb{Z}_2)$ so that
$$i^\ast (u)=\mathrm{Sq}^2l_5\quad i^\ast (\mm{Sq}^1u)=\mathrm{Sq}^3l_5\quad i^\ast (v)=\mathrm{Sq}^4l_5.$$
Combining Steenrod operations with the multiplication of the Serre spectral sequence, we obtain Table \eqref{HBZ2}.

Since $\mathrm{Sq}^1 (ux)=x\mathrm{Sq}^1u$ and $i^\ast(\mathrm{Sq}^1v)=l_5^2=i^\ast(\mathrm{Sq}^2\mathrm{Sq}^1u)$, it follows that
$\mathrm{Sq}^1:H^9(B,\mathbb{Z}_2)\to H^{10}(B,\mathbb{Z}_2)$ is an  isomorphism. Therefore we may choose  the element  $v\in H^9(B;\mb{Z}_2)$ so that
 $i^\ast(v)=\mm{Sq}^4l_5$ and in additional satisfying that
	$\mathrm{Sq}^1v=\mathrm{Sq}^2\mathrm{Sq}^1u$.
\end{proof}

We now determine the $\mathscr A_2$-module structure of
$H^*(B;\mathbb Z_2)$ in the range up to dimension $15$.
It is enough in the
range under consideration to determine the actions of
$\Sq^1,\Sq^2$, and $\Sq^4$ on the generators $u$ and $v$.

\begin{proposition}\label{Sq2uSq2v}
$\Sq^2u=ux,$ $\Sq^2v\neq0.$
\end{proposition}

\begin{proof}
Consider the Serre spectral sequence for the homotopy
fibration
\[
B\longrightarrow \CP^\infty
\xrightarrow{\,r_s\,} \K(\mathbb Z,6),
\]
whose $E_2$-term is
$
E_2^{p,q}
=
H^p\bigl(\K(\mathbb Z,6);
H^q(B;\mathbb Z_2)\bigr).
$

First consider $u\in H^7(B;\mathbb Z_2)$. By degree reasons,
  we have
\[
E_6^{0,7}=E_2^{0,7}
=\mathbb Z_2\langle u\rangle.
\]
We claim that
\[
d_6(u)=l_6x\in E_6^{6,2}.
\]
It is enough to show that $d_6(u)$ is nonzero.  Suppose otherwise.
Since
\[
H^7(\CP^\infty;\mathbb Z_2)=0,
\]
the class $u$ cannot survive to $E_\infty$.  By degree reasons, its
only remaining possible nonzero differential is
\[
d_8:E_8^{0,7}\longrightarrow E_8^{8,0}.
\]
Thus $d_8(u)$ would have to be nonzero.

On the other hand,
$r_s^*(l_6)=x^3,$
and hence
\[
r_s^*(\Sq^2l_6)
=
\Sq^2(x^3)
=
x^4.
\]
It follows from \cite[15.29]{Sw} that the class
$\Sq^2l_6\in E_2^{8,0}$ survives to $E_\infty$, so that
$E_\infty^{8,0}=E_2^{8,0}.$
This contradicts 
$d_8(u)\ne 0$.  Therefore
\[
d_6(u)=l_6x.
\]

The compatibility of Steenrod operations with the Serre spectral
sequence \cite[Theorem~2.17 and (2.79)]{Si} now gives
\[
d_6(\Sq^2u)
=
\Sq^2d_6(u)
=
l_6x^2\neq0,
\]
where only the component in bidegree $(6,4)$ is relevant.
Thus,
\[
\Sq^2u\neq0.
\]

Since
$H^9(B;\mathbb Z_2)
=
\mathbb Z_2\langle ux,v\rangle,$
write
\[
\Sq^2u=a\,ux+b\,v,
\quad
a,b\in\mathbb Z_2.
\]
Applying
$
i^*:H^*(B;\mathbb Z_2)
\longrightarrow
H^*(\K(\mathbb Z,5);\mathbb Z_2)$
and 
\[
i^*(u)=\Sq^2l_5,
\qquad
i^*(v)=\Sq^4l_5,
\qquad
i^*(x)=0,
\]
we obtain
\[
i^*(\Sq^2u)
=
\Sq^2\Sq^2l_5
=
\Sq^3\Sq^1l_5
=
0,
\]
where $\Sq^1l_5=0$.  Hence
 $b=0$.  Since $\Sq^2u\neq0$, we must have 
$\Sq^2u=ux.$

We next prove that $\Sq^2v\neq0$.  Since
\[
E_6^{0,9}
=
\mathbb Z_2\langle ux,v\rangle
\]
and
\[
d_6(ux)=d_6(u)x=l_6x^2\neq0,
\]
the one-dimensional group $E_8^{0,9}=E_7^{0,9}$ is generated by
\[
w=v+\epsilon ux
\]
for some $\epsilon\in\mathbb Z_2$.

As above, $w$ cannot survive to $E_\infty$, because
$H^9(\CP^\infty;\mathbb Z_2)=0$.  Its first possible nonzero
differential is therefore
\[
d_8:E_8^{0,9}\longrightarrow E_8^{8,2},
\]
and the same argument using
\[
r_s^*(\Sq^4l_6)=\Sq^4x^3=x^5
\]
shows that
\[
d_8(w)=\Sq^2l_6\cdot x\neq0.
\]
Applying $\Sq^2$ on the spectral sequence gives
\[
d_8(\Sq^2w)
=
\Sq^2d_8(w)
=
\Sq^2l_6\cdot x^2\neq0.
\]
Hence
$\Sq^2w\neq0.$
Finally, since
\[
\Sq^2(ux)
=
(\Sq^2u)x+\Sq^1u\,\Sq^1x+u\,\Sq^2x
=
ux^2+ux^2
=
0,
\]
we have
\[
\Sq^2w
=
\Sq^2(v+\epsilon ux)
=
\Sq^2v.
\]
Consequently
$\Sq^2v\neq0.$
\end{proof}

To determine $\Sq^2v$ and $\Sq^4v$ explicitly, we restrict the
fibration
\[
\K(\mathbb Z,5)\longrightarrow B\longrightarrow\CP^\infty
\]
to the $4$-skeleton $\CP^2\subset\CP^\infty$.  Since the composite
\[
\CP^2\longrightarrow\CP^\infty
\longrightarrow \K(\mathbb Z,6)
\]
is null-homotopic, the restricted fibration is trivial.  Thus we have a
homotopy-commutative diagram
\begin{equation}\label{fib2}
\begin{tikzcd}
\K(\mathbb Z,5) \ar[r] \ar[d]
&
\K(\mathbb Z,5) \ar[r] \ar[d,"i"]
&
\K(\mathbb Z,5) \ar[d]
\\
\K(\mathbb Z,5)\times\CP^2
\ar[r,"h"] \ar[d]
&
B \ar[d] \ar[r]
&
\mathrm{PK}(\mathbb Z,6) \ar[d]
\\
\CP^2 \ar[r]
&
\CP^\infty \ar[r]
&
\K(\mathbb Z,6)
\end{tikzcd}
\end{equation}
For the product fibration, the spectral sequence 
$$
E_2^{p,q}(\times)=H^{p}(\CP^2;H^{q}(\K(\mathbb Z,5);\mathbb Z_2))
\Longrightarrow
H^{p+q}(\K(\mathbb Z,5)\times\CP^2;\mathbb Z_2)$$
collapses at
$E_2$.
The map $h$ induces a morphism between the corresponding Serre
spectral sequences.

\begin{lemma}\label{injective}
The homomorphism
\[
h^*:H^*(B;\mathbb Z_2)
\longrightarrow
H^*(\K(\mathbb Z,5)\times\CP^2;\mathbb Z_2)
\]
is injective.
\end{lemma}

\begin{proof}
It is enough to show that the induced homomorphism on the associated graded
groups
$
h^*:E_\infty^{p,q}(B)
\longrightarrow
E_\infty^{p,q}(\times)$
is injective.

By Proposition~\ref{convergence},
\[
E_7^{p,q}(B)=E_\infty^{p,q}(B),
\]
and by Lemma~\ref{vanish0} we only need to consider $p\le4$.
In this range,
\[
E_7^{p,q}(B)
=
\ker\!\left(
d_6:E_6^{p,q}(B)
\longrightarrow E_6^{p+6,q-5}(B)
\right),
\]
and hence
\[
E_7^{p,q}(B)
\subseteq
E_6^{p,q}(B)
=
E_2^{p,q}(B).
\]
Moreover, for $p\le 4$,
\[
E_2^{p,q}(B)
\cong
H^p(\CP^2;\mathbb Z_2)
\otimes
H^q(\K(\mathbb Z,5);\mathbb Z_2)=
E_2^{p,q}(\times).
\]
Under this identification, by $
E_\infty^{p,q}(\times)=E_2^{p,q}(\times),$ we have 
$$h^\ast:E_7^{p,q}(B)\to E_7^{p,q}(\times)$$
 is the inclusion. Thus the homomorphism on $E_\infty$-terms is injective.  This finishes the proof.
\end{proof}

\begin{proposition}\label{huv}
The classes $u$ and $v$ satisfy
\[
h^*(u)=\Sq^2l_5+xl_5,
\qquad
h^*(v)=\Sq^4l_5+x^2l_5.
\]
\end{proposition}

\begin{proof}
Since
$i^*(u)=\Sq^2l_5$ and $i^*(v)=\Sq^4l_5,$
degree considerations give
\[
h^*(u)=\Sq^2l_5+cxl_5
\]
and
\[
h^*(v)
=
\Sq^4l_5
+a x^2l_5
+b x\Sq^2l_5
\]
for some $a,b,c\in\mathbb Z_2$.

We first determine $c$.  By Proposition~\ref{Sq2uSq2v},
$\Sq^2u\neq0$, and hence Lemma~\ref{injective} implies
\[
h^*(\Sq^2u)\neq0.
\]
Using the Adem relation
$
\Sq^2\Sq^2=\Sq^3\Sq^1$
and $\Sq^1l_5=0$, we obtain
\[
\begin{aligned}
h^*(\Sq^2u)
&=
\Sq^2h^*(u)\\
&=
\Sq^2(\Sq^2l_5+cxl_5)\\
&=
c\,\Sq^2(xl_5).
\end{aligned}
\]
Therefore $c=1$, and
$h^*(u)=\Sq^2l_5+xl_5.$

Next we determine $b$.  Using
$$
\Sq^1v=\Sq^2\Sq^1u,$$
together with the Cartan and Adem relations, we have
\[
\begin{aligned}
h^*(\Sq^1v)
&=
\Sq^1h^*(v)\\
&=
l_5^2+b\,x\Sq^3l_5,
\end{aligned}
\]
whereas
\[
\begin{aligned}
h^*(\Sq^2\Sq^1u)
&=
\Sq^2\Sq^1h^*(u)\\
&=
\Sq^2\Sq^3l_5\\
&=
(\Sq^5+\Sq^4\Sq^1)l_5\\
&=
l_5^2.
\end{aligned}
\]
Thus
$b\,x\Sq^3l_5=0,$
and hence $b=0$.

It remains to determine $a$.  Since $\Sq^2v\neq0$,
Lemma~\ref{injective} gives
\[
h^*(\Sq^2v)\neq0.
\]
With $b=0$,
\[
\begin{aligned}
h^*(\Sq^2v)
&=
\Sq^2(\Sq^4l_5+a x^2l_5)\\
&=
a\,x^2\Sq^2l_5.
\end{aligned}
\]
Here
$\Sq^2\Sq^4l_5
=
(\Sq^6+\Sq^5\Sq^1)l_5
=
0$.  Hence $a=1$.
\end{proof}

We can now identify the remaining Steenrod operations.

\begin{proposition}\label{Sq2vSq4v}
$\Sq^2v=ux^2,$ $\Sq^4v=\Sq^6u.$
\end{proposition}

\begin{proof}
By Proposition~\ref{huv} and the Cartan formula,
\[
\begin{aligned}
h^*(\Sq^2v)
&=
\Sq^2(\Sq^4l_5+x^2l_5)\\
&=
x^2\Sq^2l_5.
\end{aligned}
\]
On the other hand, since $x^3=0$ in $H^*(\CP^2;\mathbb Z_2)$,
\[
\begin{aligned}
h^*(ux^2)
&=
(\Sq^2l_5+xl_5)x^2\\
&=
x^2\Sq^2l_5.
\end{aligned}
\]
The injectivity of $h^*$ therefore yields
$\Sq^2v=ux^2.$

Similarly,
\[
\begin{aligned}
h^*(\Sq^4v)
&=
\Sq^4(\Sq^4l_5+x^2l_5)\\
&=
\Sq^6\Sq^2l_5+x^2\Sq^4l_5,
\end{aligned}
\]
where we used the Adem relation
$
\Sq^4\Sq^4
=
\Sq^7\Sq^1+\Sq^6\Sq^2$
and $\Sq^1l_5=0$.
On the other hand,
\[
\begin{aligned}
h^*(\Sq^6u)
&=
\Sq^6(\Sq^2l_5+xl_5)\\
&=
\Sq^6\Sq^2l_5+x^2\Sq^4l_5.
\end{aligned}
\]
Hence, again by the injectivity of $h^*$,
$\Sq^4v=\Sq^6u.$
\end{proof}

\begin{proposition}
$
\Sq^2\Sq^4u=\Sq^6u,
$ $\Sq^5\Sq^1u=0$.
\end{proposition}

\begin{proof}
By the Adem relation,
\[
\Sq^2\Sq^4=\Sq^5\Sq^1+\Sq^6.
\]
It therefore suffices to show that $\Sq^5\Sq^1u=0$.
By Proposition~\ref{huv},
\[
h^*(\Sq^5\Sq^1u)
=\Sq^5\Sq^1(\Sq^2l_5+xl_5)
=\Sq^5\Sq^3l_5\\
=0
\]
where the last equality follows from the Adem relation
$\Sq^5\Sq^3=0$.
By the injectivity of $h^*$, we conclude 
$\Sq^5\Sq^1u=0$.
\end{proof}

Consequently, the relevant Steenrod operations on the generators
$x$, $u$ and $v$ are summarized as follows:
\begin{table}[!h]
\centering
\caption{The actions of $\Sq^1$, $\Sq^2$, and $\Sq^4$ on $u$ and $v$.}
\label{SqBtable}
\begin{tabular}{c|ccc}
\hline
 & $\Sq^1$ & $\Sq^2$ & $\Sq^4$\\
\hline
$x$
& $0$
& $x^2$
& $0$
\\
$u$
& $\Sq^1u$
& $ux$
& $\Sq^4u$
\\
$v$
& $\Sq^2\Sq^1u$
& $ux^2$
& $\Sq^6u$
\\
\hline
\end{tabular}
\vspace{0.5em}
\parbox{\linewidth}{\small \centering In particular, $\Sq^2\Sq^4u=\Sq^6u$, $\Sq^5\Sq^1u=0$.}
\end{table}

\subsection{$H^\ast(B;\mathbb Z_3)$ and its module over the Steenrod
algebra}\quad

  As in the
mod $2$ case, we work directly with the Serre spectral sequence
of the fibration
\[
\K(\mathbb Z,5)\longrightarrow B\longrightarrow \CP^\infty,
\]
to determine the mod $3$ cohomology of $B$ and its module
structure over the mod $3$ Steenrod algebra $\mathcal A$.

We first recall the mod $3$ cohomology of the fiber.  Let
\[
l_5\in H^5(\K(\mathbb Z,5);\mathbb Z_3)
\]
denote the fundamental class, let $\mathcal P^i$ denote the
$i$th reduced power, and let $\beta$ denote the mod $3$
Bockstein.

By the Cartan--Serre description, in degrees at most $14$ the
nonzero reduced cohomology groups of $\K(\mathbb Z,5)$ are
\[
\begin{aligned}
H^5(\K(\mathbb Z,5);\mathbb Z_3)
 &=\langle l_5\rangle,\\
H^9(\K(\mathbb Z,5);\mathbb Z_3)
 &=\langle \mathcal P^1l_5\rangle,\\
H^{10}(\K(\mathbb Z,5);\mathbb Z_3)
 &=\langle \beta\mathcal P^1l_5\rangle,\\
H^{13}(\K(\mathbb Z,5);\mathbb Z_3)
 &=\langle \mathcal P^2l_5\rangle,\\
H^{14}(\K(\mathbb Z,5);\mathbb Z_3)
 &=
 \langle
 \beta\mathcal P^2l_5,\,
 l_5\mathcal P^1l_5
 \rangle .
\end{aligned}
\]
All other reduced cohomology groups in degrees at most $14$
vanish.  

We shall also use the mod $3$ Adem relation
\begin{equation}\label{P1P1}
\mathcal P^1\mathcal P^1=-\mathcal P^2.
\end{equation}

Consider now the Serre spectral sequence
\[
E_2^{p,q}(B)
=
H^p(\CP^\infty;\mathbb Z_3)
\otimes
H^q(\K(\mathbb Z,5);\mathbb Z_3)
\Longrightarrow
H^{p+q}(B;\mathbb Z_3).
\]
with the transgression 
\[
d_6(l_5)=x^3.
\]

\begin{lemma}\label{vanish1}
$
E_6^{p,q}(B)=E_2^{p,q}(B).$ In total degrees at most $14$, the only nonzero $d_6$-differentials
are those forced by
\[
d_6(l_5)=x^3
\]
and the derivation property.  In particular,
\[
d_6(\mathcal P^1l_5)
=
d_6(\beta\mathcal P^1l_5)
=
d_6(\mathcal P^2l_5)
=
d_6(\beta\mathcal P^2l_5)
=0,
\]
whereas
\[
d_6(l_5\mathcal P^1l_5)
=
x^3\mathcal P^1l_5.
\]
Furthermore,
\[
E_\infty^{p,q}(B)=E_7^{p,q}(B),
\qquad
E_7^{p,q}(B)=0
\quad\text{for }p\ge5.
\]
\end{lemma}

\begin{proof}
Consider first the universal path-loop fibration
\[
\K(\mathbb Z,5)
\longrightarrow
\mathrm{PK}(\mathbb Z,6)\simeq *
\longrightarrow
\K(\mathbb Z,6).
\]
Its fundamental class satisfies
$\bar d_6(l_5)=l_6.$
By the transgression theorem \cite[Corollary 6.9]{McCleary}, the classes obtained from $l_5$
by nontrivial Steenrod operations survive beyond the $E_6$-page.
Hence, in the range under consideration,
\[
\bar d_6(\mathcal P^1l_5)
=
\bar d_6(\beta\mathcal P^1l_5)
=
\bar d_6(\mathcal P^2l_5)
=
\bar d_6(\beta\mathcal P^2l_5)
=0.
\]
Naturality gives the corresponding equalities in the spectral
sequence for $B$.
Since $d_6$ is a derivation,
$
d_6(l_5\mathcal P^1l_5)
=
x^3\mathcal P^1l_5.$

As in the mod $2$ case, through degree $14$, we have
\[
E_7^{*,*}(B)
\cong
\frac{\mathbb Z_3[x]}{(x^3)}
\otimes
\left\langle
1,\,
\mathcal P^1l_5,\,
\beta\mathcal P^1l_5,\,
\mathcal P^2l_5,\,
\beta\mathcal P^2l_5
\right\rangle.
\]
  Hence only the columns
$p=0,2,4$ remain nonzero.

For $r\ge7$, every differential $d_r$ starting in one of these
columns has target in a column $p\ge7$, which is already zero,
and no differential can enter the $p=0,2,4$ columns.  Therefore
$E_\infty=E_7.$
\end{proof}

This immediately gives the mod $3$ cohomology of $B$ in the range.

\begin{lemma}\label{3cohomology}
Through degree $14$, the nonzero positive-dimensional
cohomology groups of $B$ are
\begin{center}
\begin{tabular}{|c|c|c|c|c|c|c|c|c|c|c|c|c|c|c|c|}
\hline
 i   &2 &4&9&10 & 11& 12 & 13& 14\\
\hline
$H^i(B;\mathbb{Z}_3)$ &$x$ &$x^2$& $ z_9 $ &$\beta z_9 $& $ xz_9$& $x\beta z_9$ &$z_{13}$, $x^2z_{9}$& $\beta z_{13} $, $x^2\beta z_9$\\
\hline
\end{tabular}
\end{center}
where $z_9$ satisfies
\[
i^*(z_9)=\mathcal P^1l_5,
\]
and we set
$z_{13}:=\mathcal P^1z_9.$
In particular,
\[
i^*(z_{13})
=
\mathcal P^1\mathcal P^1l_5
=
-\mathcal P^2l_5.
\]
\end{lemma}

\begin{proof}
The generator description follows directly from
Lemma~\ref{vanish1}.

Choose $z_9$ representing the class
$\mathcal P^1l_5\in E_\infty^{0,9}$.  Then
\[
i^*(z_9)=\mathcal P^1l_5.
\]
Since
$\mathcal P^1\mathcal P^1=-\mathcal P^2$
and $\mathcal P^2l_5\neq0$, we have
\[
i^*(\mathcal P^1z_9)
=
-\mathcal P^2l_5\neq0.
\]
Thus
$z_{13}:=\mathcal P^1z_9$
is nonzero and is linearly independent of $x^2z_9$.

Note that
\[
i^*(\beta z_9)
=
\beta\mathcal P^1l_5\neq0 \quad\text{and}\quad i^*(\beta z_{13})
=
-\beta\mathcal P^2l_5\neq0.
\]
Hence the displayed classes give the asserted bases.
\end{proof}

We next determine a Steenrod operation:
\begin{lemma}\label{P1betaz9}
There exists $m\in\mathbb Z_3$ such that
\[
\mathcal P^1\beta z_9
=
-\beta z_{13}+mx^2\beta z_9.
\]
\end{lemma}

\begin{proof}
By the mod $3$ Adem relation
$
\mathcal P^1\beta\mathcal P^1
=
\beta\mathcal P^2+\mathcal P^2\beta,$
we obtain
\[
i^*(\mathcal P^1\beta z_9)
=\mathcal P^1\beta\mathcal P^1l_5
=\beta\mathcal P^2l_5
=-i^*(\beta z_{13})
\]
where we used $\beta l_5=0$.
Since
$
H^{14}(B;\mathbb Z_3)
=
\langle\beta z_{13},x^2\beta z_9\rangle$
and
\[
\ker\bigl(
i^*:H^{14}(B;\mathbb Z_3)
\to H^{14}(\K(\mathbb Z,5);\mathbb Z_3)
\bigr)
=
\langle x^2\beta z_9\rangle,
\]
the desired identity follows.
\end{proof}

\subsection{The integral homology of $B$}
\quad

We now determine the integral homology of $B$ directly from the
mod-$p$ calculations above, together with the Serre spectral
sequence of
\[
\K(\mathbb Z,5)\longrightarrow B\longrightarrow\CP^\infty .
\]

Recall that $B$ is the homotopy fibre of the map
\[
r_s:\CP^\infty\longrightarrow \K(\mathbb Z,6)
\]
satisfying
\[
r_s^*(l_6)=s x^3,
\qquad
\gcd(s,6)=1.
\]

\begin{proposition}\label{HBZ}
Up to degree $15$, the nonzero integral cohomology groups of $B$ are
\begin{center}
\begin{tabular}{|c|c|c|c|c|c|c|c|c|c|c|c|c|c|c|c|}
\hline
 i & 2& 4 &6 &8 & 10  \\
\hline
$H^i$ &$\mathbb{Z}\langle x \rangle$ &$\mathbb{Z}\langle x^2 \rangle$& $\mathbb{Z}_s\langle x^3 \rangle $& $\mathbb{Z}_2 \oplus \mathbb{Z}_s\langle x^4\rangle$ &$\mathbb{Z}_s\langle x^5\rangle \oplus \mathbb{Z}_2^2 \oplus \mathbb{Z}_3$\\
\hline
\end{tabular}
\begin{tabular}{|c|c|c|c|c|c|c|c|c|c|c|c|c|c|c|c|}
\hline
 i & 12& 14 &15   \\
\hline
$H^i$ &$\mathbb{Z}_s\langle x^6\rangle \oplus \mathbb{Z}^3_2 \oplus \mathbb{Z}_3$ &$\mathbb{Z}_s\langle x^7\rangle  \oplus  \mathbb{Z}^3_2 \oplus \mathbb{Z}^2_3 \oplus \mathbb{Z}_5$& $\mathbb{Z}_2\oplus \cdots $\\
\hline
\end{tabular}
\end{center}
\end{proposition}

\begin{proof}
We determine separately the $p$-primary
components.

\medskip
\noindent
\textbf{Step 1: the free part and the $\mathbb Z_s$-summands.}

Consider the integral cohomological Serre spectral sequence
\[
E_2^{p,q}
=
H^p\bigl(\CP^\infty;
H^q(\K(\mathbb Z,5);\mathbb Z)\bigr)
\Longrightarrow
H^{p+q}(B;\mathbb Z).
\]
By naturality with respect to the universal path-loop fibration
\[
\K(\mathbb Z,5)
\longrightarrow
\mathrm{PK}(\mathbb Z,6)\simeq *
\longrightarrow
\K(\mathbb Z,6),
\]
and the defining relation
$r_s^*(l_6)=s x^3,$
we have
\[
d_6(l_5)=s x^3,
\]
and thus
\[
d_6(x^k l_5)
=
x^k d_6(l_5)
=
s x^{k+3}
\qquad (k\ge0).
\]
 For $k\ge0$, the differential
$d_6:
E_6^{2k,5}\longrightarrow E_6^{2k+6,0}$
is multiplication by $s$:
\[
\mathbb Z\langle x^k l_5\rangle
\xrightarrow{\;\times s\;}
\mathbb Z\langle x^{k+3}\rangle .
\]
Since $s\neq0$, this map is injective.  Consequently,
$
E_7^{2k,5}=0,$
while
\[
E_7^{2k+6,0}
\cong
\mathbb Z_s
\]
with generator represented by $x^{k+3}$.

In particular, in total degrees at most $14$,
\[
E_\infty^{2,0}=E_7^{2,0}\cong\mathbb Z\langle x\rangle,
\qquad
E_\infty^{4,0}=E_7^{4,0}\cong\mathbb Z\langle x^2\rangle,
\]
and
\[
E_7^{2j,0}
\cong
\mathbb Z_s\langle x^j\rangle,
\qquad
3\le j\le7.
\]
The source of an incoming
differential for $E_r^{2j,0}$ is $E_r^{2j-r,r-1}$.

For $7\le r\le14$, 
$H^{r-1}(\K(\mathbb Z,5))$ has no
$p$-primary torsion for any prime $p$ coprime to $6$.
Consequently, $E_2^{2j-r,r-1}$, and hence
$E_r^{2j-r,r-1}$, has only $2$- and $3$-primary
components. Since
$
E_7^{2j,0}=\Z_s$
 and $\gcd(s,6)=1$, every incoming
differential for
$ E_r^{2j,0}$
vanishes.
 Hence
$$E_\infty^{2j,0}=E_7^{2j,0}
\cong
\mathbb Z_s\langle x^j\rangle,
\qquad
3\le j\le7.$$
Therefore, for $3\le j\le7$,
\[
H^{2j}(B)
\supseteq
\mathbb Z_s\langle x^j\rangle.
\]

Moreover, for rational cohomology,  the  differential
\[
d_6(l_5)=s x^3
\]
 kills both $l_5$ and all classes $x^j$ with $j\ge3$.
Thus
\[
H^*(B;\mathbb Q)
\cong
\mathbb Q\{1,x,x^2\},
\]
so the only free integral cohomology groups are $H^0$, $H^2$ and $H^4$.

\medskip
\noindent
\textbf{Step 2: the $2$-primary component.}

By Table \eqref{HBZ2}, the operation
$$\Sq^1:H^i(B;\mb{Z}_2)\to H^{i+1}(B;\mb{Z}_2)$$
is isomorphic for $i=7,9,11$, injective for $i=13$.  For $i=14$, we have 
$$\Sq^1(\Sq^6\Sq^1u)=\Sq^7\Sq^1u.$$
Furthermore, it is also injective that
$$\Sq^1:H^{15}(B;\mb{Z}_2)/\{\Sq^7\Sq^1u\}\to H^{16}(B;\mb{Z}_2).$$

The operation $\Sq^1$ is the mod $2$  Bockstein operation.
From \cite[Prop. 3E.3]{Hatcher}, we obtain that a summand $\mb{Z}_{2^k}\subset H^n(B;\mb{Z})$ corresponds an element $\alpha\in H^n(B;\mb{Z}_2)$ such that 
$$\text{$\alpha\notin\mathrm{im}\{\Sq^1\}$ for $k>1$; $\alpha\in\mathrm{im}\{\Sq^1\}$ for $k=1$.}$$   

Therefore, Table \eqref{HBZ2} gives the $2$-local cohomology groups.

\medskip
\noindent
\textbf{Step 3: the $3$-primary component.}

 Lemma \ref{3cohomology} and \cite[Prop. 3E.3]{Hatcher} gives the $3$-local cohomology groups.

\medskip
\noindent
\textbf{Step 4: primes $p\ge5$.}

For $p\ge7$, there is no additional $p$-primary class in the
cohomology of $\K(\mathbb Z,5)$ in total degrees at most $14$.
Indeed, the first mod-$p$ Steenrod class obtained from $l_5$
is
$
\mathcal P^1l_5$ with
$|\mathcal P^1l_5|
=
5+2(p-1)\ge17.$

It remains to consider $p=5$.  In this case,
\[
\mathcal P^1l_5
\in
H^{13}(\K(\mathbb Z,5);\mathbb Z_5),\quad \beta\mathcal P^1l_5
\in
H^{14}(\K(\mathbb Z,5);\mathbb Z_5)
\]
and
\[
H^{13}(\K(\mathbb Z,5);\mathbb Z_{(5)})=0.
\]

We claim that $\mathcal P^1l_5$ survives in the mod-$5$ Serre
spectral sequence for $B$.  If $5\mid s$, this follows immediately
from
\[
d_6(l_5)= s\,x^3=0.
\]
The transgression theorem \cite[Corollary 6.9]{McCleary} gives
$$d_{14}(\mathcal P^1l_5)
=0.$$
If $5\nmid s$, then
\[
d_6(x^4l_5)= s\,x^7,
\]
so $x^7=0$ on the $E_7$-page.  The transgression theorem 
shows that the only possible differential on
$\mathcal P^1l_5$ is
\[
d_{14}(\mathcal P^1l_5)
=
\mathcal P^1( s\,x^3)
=
3 s\,x^7,
\]
which is already zero on the $E_{14}$-page.  Hence
$\mathcal P^1l_5$ survives in either case.

Let $e_{13}\in H^{13}(B;\Z_5)$ satisfy $\mathcal P^1l_5=i^\ast(e_{13})\in H^{13}(K(\mathbb Z,5);\mathbb Z_5)$. Since
$
\beta\mathcal P^1l_5
\neq0$, we have 
$$\beta e_{13}\ne 0.$$
Thus, the Bockstein detects an order-$5$ summand
$\mathbb Z_5\subset H^{14}(B;\mathbb Z).$

When $5\nmid s$, this is the entire $5$-primary component apart
from the possible $5$-part of $\mathbb Z_s$, which is trivial.

Suppose now that $5\mid s$, and write
$
5^a\Vert s.$ 
Assume there is a nontrivial extension, i.e., 
$
\mathbb Z_{5^{a+1}}\cong H^{14}(B;\mathbb Z)_{(5)}.$ By the universal coefficient theorem,
$$H_{13}(B;\Z_{(5)})\cong H^{14}(B;\Z_{(5)})$$ and
$$ H^{13}(B;\Z_5)\cong \mm{Hom}( H_{13}(B;\Z_{(5)}), \Z_5)\cong\Z_5.$$
Here we note that $H^{13}(B;\Z_{(5)})=0$.

However, by the mod $5$ Serre spectral sequence,
$$ H^{13}(B;\Z_5)\cong E_\infty^{8,5}\oplus E_\infty^{0,13}\cong \Z_5\oplus \Z_5.$$
This is a contradiction. Therefore,
$H^{14}(B)_{(5)}
\cong
\mathbb Z_5^a\oplus\mathbb Z_5.$
\end{proof}

\section{Fundamental classes of Bazaikin spaces in $H_*(B)$}
\label{oriented}

Recall 
$
\mathcal B_{\bold q}
=
S^1\backslash \SU(6)/\Sp(3)$.
Put
\[
P=\SU(6)/\Sp(3).
\]
The $S^1$-action on $P$ has kernel $\mathbb Z_2$.
Hence the induced action of
$S^1/\mathbb Z_2\cong S^1$ is free and determines a principal
circle bundle
\[
S^1\longrightarrow P
\xrightarrow{\pi }\mathcal B_{\bold q}.
\]
Let
$
c:\mathcal B_{\bold q}\longrightarrow BS^1=\CP^\infty\simeq \K(\Z,2)$
be its classifying map, and
$x$
be the corresponding generator in $H^2(\mathcal B_{\bold q})$.

As discussed in Section~\ref{Normalbg}, the map $c$ admits a
canonical lift
\[
f:\mathcal B_{\bold q}\longrightarrow B
\]
which is a $6$-equivalence.

The Gysin sequence of the circle bundle  contains the
exact sequences
\[
0\longrightarrow H^5(P)
\xrightarrow{\Phi}
H^4(\mathcal B_{\bold q})
\longrightarrow
H^6(\mathcal B_{\bold q})
\longrightarrow0
\]
and
\begin{equation}\label{beta9top9}
0\longrightarrow
H^9(\mathcal B_{\bold q};\mathbb Z)
\xrightarrow{\pi^*}
H^9(P;\mathbb Z)
\longrightarrow
H^8(\mathcal B_{\bold q};\mathbb Z)
\longrightarrow0.
\end{equation}
By \cite{FL},
$
H^*(P)\cong H^*(S^5\times S^9).$
Since
$
H^8(\mathcal B_{\bold q})
\cong\mathbb Z_s$
and $\gcd(s,6)=1$, reduction modulo $2$ or $3$ shows that
\[
\pi^*:
H^9(\mathcal B_{\bold q};\mathbb Z_i)
\xrightarrow{\cong}
H^9(P;\mathbb Z_i),
\qquad i=2,3.
\]

Choose a generator
$
\alpha_5\in H^5(P;\mathbb Z)$
such that
\[
\Phi(\alpha_5)=sx^2\in H^4(\mathcal B_{\bold q}).
\]
Fix the orientation class
\[
[P]\in H_{14}(P;\mathbb Z),
\]
and choose the generator
$
\alpha_9\in H^9(P)$
so that
\[
\langle\alpha_5\alpha_9,[P]\rangle=1.
\]
By \eqref{beta9top9}, there is a generator
$
\beta_9\in H^9(\mathcal B_{\bold q};\mathbb Z)$
satisfying
\[
\pi^*(\beta_9)=s\alpha_9.
\]

Therefore, the orientation of $P$, together with the free circle action,
induces a canonical orientation
\begin{equation}\label{canon_orien}
	[\mathcal B_{\bold q}]
\in H_{13}(\mathcal B_{\bold q};\mathbb Z)
\end{equation}
which duals to $x^2\beta_9\in H^{13}(\mathcal B_{\bold q};\mathbb Z)$.

\subsection{The $9$-skeleton of $P$ and Steenrod operations}\quad

The homotopy fibration
\[
\K(\mathbb Z,5)\longrightarrow B\longrightarrow\CP^\infty
\]
also gives a principal circle fibration
\[
S^1\longrightarrow \K(\mathbb Z,5)
\xrightarrow{i}B
\]
with Euler class $x$.  Pulling it back along $f$ gives the
commutative diagram
\[
\begin{tikzcd}
S^1 \ar[r] \ar[d,equal]
&
P \ar[r,"\pi"] \ar[d,"\bar f"]
&
\mathcal B_{\bold q} \ar[d,"f"]
\\
S^1 \ar[r]
&
\K(\mathbb Z,5) \ar[r,"i"]
&
B 
\end{tikzcd}
\]
Naturality of the Gysin sequence shows that
\begin{equation}\label{barf}
\bar f^*:
H^5(\K(\mathbb Z,5);\mathbb Z)
\xrightarrow{\cong}
H^5(P;\mathbb Z).
\end{equation}
After choosing the sign of the fundamental class
$l_5\in H^5(\K(\mathbb Z,5);\mathbb Z)$ appropriately, we may
assume that
\[
\bar f^*(l_5)=\alpha_5.
\]

Since $P$ is a $4$-connected closed $14$-manifold with
$
H^*(P;\mathbb Z)
\cong H^*(S^5\times S^9;\mathbb Z),$
\cite[Proposition~4.1]{Wall1965} gives a minimal CW decomposition
\[
P\simeq
\bigl(S^5\cup_a e^9\bigr)\cup e^{14},
\]
where
$
a\in\pi_8(S^5)\cong\mathbb Z_{24}.$
Since the top cell does not affect $\pi_8$, we have
\[
\pi_8(P)
\cong
\mathbb Z_{24}/\langle a\rangle.
\]
By \cite{FL}, $\pi_8(P)=0$, and hence $a$ is a generator of
$\pi_8(S^5)$.  

\begin{lemma}\label{4module}
The operations
\[
\Sq^4:
H^5(P;\mathbb Z_2)
\longrightarrow
H^9(P;\mathbb Z_2)
\]
and
\[
\mathcal P^1:
H^5(P;\mathbb Z_3)
\longrightarrow
H^9(P;\mathbb Z_3)
\]
are isomorphisms.
\end{lemma}

\begin{proof}
Recall $\pi_7(S^4)=\mb{Z}\oplus \mb{Z}_{12}$. The suspension homomorphism
\[
\Sigma:\pi_7(S^4)\longrightarrow\pi_8(S^5)=\mb{Z}_{24}
\]
is surjective. Let
$
\psi:S^7\longrightarrow S^4$
be the quaternionic Hopf map, which generates the free summand
of $\pi_7(S^4)$ and has Hopf invariant one.  Its suspension
$\Sigma\psi$ generates $\pi_8(S^5)\cong\mathbb Z_{24}$.

 Since the attaching map $a\in \pi_8(S^5)$ is a generator, there exists $n\in \mb{Z}$ such that $\mathrm{gcd}(n,24)=1$ and $\Sigma n\psi=\alpha$. 

Let
$C_\psi=S^4\cup_{n\psi} e^8$, and let $w_i$ generate $H^i(C_\psi)$ for $i=4,8$. Since the Hopf invariant of $n\psi$ equals $n$, we have 
$$w_4^2=nw_8,\quad \gcd(n,24)=1.$$  
This implies that the following oprerations $$\Sq^4:H^4(C_\psi;\mb{Z}_2)\to H^8(C_\psi;\mb{Z}_2)$$ 
$$\mathcal P^1:H^4(C_\psi;\mb{Z}_3)\to H^8(C_\psi;\mb{Z}_3)$$
are nontrivial. Their sources and targets are one-dimensional,
hence both homomorphisms are isomorphisms.

The $9$-skeleton of $P$ satisfies
\[
P^{(9)}= S^5\cup_\alpha e^9\simeq \Sigma C_\psi. 
\]
Since the Steenrod operations commute with cohomology suspension, the desired results follow.
\end{proof}

\subsection{The image of the fundamental class}
\quad

We next record the Steenrod operations on
$\mathcal B_{\bold q}$ needed below.
By Section~\ref{Normalbg}, the stable difference bundle \eqref{eq:difference-bundle}
$$\delta_{\bold q}
=
\mathscr N_{\bold q}
-
f^{*}\eta$$
satisfies $w_2(\delta_{\bold q})=w_4(\delta_{\bold q})=0$. This implies 
$$w_2(\mathscr N_{\bold q})=w_2(f^{*}\eta)\quad \text{and} $$
$$w_4(\mathscr N_{\bold q})+w_4(f^{*}\eta)+w_2(\mathscr N_{\bold q})w_2(f^{*}\eta)=0.$$
By \eqref{eq:eta-char}, we have
$$w_2(\mathscr N_{\bold q})\ne 0\quad \text{and} \quad w_4(\mathscr N_{\bold q})=0.$$
Since $\mathscr N_{\bold q}\oplus T\mathcal B_{\bold q}$ is stably trivial, it follows that
\begin{equation}\label{Whiteney_Bq}
	w_2(\mathcal B_{\bold q})\neq0
\quad \text{and} \quad
w_4(\mathcal B_{\bold q})\neq0.
\end{equation}
Thus the mod $2$ Wu classes satisfy
\[
v_2(\mathcal B_{\bold q})\neq0 \quad \text{and} \quad
v_4(\mathcal B_{\bold q})=0.
\]
The mod $2$ Wu formula therefore gives
\begin{equation}\label{SqBq}
\Sq^2(\beta_{11})=\beta_{13},
\qquad
\Sq^2(\beta_9)=0,
\qquad
\Sq^4(\beta_9)=0,
\end{equation}
where the same symbols denote the mod-$2$ reductions of the
corresponding integral generators.

The first mod-$3$ Wu class is the reduction
of the first Pontryagin class. Therefore,
by \eqref{eq:p1Bq} and the mod $3$ Wu formula,
\begin{equation}\label{P1Bq}
\mathcal P^1(\beta_9)
\equiv
k(x^2\beta_9)\pmod 3.
\end{equation}

The following identification of the image of the fundamental
class will be used in the computation of the normal bordism
group.

\begin{theorem}\label{image}
Under the canonical lift
$
f:\mathcal B_{\bold q}\longrightarrow B,$
the mod-$2$ fundamental class satisfies
\[
f_*([\mathcal B_{\bold q}]_2)
=
(vx^2)^\vee,
\]
where $-^\vee$ denotes a homology class dual to the cohomology class
$
-.$

In the mod-$3$ case, 
the coefficient of $(x^2z_9)^\vee$ in
$f_*([\mathcal B_{\bold q}]_3)$ is nonzero.

\end{theorem}

\begin{proof}
By Lemma~\ref{4module} and \eqref{barf},
\[
\bar f^*(\Sq^4l_5)
=
\Sq^4\bar f^*(l_5)
=
\Sq^4\alpha_5\ne 0\in 
H^9(P;\mathbb Z_2).\]
Since
$
i^*(v)=\Sq^4l_5,$
the pullback square gives
\[
\pi^*f^*(v)
=
\bar f^*i^*(v)
=
\Sq^4\alpha_5.
\]
On the other hand, since
$
\pi^*:
H^9(\mathcal B_{\bold q};\mathbb Z_2)
\xrightarrow{\cong}
H^9(P;\mathbb Z_2),$
we obtain 
\begin{equation}\label{fastvbeta9}
f^*(v)=\beta_9.
\end{equation}
Therefore, it follows that 
\[
\begin{aligned}
\left\langle
x^2v,\,
f_*([\mathcal B_{\bold q}]_2)
\right\rangle
&=
\left\langle
f^*(x^2v),\,
[\mathcal B_{\bold q}]_2
\right\rangle\\
&=
\left\langle
x^2\beta_9,\,
[\mathcal B_{\bold q}]_2
\right\rangle
=1.
\end{aligned}
\]

For mod-$3$ coefficients, Lemma~\ref{4module} similarly shows
that the nonzero class
\[
\mathcal P^1l_5\in
H^9(K(\mathbb Z,5);\mathbb Z_3)
\]
pulls back to a generator of $H^9(P;\mathbb Z_3)$.
The pullback square also gives
\begin{equation}\label{fastz_9_beta_9}
	f^*(z_9)=\mathfrak m \beta_9
\end{equation}
for some $0\ne \mathfrak m\in \Z_3$. Then
we obtain
\[
\left\langle
x^2z_9,\,
f_*([\mathcal B_{\bold q}]_3)
\right\rangle
=
\left\langle
\mathfrak m x^2\beta_9,\,
[\mathcal B_{\bold q}]_3
\right\rangle
=\mathfrak m.
\]
Hence the coefficient of $(x^2z_9)^\vee$ in
$f_*([\mathcal B_{\bold q}]_3)$ is nonzero.
\end{proof}

\section{Preliminary computation for $\Omega_{13}^{\mm{O}\langle 8\rangle}(\eta)$}

In this section we use the AHSS to give
the basic group structure of $\Omega_{13}^{\mm{O}\langle 8 \rangle}(\eta)$.
We then use the Hurewicz homomorphism to
identify the resulting summands by geometric invariants.

Recall first the homological AHSS
\[
E^2_{p,q}
=
H_p(X;\pi_q(E))
\Longrightarrow E_{p+q}(X)\cong
\pi_{p+q}(E\wedge X),
\]
with differential
$
d^r:E^r_{p,q}\longrightarrow E^r_{p-r,q+r-1}.$

We then recall two standard facts that will be used repeatedly.

\begin{lemma}\label{d2}
Let $E$ and $X$ be connective spectra.  Suppose that the unit
$\iota:S^0\longrightarrow E$
induces isomorphisms
$
\pi_i(S^0)\xrightarrow{\cong}\pi_i(E)$ for $0\le i\le3.$
Let $(\Sq^2)^\vee$ denote the homology operation dual to $\Sq^2$.
Then:
\begin{enumerate}
\item[(1)]
$
d_2:E^2_{p+2,0}\longrightarrow E^2_{p,1}$
is the composite
\[
H_{p+2}(X;\mathbb Z)
\xrightarrow{\rho_{2}}
H_{p+2}(X;\mathbb Z_2)
\xrightarrow{(\Sq^2)^\vee}
H_p(X;\mathbb Z_2),
\]
where $\rho_{2}$ is induced by reduction modulo $2$.

\item[(2)]
$
d_2:E^2_{p+2,1}\longrightarrow E^2_{p,2}
$
is
$
(\Sq^2)^\vee:
H_{p+2}(X;\mathbb Z_2)
\longrightarrow H_p(X;\mathbb Z_2).
$

\item[(3)]
$
d_2:E^2_{p+2,2}\longrightarrow E^2_{p,3}$
is the composite
\[
H_{p+2}(X;\mathbb Z_2)
\xrightarrow{(\Sq^2)^\vee}
H_p(X;\mathbb Z_2)
\xrightarrow{\gamma_*}
H_p(X;\mathbb Z_{24}),
\]
where $\gamma_*$ is induced by the coefficient homomorphism
\[
\gamma:\mathbb Z_2\longrightarrow\mathbb Z_{24},
\quad \gamma(1)=12.
\]
\end{enumerate}
\end{lemma}

\begin{proof}
Items~(1) and~(2) follow from
\cite{Shen2026Differentials}.  We deduce~(3) from~(2) by
naturality in the coefficient spectrum.

Write $\eta_{\mathrm{st}}\in\pi_1^S$ for the stable Hopf element.  It defines a map
\[
\mu_{\eta_{\mathrm{st}}}:
\Sigma E
\cong S^1\wedge E
\xrightarrow{\eta_{\mathrm{st}}\wedge1_E}
S^0\wedge E
\cong E.
\]
On homotopy groups, this map induces multiplication by
$\eta_{\mathrm{st}}$:
\[
\pi_{q-1}(E)=\pi_q(\Sigma E)
\longrightarrow\pi_q(E).
\]
Using the isomorphisms induced by the unit, the standard
relations in the stable stems give
\[
\eta_{\mathrm{st}}\cdot-:
\pi_1(E)\xrightarrow{\cong}\pi_2(E)
\]
and
\[
\eta_{\mathrm{st}}\cdot\eta_{\mathrm{st}}^2
=
\eta_{\mathrm{st}}^3
=
12\nu
\in\pi_3(E)\cong\mathbb Z_{24},
\]
where $\nu$ is a generator.

The map $\mu_{\eta_{\mathrm{st}}}$ induces a morphism of AHSSs.
After identifying the AHSS for $\Sigma E$ with the
coefficient-degree shift of that for $E$, naturality gives
\[
\begin{tikzcd}[column sep=large,row sep=large]
H_{p+2}(X;\pi_1(E))
  \arrow[r,"d_2"]
  \arrow[d,"(\eta_{\mathrm{st}}\cdot-)_*"',"\cong"]
&
H_p(X;\pi_2(E))
  \arrow[d,"(\eta_{\mathrm{st}}\cdot-)_*"]
\\
H_{p+2}(X;\pi_2(E))
  \arrow[r,"d_2"]
&
H_p(X;\pi_3(E)).
\end{tikzcd}
\]
By~(2), the upper horizontal map is $(\Sq^2)^\vee$.
Under the coefficient identifications, the left vertical
map is the identity and the right vertical map is $\gamma_*$. Thus we finish the proof.
\end{proof}

Let $E$ be a connective spectrum with $\pi_0(E)=\mathbb Z$.  Its
$0$-th Postnikov section gives a map
$
\chi:E\longrightarrow  H\Z,
$
where $H\mathbb Z$ is the integral Eilenberg--Mac Lane
spectrum. The induced homomorphism
$$
\chi_*:E_n(X)\longrightarrow H_n(X;\mathbb Z)$$
is indeed an edge homomorphism of the AHSS. Therefore, 

\begin{lemma}\label{chi}
Let $E$ and $X$ be connective spectra with $\pi_0(E)=\mathbb Z$.
A class
$z\in E_n(X)$
is detected by
$
x\in E^2_{n,0}$
in the AHSS if and only if
$
\chi_*(z)=x
$. 
\end{lemma}

\subsection{Preliminary computation for AHSS}\label{pprimary}\quad

Let $\mm{M}\eta$ and $\mm{MO}\langle 8\rangle$ denote the Thom spectra of the bundles $\eta$ and $\gamma_8$ (see Section \ref{Normalbg}) respectively.
Put
$$\mm M_{\mathscr B}=\mm{MO}\langle 8\rangle\wedge \mm M\eta.$$
 By the Pontryagin-Thom isomorphism, 
 $\Omega_{\ast}^{\mm{O}\langle 8 \rangle}(\eta)\cong \pi_\ast(\mm M_{\mathscr B}).$

We next consider the AHSS
\begin{equation}\label{AHSSforMO8Meta-odd}
E^2_{p,q}
=
H_p(\mm{M}\eta;
\pi_q(\mm{MO}\langle8\rangle))
\Longrightarrow
\pi_{p+q}(\mm{M}_{\mathscr B}).
\end{equation}

We recall the low-dimensional homotopy groups of
$\mm{MO}\langle8\rangle$ \cite{Giambalvo,HoRa1995}.
\begin{table}[!h]
\centering
\caption{$\pi_*(\mm{MO}\langle8\rangle)$}
\label{piMO8}
\begin{tabular}{c|cccccccccc}
\hline
$i$
&0&1&2&3&4&5&6&7&8&9\\
\hline
$\pi_i(\mm{MO}\langle8\rangle)$
&$\mathbb Z$
&$\mathbb Z_2$
&$\mathbb Z_2$
&$\mathbb Z_{24}$
&$0$
&$0$
&$\mathbb Z_2$
&$0$
&$\mathbb Z\oplus\mathbb Z_2$
&$\mathbb Z_2^2$
\\
\hline
\end{tabular}

\medskip

\begin{tabular}{c|ccccc}
\hline
$i$&10&11&12&13&14\\
\hline
$\pi_i(\mm{MO}\langle8\rangle)$
&$\mathbb Z_6$
&$0$
&$\mathbb Z$
&$\mathbb Z_3$
&$\mathbb Z_2$
\\
\hline
\end{tabular}
\end{table}

Combining Proposition~\ref{HBZ} with the Thom isomorphism gives
the following relevant $E^2$-terms.

\begin{table}[!h]
\centering
\caption{Relevant terms of
$E^2_{p,q}
=H_p(\mm{M}\eta;\pi_q(\mm{MO}\langle8\rangle))$.}
\label{E2AHSSforMO8Meta}

\begin{tabular}{c|cccccccc}
\hline
$p+q=13$, $p$
&0&4&5&7&10&11&12&13\\
\hline
$E^2_{p,q}$
&$\mathbb Z_3$
&$\mathbb Z_2^2$
&$\mathbb Z_s$
&$\mathbb Z_2$
&$\mathbb Z_2^2\oplus\mathbb Z_3$
&$\mathbb Z_2^3$
&$\mathbb Z_2^3$
&$\mathbb Z_2^3\oplus
  \mathbb Z_3^2\oplus\mathbb Z_5\oplus\mathbb Z_s$
\\
\hline
\end{tabular}

\medskip

\begin{tabular}{c|cccccccc}
\hline
$p+q=14$, $p$
&0&2&4&8&11&12&13&14\\
\hline
$E^2_{p,q}$
&$\mathbb Z_2$
&$\mathbb Z$
&$\mathbb Z_6$
&$\mathbb Z_2$
&$\mathbb Z_2^3\oplus\mathbb Z_3$
&$\mathbb Z_2^3$
&$\mathbb Z_2^3$
&$\mathbb Z_2$
\\
\hline
\end{tabular}

\medskip

\begin{tabular}{c|cccccc}
\hline
$p+q=12$, $p$
&0&2&4&9&10&11\\
\hline
$E^2_{p,q}$
&$\mathbb Z$
&$\mathbb Z_6$
&$\mathbb Z_2\oplus\mathbb Z$
&$\mathbb Z_2^2\oplus\mathbb Z_3$
&$\mathbb Z_2^2$
&$\mathbb Z_2^3$
\\
\hline
\end{tabular}
\end{table}

\begin{lemma}\label{E_infty_012}
In the AHSS \eqref{AHSSforMO8Meta-odd},
$
E^\infty_{0,12}
=
E^2_{0,12}
\cong\mathbb Z.$
\end{lemma}

\begin{proof}
For
$
E^2_{0,12}
=\mathbb Z,$
 there is no outgoing
differential.

A differential $d_r$ entering $E^r_{0,12}$ has the source 
$
E^r_{r,\,13-r}
$.
Every such source is a subquotient of
$E^2_{r,\,13-r},$
which is finite. It follows inductively that no differential enters the
term $E^r_{0,12}$. This implies the desired equation.
\end{proof}

\begin{lemma}\label{E3103}
In the AHSS \eqref{AHSSforMO8Meta-odd},
$E^3_{10,3}
\cong
\mathbb Z_2\oplus\mathbb Z_3.$
\end{lemma}

\begin{proof}
By Proposition \ref{HBZ} and the universal theorem,
\begin{equation}\label{E2103_Tor}
	E^2_{10,3}\cong \Tor_1^{\Z}(\Z_{24},H_9(\M \eta))\cong \Z_2^2\oplus \Z_3.
\end{equation}
Let $N=H_9(\M \eta)$. The short exact sequence
\[
0\longrightarrow\Z_2
\xrightarrow{\gamma}\Z_{24}
\longrightarrow\Z_{12}
\longrightarrow0
\]
induces an exact sequence
\[
\begin{aligned}
0\longrightarrow&
\Tor_1^\Z(\Z_2,N)_{(2)}
\xrightarrow{\gamma_*}
\Tor_1^\Z(\Z_{24},N)_{(2)}
\\
\longrightarrow&
\Tor_1^\Z(\Z_{12},N)_{(2)}
\xrightarrow{\partial}
(\Z_2\otimes N)_{(2)}
\longrightarrow
(\Z_{24}\otimes N)_{(2)}
\end{aligned}
\]
for the $2$-local groups.
Since $N_{(2)}\cong\Z_2^2$, the last
homomorphism is zero. Thus $\partial$
is surjective. As both its source and
target are isomorphic to $\Z_2^2$,
it is an isomorphism. By exactness,
$\gamma_*$ is therefore an isomorphism.
The naturality of the universal coefficient
theorem now gives
\[
\gamma_*:
H_{10}(\M\eta;\Z_2)
\xrightarrow{\cong}
H_{10}(\M\eta;\Z_{24})_{(2)}
=(E^2_{10,3})_{(2)}.
\]

Let
$a=(\Sq^2\Sq^1u)U$, $b=(x\Sq^1u)U$
be the generators of
$H^{10}(\mathrm M\eta;\mathbb Z_2)$ where $U$ denotes the stable Thom class, and let
$\{a^\vee,b^\vee\}$ denote the dual basis of $H_{10}(\mathrm M\eta;\mathbb Z_2)$. 

By Table \eqref{SqBtable} and $\Sq^2U=w_2(\eta)U$,
$$\Sq^2(a)=\Sq^2(b)=(x\Sq^2\Sq^1u)U.$$
Then Lemma \ref{d2}(3) gives
\[
\begin{aligned}
E^3_{10,3}
&\cong
\frac{
\mathbb Z_2\langle \gamma_\ast(a^\vee) ,\gamma_\ast(b^\vee)\rangle
}{
\mathbb Z_2\langle \gamma_\ast(a^\vee+b^\vee)\rangle
}
\oplus\mathbb Z_3\\
&\cong
\mathbb Z_2\oplus\mathbb Z_3.
\end{aligned}
\]
This finishes the proof.
\end{proof}

\begin{lemma}\label{E_58130_odd}
Let $G$ denote the direct sum of the $p$-primary components of
$
\Omega_{13}^{\mm{O}\langle8\rangle}(\eta)$
over all primes $p\ge5$. Then $G$ fits into the exact sequence
\begin{equation}\label{odd-extension}
0\longrightarrow
\mathbb Z_s
\longrightarrow
G
\longrightarrow
\mathbb Z_s\oplus\mathbb Z_5
\longrightarrow0.
\end{equation}	
\end{lemma}

\begin{proof}
For primes $p\ge5$, Table~\ref{E2AHSSforMO8Meta} shows that the
only contributions to $p$-primary part in total degree $13$ occur in bidegrees
$
(5,8)$ and $(13,0).$
The first contains $\mathbb Z_s$, while the second contains
$\mathbb Z_s\oplus\mathbb Z_5$.

By bidegree considerations in the AHSS
\eqref{AHSSforMO8Meta-odd}, $E_2^{5,8}\cong\mathbb Z_s$
survives to $E_\infty$.  Likewise,
$
\mathbb Z_s\oplus\mathbb Z_5
\subset E_2^{13,0}
$
survives to
$E_\infty^{13,0}$.
Consequently the Atiyah-Hirzebruch (AH) filtration gives the exact sequence.
\end{proof}

\subsection{The Hurewicz homomorphism and the  extension}
\quad 

To understand
\eqref{odd-extension},  we first study the integral Hurewicz homomorphism for the free summand of
$\pi_8(\mm{MO}\langle8\rangle)$.

For a spectrum $X$, the integral Hurewicz homomorphism is induced
by the unit
$
\iota:S^0\longrightarrow H\mathbb Z:$
\begin{equation}\label{hurewicziota}
\iota_*:
\pi_n(X)
\longrightarrow
H_n(X;\mathbb Z).
\end{equation}

\begin{lemma}\label{240MO8Hure}
After compatible choices of generators, the Hurewicz homomorphism
on the free summand is multiplication by $240$:
\[
\iota_*:
\mathbb Z
\subset
\pi_8(\mm{MO}\langle8\rangle)
\longrightarrow
H_8(\mm{MO}\langle8\rangle)
\cong\mathbb Z,
\quad
a\longmapsto240a.
\]
\end{lemma}

\begin{proof}
By Thom isomorphism,
\[
H_8(\mm{MO}\langle8\rangle)
\cong
H_8(\mm{BO}\langle8\rangle)
\cong\mathbb Z,\quad
\text{
and}
\quad
H_i(\mm{MO}\langle8\rangle)=0,
\quad 0<i\le7.
\]
Consider the AHSS
\[
E^2_{p,q}
=
H_p(\mm{MO}\langle8\rangle;\pi_q(S^0))
\Longrightarrow
\pi_{p+q}(\mm{MO}\langle8\rangle).
\]
 The
only possible nontrivial differential for $E^r_{0,7}$ is
\[
d_8:
E^8_{8,0}\cong\mathbb Z
\longrightarrow
E^8_{0,7}
\cong\pi_7^S
\cong\mathbb Z_{240}.
\]
Since
$\pi_7(\mm{MO}\langle8\rangle)=0,$
this differential is surjective.  Hence
\[
E_\infty^{8,0}
=
\ker(d_8)
=
240\mathbb Z
\subset H_8(\mm{MO}\langle8\rangle).
\]
Note that $\chi:S^0\to H\Z$ is homotopic to the unit $\iota$.
By Lemma~\ref{chi}, this edge subgroup is precisely the image of
the Hurewicz homomorphism on the free summand of
$\pi_8(\mm{MO}\langle8\rangle)$.
\end{proof}

By Thom isomorphism and the K\"unneth formula, the prime-to-$6$
part of
$
H_{13}(\mm{M}_{\mathscr B})$
contains
$
\mathbb Z_s\oplus\mathbb Z_s\oplus\mathbb Z_5.$
Choose generators
\[
m_8\in H_8(\mm{MO}\langle8\rangle),
\qquad
y_5\in\mathbb Z_s\subset H_5(\mm{M}\eta),
\]
and
\[
y_{13,s}\in\mathbb Z_s\subset H_{13}(\mm{M}\eta),
\qquad
y_{13,5}\in\mathbb Z_5\subset H_{13}(\mm{M}\eta).
\]
The three corresponding homology generators of $H_{13}(\mm{M}_{\mathscr B})$  are
\[
m_8\otimes y_5,\qquad
y_{13,s},\qquad
y_{13,5}.
\]

\begin{proposition}\label{ZsHure}
Assume that $\gcd(s,5)=1$.  Then the extension
\eqref{odd-extension} splits.  More precisely,
$
G
\cong
\mathbb Z_s\oplus\mathbb Z_s\oplus\mathbb Z_5$.
Furthermore, the Hurewicz homomorphism restricts to an isomorphism
\[
G
\xrightarrow{\cong}
\mathbb Z_s\langle m_8\otimes y_5\rangle
\oplus
\mathbb Z_s\langle y_{13,s}\rangle
\oplus
\mathbb Z_5\langle y_{13,5}\rangle
\subset
H_{13}(\mm{M}_{\mathscr B}).
\]
\end{proposition}

\begin{proof}
Consider the AHSS for ordinary homology
\[
E^2_{p,q}(H)
=
H_p(\mm{M}\eta;
H_q(\mm{MO}\langle8\rangle))
\Longrightarrow
H_{p+q}(\mm{M}_{\mathscr B}).
\]
By \cite[Lemma~4.1]{Giambalvo},
$H_q(\mm{MO}\langle8\rangle)$ has no odd torsion in the relevant
degrees.  Hence
$
m_8\otimes y_5$
is provided by
\[
E_\infty^{5,8}(H)
=
E_2^{5,8}(H)
\cong\mathbb Z_s,
\]
whereas
$
y_{13,s}, y_{13,5}$
are provided by
\[
\mathbb Z_s\oplus\mathbb Z_5
\subset
E_\infty^{13,0}(H).
\]

The Hurewicz homomorphism induces a morphism from the bordism AHSS
to this homology AHSS.  On the two relevant associated-graded
pieces, 
\begin{equation}\label{E_infty_130_ZsZ5_to_homology}
	\text{$\mathbb Z_s\oplus\mathbb Z_5\subset E_\infty^{13,0}
\longrightarrow
\mathbb Z_s\oplus\mathbb Z_5\subset E_\infty^{13,0}(H)$
is the identity}.
\end{equation}

And by Lemma~\ref{240MO8Hure},
$E_\infty^{5,8}
\longrightarrow
E_\infty^{5,8}(H)$ is the
multiplication by 240.
Since
$\gcd(s,6)=1$ and
$\gcd(s,5)=1$,
we have
$
\gcd(s,240)=1.$
Thus multiplication by $240$ is an automorphism of
$\mathbb Z_s$.  

The corresponding homology filtration is split by the K\"unneth
decomposition for $H_\ast(\BO\langle 8\rangle \times B)\cong H_\ast(\MO\langle 8\rangle \wedge \M\eta)$:
\[
\mathbb Z_s\langle m_8\otimes y_5\rangle
\oplus
\bigl(
\mathbb Z_s\langle y_{13,s}\rangle
\oplus
\mathbb Z_5\langle y_{13,5}\rangle
\bigr).
\]
Comparing the two short exact sequences and applying the short
five lemma shows that the Hurewicz homomorphism is an isomorphism
on $G$.  Consequently the extension \eqref{odd-extension} is
split.
\end{proof}

\subsection{Geometric invariants detecting the two
$\mathbb Z_s$-components}\label{Zs}\quad


We record the form of the Thom--Pontryagin principle needed below.

\begin{lemma}\label{detectbyorientableclass}
Let $M^n$ be a closed manifold equipped with a normal
$\mathscr B$-structure
$
\bar{\nu}:M\longrightarrow\mathscr B.
$
Then:
\begin{enumerate}
\item[(1)]
if the Hurewicz homomorphism is injective on a direct summand
$
G\subset\pi_n(\mm{M}_{\mathscr B}),
$
the component of $[M,\bar \nu]$ in $G$ is determined by the image of the fundamental class
$
\bar{\nu}_*[M]
\in H_n(\mathscr B);$

\item[(2)]
every spherical class in
$
H_n(\mm{M}_{\mathscr B})
\cong H_n(\mathscr B)
$
is represented by the image of the fundamental class of a closed
$n$-manifold with normal $\mathscr B$-structure;

\item[(3)]
if
$\bar{\nu}_*[M]=y\in H_n(\mathscr B),$
then the Thom class
$
yU\in H_n(\mm{M}_{\mathscr B})$
is spherical.
\end{enumerate}
\end{lemma}

Recall the normal $5$-smoothing \eqref{eq:normal-smoothing} of a Bazaikin space $\mathcal B_{\bold q}$ 
\[
\bar{\nu}_{{\bold q}}=(g,f):
\mathcal B_{\bold q}\longrightarrow
\mathscr B
=
\mm{BO}\langle8\rangle\times B.
\]
   Under the assumption $(s,5)=1$,
Proposition~\ref{ZsHure} implies that the
\[
\mathbb Z_s\oplus\mathbb Z_s\oplus\mathbb Z_5
\]
component of $[\mathcal B_{\bold q},(g,f)]$ is completely determined by
\[
\bar{\nu}_{{\bold q}*}
[\mathcal B_{\bold q}]
\in H_{13}(\mathscr B).
\]

We now identify the two $\mathbb Z_s$-coordinates.

\subsubsection{The filtration-$13$ summand and the linking form}

Let $M$ have the cohomology ring of a Bazaikin space.  The
Bockstein associated with
\[
0\longrightarrow\mathbb Z
\xrightarrow{\times s}
\mathbb Z
\longrightarrow\mathbb Z_s
\longrightarrow0
\]
is an isomorphism
\[
\beta:
H^5(M;\mathbb Z_s)
\xrightarrow{\cong}
H^6(M),
\]
since $H^5(M)=0$.

The linking form can therefore be written as
\[
L:
H^6(M)\times H^8(M)
\longrightarrow\mathbb Z_s,
\]
\[
L(a,b)
=
\left\langle
\beta^{-1}(a)\cup \rho_s(b),
[M]_s
\right\rangle,
\]
where $\rho_s$ denotes reduction modulo $s$.
For a generator $x\in H^2(M)$, define
\[
\ell k(M)
=
L(x^3,x^4)
\in\mathbb Z_s.
\]

\begin{lemma}\label{130Zs}
The component of
$[\mathcal B_{\bold q},\bar \nu_{\bold q}]$ in the summand
$$
\mathbb Z_s
\subset
\Omega_{13}^{\mm{O}\langle8\rangle}(\eta)/F_{5,8}$$
 is determined by the linking form
$\ell k(\mathcal B_{\bold q})$ where $F_{5,8}$ is the filtration associated with the AHSS \eqref{AHSSforMO8Meta-odd}.
\end{lemma}

\begin{proof}
Since
$
\beta:
H^5(B;\mathbb Z_s)
\xrightarrow{\cong}
H^6(B)$
is an isomorphism, 
$
\beta^{-1}(x^3)$
generates $H^5(B;\mathbb Z_s)$.
Let
\[
\beta_s:
H^i(B;\mathbb Z_s)
\longrightarrow
H^{i+1}(B;\mathbb Z_s)
\]
be the Bockstein associated with
\[
0\longrightarrow\mathbb Z_s
\longrightarrow\mathbb Z_{s^2}
\longrightarrow\mathbb Z_s
\longrightarrow0.
\]
By the standard comparison of Bockstein homomorphisms and the
Cartan formula \cite[pp.~303--304]{Hatcher},
$\beta_s$ restricts to an isomorphism on the relevant
$\mathbb Z_s$-summands.  Thus, for some unit
$m\in\mathbb Z_s^\times$,
\[
\beta_s\bigl(\beta^{-1}(x^3)\bigr)
=
m\,\rho_s(x^3).
\]
Consequently, by $\beta_s
\left( \rho_s(x^4)\right)=0$,
\[
\beta_s
\left(
\beta^{-1}(x^3)\cup \rho_s(x^4)
\right)
=
m\,\rho_s(x^7).
\]
Since $\rho_s(x^7)$ is a generator of $H^{14}(B;\Z_s)$,
\[
\theta_s
:=
\beta^{-1}(x^3)\cup \rho_s(x^4)
\]
is a generator of the relevant
$\mathbb Z_s$-summand of
$H^{13}(B;\mathbb Z_s)$.

Recall that the canonical lift
$
f:\mathcal B_{\bold q}\longrightarrow B$
satisfies
$f^*(x)=x.$
Hence
\[
f^*(\theta_s)
=
\beta^{-1}(x^3)\cup \rho_s(x^4),
\]
and thus
\[
\begin{aligned}
\left\langle
\theta_s,
f_*[\mathcal B_{\bold q}]_s
\right\rangle
&=
\left\langle
f^*(\theta_s),
[\mathcal B_{\bold q}]_s
\right\rangle\\
&=
L(x^3,x^4)
=
\ell k(\mathcal B_{\bold q}).
\end{aligned}
\]
This completes the proof by the identity \eqref{E_infty_130_ZsZ5_to_homology}.
\end{proof}

\begin{remark}
	In Lemma \ref{130Zs}, $s$ may be divided by $5$. 
\end{remark}

\subsubsection{The filtration-$5$ summand and the second
Pontryagin class}

Let
\[
m_8^*\in H^8(\mm{BO}\langle8\rangle;\mathbb Z)
\cong\mathbb Z
\]
denote the dual of $m_8\in H_8(\mm{BO}\langle8\rangle;\mathbb Z)$, and let $\gamma_8$ denote the universal bundle
over $\mm{BO}\langle8\rangle$.

\begin{proposition}\label{secondPgamma8}
There exist integers $a,b\ge0$ such that
\[
p_2(\gamma_8)
=
2^a3^b\,m_8^*.
\]
In particular, for every integer $s$ with $(s,6)=1$, the
coefficient of $m_8^*$ in $p_2(\gamma_8)$ is a unit modulo $s$.
\end{proposition}

\begin{proof}
Consider the standard fibrations
\[
\K(\mathbb Z_2,1)
\longrightarrow
\mm{BSpin}
\longrightarrow
\mm{BSO}
\]
and
\[
\K(\mathbb Z,3)
\longrightarrow
\mm{BO}\langle8\rangle
\longrightarrow
\mm{BSpin}.
\]
The universal bundle over $\mm{BO}\langle8\rangle$ is the
pullback of the universal oriented bundle, so
$
p_2(\gamma_8)
$
is the pullback of
$
p_2\in H^8(\mm{BSO};\mathbb Z).
$

For every prime $p>3$, the Serre spectral sequences of the above
fibrations show that, in total degree $8$, all terms away from the
base row have order prime to $p$.  Therefore, after localization
at $p$, the pullback of $p_2$ generates
\[
H^8(\mm{BO}\langle8\rangle;\mathbb Z_{(p)}).
\]
Hence the coefficient of $m_8^*$ in $p_2(\gamma_8)$ is not
divisible by any prime $p>3$.  It is therefore of the form
$2^a3^b$.
\end{proof}

Recall the normal $5$-smoothing of a Bazaikin space:
\[
\bar{\nu}_{{\bold q}}
=
(g,f):
\mathcal B_{\bold q}
\longrightarrow
\mm{BO}\langle8\rangle\times B.
\]
Its stable normal bundle satisfies
\[
\mathscr N_{\bold q}
\cong
g^*\gamma_8\oplus f^*\eta.
\]
Since
$
p_1(g^*\gamma_8)=0
$ and $\gcd(s,2)=1$,
we obtain
\[
p_2(\mathscr N_{\bold q})
=
g^*p_2(\gamma_8)
+
f^*p_2(\eta)\in H^8(\mathcal B_{\bold q})=\mb{Z}_s.
\]
The second summand is fixed by the map $f$.  Moreover,
$p_1(\mathcal B_{\bold q})$ is fixed in our comparison, and
\[
p(\mathscr N_{\bold q})
\,p(T\mathcal B_{\bold q})=1.
\]
Hence $p_2(\mathcal B_{\bold q})$ determines
$p_2(\mathscr N_{\bold q})$ and thus determines
\[
g^*p_2(\gamma_8)
\in H^8(\mathcal B_{\bold q})
\cong\mathbb Z_s.
\]
By Proposition~\ref{secondPgamma8} and $(s,6)=1$, this in turn
determines
\[
g^*(m_8^*)
\in H^8(\mathcal B_{\bold q}).
\]

\begin{lemma}
	Assume $\gcd(\sigma_3(\bold q),5)=1$. The $\mathbb Z_s\oplus\mathbb Z_s$-component of $[\mathcal B_\bold q,\bar{\nu}_\bold q]$ in $\Omega_{13}^{\mathrm{O}\langle 8\rangle}(\eta)$ is
determined by
$\ell k(\mathcal B_{\bold q})$
and
$p_2(\mathcal B_{\bold q}).$
\end{lemma}

\begin{proof}
By Proposition~\ref{ZsHure}, the two $\mathbb Z_s$-summands in
$\Omega_{13}^{\mm{O}\langle8\rangle}(\eta)$
are provided respectively by
$
E_\infty^{13,0}$
and
$E_\infty^{5,8}$.
 The Hurewicz homomorphism identifies them with
$
\mathbb Z_s\langle y_{13,s}\rangle
$ and
$\mathbb Z_s\langle m_8\otimes y_5\rangle$
in $H_{13}(\mathscr B)$.

By Lemma~\ref{130Zs}, the first coordinate is determined by
$\ell k(\mathcal B_{\bold q}).$

It remains to identify the filtration-$5$ coordinate.
Write
\[
g^*(m_8^*)
=
nx^4
\in H^8(\mathcal B_{\bold q})
\cong\mathbb Z_s.
\]
As observed above, $n$ modulo $s$ is
determined by $p_2(\mathcal B_{\bold q})$.

The class
\[
m_8^*
\cup
\beta^{-1}(x^3)
\in
H^{13}(\mathscr B;\mathbb Z_s)
\]
is dual, up to multiplication by a unit of $\mathbb Z_s$, to
the generator
\[
m_8\otimes y_5
\in
H_{13}(\mathscr B;\mathbb Z_s).
\]
Pulling it back along the normal $5$-smoothing gives
\[
g^*(m_8^*)
\cup
\beta^{-1}(x^3)
=
nx^4\cup\beta^{-1}(x^3).
\]
Therefore
\begin{equation}\label{p2deterZs58}
\begin{aligned}
&
\left\langle
m_8^*\cup\beta^{-1}(x^3),
\bar{\nu}_{{\bold q}*}
[\mathcal B_{\bold q}]_s
\right\rangle
\\
&\qquad=
n
\left\langle
\beta^{-1}(x^3)\cup \rho_s(x^4),
[\mathcal B_{\bold q}]_s
\right\rangle
\\
&\qquad=
n\,\ell k(\mathcal B_{\bold q}).
\end{aligned}
\end{equation}
Hence the filtration-$5$ $\mathbb Z_s$-coordinate is determined
by the linking form together with
$p_2(\mathcal B_{\bold q})$. This completes the proof.
\end{proof}

\section{The $2$-primary part of $\Omega_{13}^{\mm{O}\langle 8 \rangle}(\eta)$}\label{2primary}

Recall $\Omega_{\ast}^{\mm{O}\langle 8 \rangle}(\eta)\cong \pi_\ast(\mm M_{\mathscr B})$ where $M_{\mathscr B}=\mm{MO}\langle 8\rangle\wedge \mm M\eta$.

 In this section, we  
 use three spectral sequences:
 The Adams spectral sequence (ASS) computes the $2$-primary homotopy $\Omega_{13}^{\mm{O}\langle 8 \rangle}(\eta)_{(2)}$.  The algebraic Atiyah--Hirzebruch spectral sequence (AAHSS) is only a computational device for the $E_2$-term of the ASS.  Finally, the Atiyah--Hirzebruch spectral sequence (AHSS) computes the same homotopy groups $\Omega_{13}^{\mm{O}\langle 8 \rangle}(\eta)_{(2)}$ from the skeletal filtration of $\mm M\eta$.  

\subsection{The Adams $E_2$-term and the algebraic AHSS}\quad

The mod-$2$ Adams spectral sequence is
\begin{equation}\label{ASSMB}
E_2^{s,t}
=\mm{Ext}_{\mathscr A}^{s,t}
\bigl(H^*(\mm M_{\mathscr B};\mathbb Z_2),\mathbb Z_2\bigr)
\Longrightarrow {_2}\pi_{t-s}(\mm M_{\mathscr B}),
\end{equation}
$$
d_r:E_r^{s,t}\longrightarrow E_r^{s+r,t+r-1}.$$
Here $\mathscr A$ is the mod-$2$ Steenrod algebra. In particular, $E_2^{0,t}$ is a subgroup of 
$$H_t(\M_{\mathscr B};\Z_2)\cong \mm{Hom}_{\Z_2}(H^t(\M_{\mathscr B};\Z_2),\Z_2).$$

 By the standard edge description of the Adams spectral sequence (see \cite[Chapter~19]{Sw}), we have
\begin{lemma}\label{E20tASS}
A class $e\in E_2^{0,t}$ survives to $E_\infty$ in the ASS if and only if
$
e\in H_t(\mm M_{\mathscr B};\mathbb Z_2)$
is spherical. 
\end{lemma}

By the Kunneth theorem,
\[
H^*(\mm M_{\mathscr B};\mathbb Z_2)
\cong
H^*(\mm M\eta;\mathbb Z_2)\otimes
H^*(\mm{MO}\langle8\rangle;\mathbb Z_2).
\]

On the other hand, in dimensions $<16$, Giambalvo's calculation gives
\[
H^*(\mm{MO}\langle8\rangle;\mathbb Z_2)
\cong \mathscr A\mathbin{/\!/}\mathscr A_2
\]
as an $\mathscr A$-module, where $\mathscr A_2$ is the Hopf subalgebra generated by $\Sq^1,\Sq^2,$ and $\Sq^4$ \cite{Giambalvo}.  Hence the change-of-rings theorem \cite{Ravenel1986} gives, for $t-s<15$,
\begin{equation}\label{changeofringsA2}
\mm{Ext}_{\mathscr A}^{s,t}
\bigl(H^*(\mm M_{\mathscr B};\mathbb Z_2),\mathbb Z_2\bigr)
\cong
\mm{Ext}_{\mathscr A_2}^{s,t}
\bigl(H^*(\mm M\eta;\mathbb Z_2),\mathbb Z_2\bigr).
\end{equation}

 Set
$
N=H^*(\mm M\eta;\mathbb Z_2).$
Let $U$ denotes the stable Thom class, then
\[
N\cong H^*(B;\mathbb Z_2)\otimes U
\]
as a graded vector space. The Steenrod operations on $H^*(\M\eta;\mathbb Z_2)$ are determined by 
 Tables~\eqref{HBZ2} and~\eqref{SqBtable}, together with $\Sq^iU=w_i(\eta)U$ where
\[
\text{$w_2(\eta)\ne0$ and
$w_4(\eta)=0$ (see \eqref{eq:eta-char}).}
\]

From \cite{Adams1960}, we have the
 standard Ext--Tor duality:
 \[
\mm{Ext}_{\mathscr A_2}^{s,t}(N,\Z_2)
\cong
\mm{Tor}_{s,t}^{\mathscr A_2}(\Z_2,N).
\]
Hence we first compute the latter.

Write
\[
N_n=\bigoplus_{i\ge n}N^i,
\quad
N^i=H^i(\mm M\eta;\mathbb Z_2).
\]
Filtering the bar complex by the submodules $N_n$ gives an exact couple and hence an algebraic Atiyah--Hirzebruch spectral sequence
\begin{equation}\label{AAHSS}
E_1^{s,m,n}
\cong
\mm{Tor}_{s,s+m}^{\mathscr A_2}(\mathbb Z_2,\mathbb Z_2)
\otimes N^n
\Longrightarrow
\mm{Tor}_{s,s+m+n}^{\mathscr A_2}(\mathbb Z_2,N),
\end{equation}
with
$
d_r:E_r^{s,m,n}
\longrightarrow
E_r^{s-1,m+1-r,n+r}.$

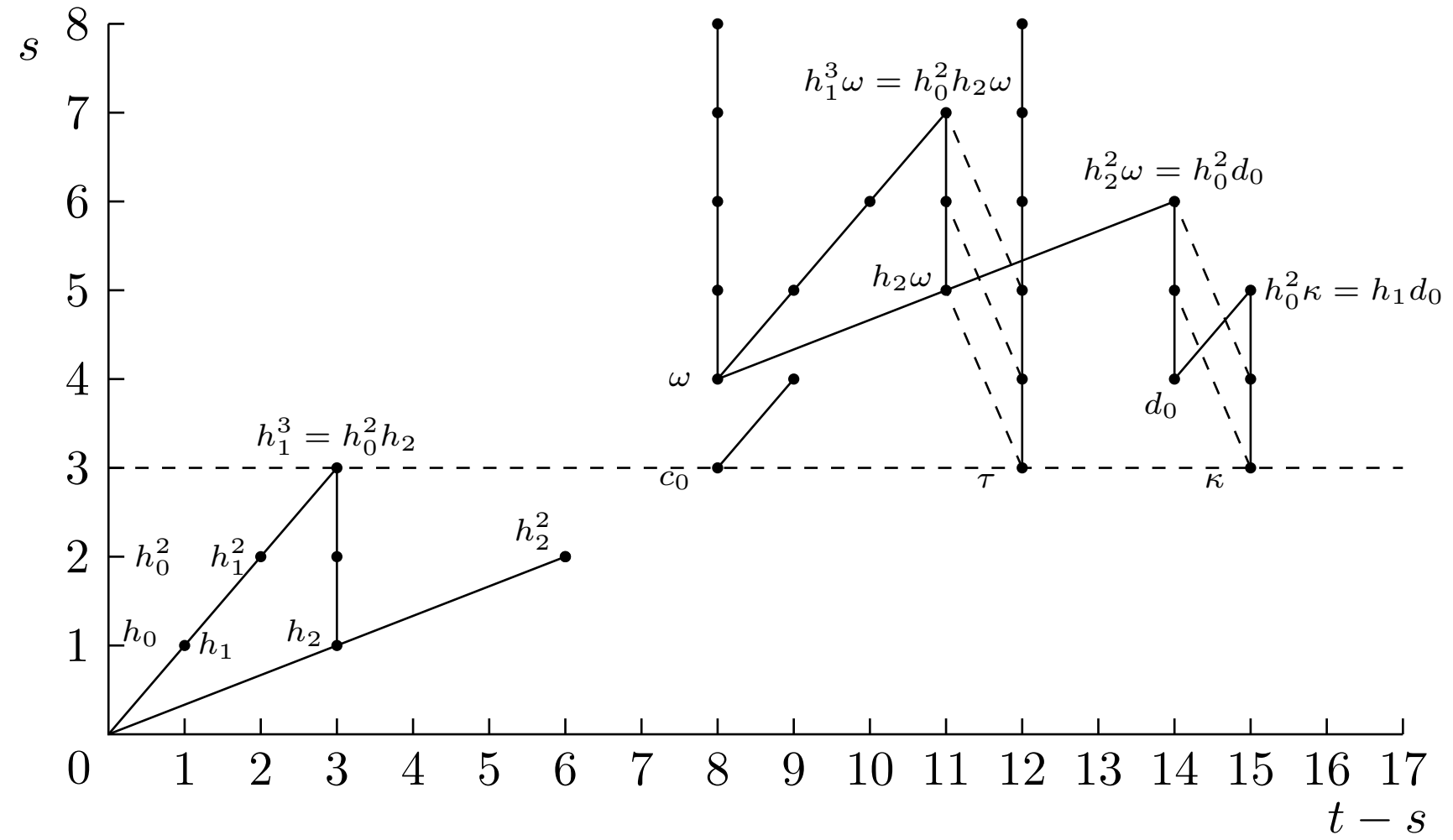
\begin{figure}[htb]
\centering
\begin{tikzpicture}[>=stealth,scale=1,line width=0.5pt]
\pgfmathsetmacro{\ticker}{0.125}

\coordinate [label=225:$0$](A) at (0,0);
\coordinate (B) at (0,5.6);
\coordinate (D) at (10.2,0);
\draw(D)--(A)--(B);

\coordinate [label=left:$s$](E) at ($(B)+(-0.4,-0.2)$);
\coordinate [label=below:$t-s$](F) at ($(D)+(-0.2,-0.4)$);

\foreach \i/\texti  in {1,2,3,4,5,6,7,8,9,10,11,12,13,14,15,16,17} {
\draw (0.6*\i,0) --(0.6*\i,\ticker) node[label=below:\texti]{};
}
\foreach \j/\textj  in {1,2,3,4,5,6,7,8} {
\draw (0,0.7*\j) --(\ticker,0.7*\j) node[label=left:\textj]{};
}
{\tiny
\coordinate[label=right:$h_0$] (I) at (0,0.8);
}
{\tiny
\coordinate[label=right:$h^2_0$] (I) at (0.1,1.4);
}
\path[name path = row1](0,0.7)--(10.2,0.7);
\path[name path = row2](0,1.4)--(10.2,1.4);
\draw[name path = row3][dashed] (0,2.1) -- (10.2,2.1);
\draw[name path = xieh1](0,0)--(0.6,0.7)--(1.2,1.4)--(1.8,2.1);
\path[draw,fill,name intersections={of = row1 and xieh1,by=h1}](h1)circle(1pt);
{\tiny
\coordinate[label=right:$h_1$] (I) at (0.6,0.7);
}
\path[draw,fill,name intersections={of = row2 and xieh1,by=h12}](h12)circle(1pt);
{\tiny
\coordinate[label=left:$h^2_1$] (I) at (1.2,1.4);
}
\path[draw,fill,name intersections={of = row3 and xieh1,by=h13}](h13)circle(1pt);
{\tiny
\coordinate[label=above:${h^3_1=h^2_0h_2}$] (I) at (1.8,2.1);
}
\draw[name path = c3] (1.8,0.7) -- (h13);
\path[draw,fill,name intersections={of = c3 and row2,by=h0h2}](h0h2)circle(1pt);
\path[draw,fill,name intersections={of = c3 and row1,by=h2}](h2)circle(1pt);
{\tiny
\coordinate[label=left:${h_2}$] (I) at (1.8,0.8);
}
\draw[name path = xieh2](0,0)--(h2)--(3.6,1.4);
\path[draw,fill,name intersections={of = xieh2 and row2,by=h22}](h22)circle(1pt);
{\tiny
\coordinate[label=left:${h^2_2}$] (I) at (3.6,1.6);
}
\path[name path = row4](0,2.8)--(10.2,2.8);
\path[name path = row5](0,3.5)--(10.2,3.5);
\path[name path = row6](0,4.2)--(10.2,4.2);
\path[name path = row7](0,4.9)--(10.2,4.9);
\path[name path = row8](0,5.6)--(10.2,5.6);
\path[name path = c8](4.8,0)--(4.8,5.6);
\path[draw,fill,name intersections={of = c8 and row3,by=c0}](c0)circle(1pt);
{\tiny
\coordinate[label=left:${c_0}$] (I) at (4.7,2.0);
}
\path[draw,fill,name intersections={of = xieh2 and row2,by=h22}](h22)circle(1pt);
\draw(4.8,2.8)--(4.8,5.6);
\path[draw,fill,name intersections={of = c8 and row4,by=omega}](omega)circle(1pt);
{\tiny
\coordinate[label=left:${\omega}$] (I) at (4.7,2.8);
}
\path[draw,fill,name intersections={of = c8 and row5,by=h0omega}](h0omega)circle(1pt);
\path[draw,fill,name intersections={of = c8 and row6,by=h02omega}](h02omega)circle(1pt);
\path[draw,fill,name intersections={of = c8 and row7,by=h03omega}](h03omega)circle(1pt);
\path[draw,fill,name intersections={of = c8 and row8,by=h04omega}](h04omega)circle(1pt);
\draw[name path = xiec0](4.8,2.1)--(5.4,2.8);
\path[draw,fill,name intersections={of = row4 and xiec0,by=h1c0}](h1c0)circle(1pt);
\draw[name path = xieomega](4.8,2.8)--(6.6,4.9);
\path[draw,fill,name intersections={of = row5 and xieomega,by=h1omega}](h1omega)circle(1pt);
\path[draw,fill,name intersections={of = row6 and xieomega,by=h12omega}](h12omega)circle(1pt);
\path[draw,fill,name intersections={of = row7 and xieomega,by=h13omega}](h13omega)circle(1pt);
{\tiny
\coordinate[label=above:${h_1^3\omega=h^2_0h_2\omega}$] (I) at (6.3,4.9);
}
\draw[name path = c11](h13omega)--(6.6,3.5);
\path[draw,fill,name intersections={of = row5 and c11,by=h2omega}](h2omega)circle(1pt);
\path[draw,fill,name intersections={of = row6 and c11,by=h0h2omega}](h0h2omega)circle(1pt);
{\tiny
\coordinate[label=left:${h_2\omega}$] (I) at (6.6,3.6);
}
\draw[name path = xie2omega](omega)--(8.4,4.2);
\path[draw,fill,name intersections={of = row6 and xie2omega,by=h22omega}](h22omega)circle(1pt);
{\tiny
\coordinate[label=above:${h^2_2\omega=h^2_0d_0}$] (I) at (h22omega);
}
\draw[name path = c12](7.2,2.1)--(7.2,5.6);
\path[draw,fill,name intersections={of = row3 and c12,by=tau}](tau)circle(1pt);
\path[draw,fill,name intersections={of = row4 and c12,by=h0tau}](h0tau)circle(1pt);
\path[draw,fill,name intersections={of = row5 and c12,by=h02tau}](h02tau)circle(1pt);
\path[draw,fill,name intersections={of = row6 and c12,by=h03tau}](h03tau)circle(1pt);
\path[draw,fill,name intersections={of = row7 and c12,by=h04tau}](h04tau)circle(1pt);
\path[draw,fill,name intersections={of = row8 and c12,by=h05tau}](h05tau)circle(1pt);
{\tiny
\coordinate[label=left:${\tau}$] (I) at (7.1,2);
}
\draw[name path = c14](h22omega)--(8.4,2.8);
\path[draw,fill,name intersections={of = row4 and c14,by=d0}](d0)circle(1pt);
\path[draw,fill,name intersections={of = row5 and c14,by=h0d0}](h0d0)circle(1pt);
{\tiny
\coordinate[label=below:${d_0}$] (I) at (8.3,2.8);
}
\draw[name path = xied0](d0)--(9,3.5);
{\tiny
\coordinate[label=right:${h^2_0\kappa=h_1d_0}$] (I) at (9,3.5);
}
\path[draw,fill,name intersections={of = row5 and xied0,by=h1d0}](h1d0)circle(1pt);
\draw[name path = c15](h1d0)--(9,2.1);
\path[draw,fill,name intersections={of = row3 and c15,by=kappa}](kappa)circle(1pt);
\path[draw,fill,name intersections={of = row4 and c15,by=h0kappa}](h0kappa)circle(1pt);
{\tiny
\coordinate[label=left:${\kappa}$] (I) at (8.9,2);
}
\draw[dashed](7.2,2.1)--(6.6,3.5);
\draw[dashed](7.2,2.8)--(6.6,4.2);\draw[dashed](7.2,3.5)--(6.6,4.9);
\draw[dashed](9,2.1)--(8.4,3.5);\draw[dashed](9,2.8)--(8.4,4.2);
\end{tikzpicture}
\renewcommand{\figurename}{Fig}
 \caption{$\mm{Ext}^{s,t}_{\mathscr{A}_2}(\mathbb{Z}_2, \mathbb{Z}_2)$}
  \label{1Fig}
\end{figure}

 Fig.~\ref{1Fig} records the relevant part of
$\mm{Tor}_{s,t}^{\mathscr A_2}(\mathbb Z_2,\mathbb Z_2)\cong \mm{Ext}^{s,t}_{\mathscr{A}_2}(\mathbb{Z}_2, \mathbb{Z}_2)$ (cf. \cite{Giambalvo}).  A vertical line denotes multiplication by $h_0$, a line of slope $1$ multiplication by $h_1$, and a line of slope $1/3$ multiplication by $h_2$. In the bar complex for $\mm{Tor}^{\mathscr A_2}_{*,*}(\mathbb Z_2,\mathbb Z_2)$, $h_0,h_1,h_2$ denote $\Sq^1,\Sq^2,\Sq^4$, respectively.  

\begin{figure}[tb]
\centering
\begin{tikzpicture}[>=stealth,scale=1,line width=0.5pt]

\pgfmathsetmacro{\ticker}{0.125}

\coordinate [label=225:$0$](A) at (0,0);
\coordinate (B) at (0,5.6);
\coordinate (D) at (10.2,0);
\draw(D)--(A)--(B);

\coordinate [label=left:$s$](E) at ($(B)+(-0.4,-0.2)$);
\coordinate [label=below:$t-s$](F) at ($(D)+(-0.2,-0.4)$);

\foreach \j/\textj  in {1,2,3,4,5,6,7,8} {
\draw (0,0.7*\j) --(\ticker,0.7*\j) node[label=left:\textj]{};
}
\path[name path = x](D)--(A);
\path[name path = row1](0,0.7)--(10.2,0.7);
\path[name path = row2](0,1.4)--(10.2,1.4);
\draw[name path = row3][dashed] (0,2.1) -- (10.2,2.1);
\path[name path = row4](0,2.8)--(10.2,2.8);
\path[name path = row5](0,3.5)--(10.2,3.5);
\path[name path = row6](0,4.2)--(10.2,4.2);
\draw[name path = row7][dashed](0,4.9)--(10.2,4.9);
\path[name path = row8](0,5.6)--(10.2,5.6);
\path[name path = c1] (0.6,0) -- (0.6,5.6);
\draw[name path = c11] (0.6,0) -- (0.6,0.15);
\path[fill,name intersections={of = c11 and x,by=12}];
{\tiny
\coordinate[label=below:$12$] (I) at (0.6,-0.1);
}
\path[name path = c2] (3.8,0) -- (3.8,5.6);
\draw[name path = c21] (3.8,0) -- (3.8,0.15);
\path[fill,name intersections={of = c21 and x,by=13}](13)circle(1pt);
{\tiny
\coordinate[label=below:$13$] (I) at (3.8,-0.1);
}
\path[name path = c3] (7,0) -- (7,5.6);
\draw[name path = c31] (7,0) -- (7,0.15);
\path[fill,name intersections={of = c31 and x,by=14}];
{\tiny
\coordinate[label=below:$14$] (I) at (7,-0.1);
}
\path[name path = c4](7.2,0)--(7.2,5.6);
\path[fill,name intersections={of = c4 and row4,by=d0U}](d0U)circle(1pt);
{\tiny
\coordinate[label=below:$d_0U$] (I) at (7.5,2.9);
}
\path[name path = c7](6.9,0)--(6.9,5.6);
\path[fill,name intersections={of = c7 and row6,by=6104}](6104)circle(1pt);
{\tiny
\coordinate[label=left:$h_1^2\omega x^2 U$] (I) at (6104);
}
\path[fill,name intersections={of = c4 and row5,by=h0d0U}](h0d0U)circle(1pt);
\path[fill,name intersections={of = c4 and row6,by=h02d0U}](h02d0U)circle(1pt);
{\tiny
\coordinate[label=right:${h_2^2\omega U=h_0^2d_0U}$] (I) at (h02d0U);
}
\draw (d0U) -- (h02d0U);

\path[fill,name intersections={of = c3 and row3,by=tauxU}](tauxU)circle(1pt);
{\tiny
\coordinate[label=below:${\tau xU}$] (I) at (tauxU);
}
\draw (tauxU) -- (7,5.6);
\path[name path = c5](4,0)--(4,5.6);
\path[fill,name intersections={of = c5 and row5,by=5112}](5112)circle(1pt);
{\tiny
\coordinate[label=right:${h_2\omega xU}$] (I) at (5112);
}
\draw[dashed] (tauxU) -- (5112);
\path[fill,name intersections={of = c3 and row4,by=1tauxU}](1tauxU)circle(1pt);
\path[fill,name intersections={of = c3 and row5,by=2tauxU}](2tauxU)circle(1pt);
\path[fill,name intersections={of = c3 and row6,by=3tauxU}](3tauxU)circle(1pt);
\path[fill,name intersections={of = c3 and row7,by=4tauxU}](4tauxU)circle(1pt);
\path[fill,name intersections={of = c3 and row8,by=5tauxU}](5tauxU)circle(1pt);
{\tiny
\coordinate[label=right:${h_0^5\tau xU}$] (I) at (5tauxU);
}
\path[fill,name intersections={of = c5 and row6,by=6112}](6112)circle(1pt);
\draw[dashed] (1tauxU) -- (6112);

\path[fill,name intersections={of = c3 and row1,by=1311}](1311)circle(1pt);
{\tiny
\coordinate[label=left:${h_2vxU}$] (I) at (1311);
}
\path[fill,name intersections={of = c4 and row1,by=1113}](1113)circle(1pt);
{\tiny
\coordinate[label=right:${h_1vx^2U}$] (I) at (1113);
}

\path[fill,name intersections={of = c5 and row5,by=5112}](5112)circle(1pt);

\path[fill,name intersections={of = c5 and row6,by=6112}](6112)circle(1pt);
{\tiny
\coordinate[label=right:${h_0h_2\omega xU}$] (I) at (6112);
}
\path[fill,name intersections={of = c5 and row7,by=7112}](7112)circle(1pt);
{\tiny
\coordinate[label=above:${(h_0^2h_2=h_1^3)\omega xU}$] (I) at (7112);
}
\draw(5112)--(7112);

\path[fill,name intersections={of = c5 and row4,by=494}](494)circle(1pt);
{\tiny
\coordinate[label=right:${h_1c_0 x^2U}$] (I) at (494);
}
\path[fill,name intersections={of = c2 and row5,by=594}](594)circle(1pt);
{\tiny
\coordinate[label=left:${h_1\omega x^2U}$] (I) at (594);
}

{\tiny
\coordinate[label=left:${v x^2U}$] (I) at (3.8, 0.2);
}
\path[fill,name intersections={of = c5 and row2,by=2211}](2211)circle(1pt);
{\tiny
\coordinate[label=right:${ h_1^2vxU}$] (I) at (2211);
}
\path[fill,name intersections={of = c1 and row3,by=0tauU}](0tauU)circle(1pt);
\path[fill,name intersections={of = c1 and row4,by=1tauU}](1tauU)circle(1pt);
\path[fill,name intersections={of = c1 and row5,by=2tauU}](2tauU)circle(1pt);
\path[fill,name intersections={of = c1 and row6,by=3tauU}](3tauU)circle(1pt);
\path[fill,name intersections={of = c1 and row7,by=4tauU}](4tauU)circle(1pt);
\path[fill,name intersections={of = c1 and row8,by=5tauU}](5tauU)circle(1pt);
\draw(0tauU)--(5tauU);
{\tiny
\coordinate[label=below:${\tau U}$] (I) at (0tauU);
}

\path[name path = c6](0.8,0)--(0.8,5.6);
\path[fill,name intersections={of = c6 and row4,by=484}](484)circle(1pt);
\path[fill,name intersections={of = c6 and row5,by=584}](584)circle(1pt);
\path[fill,name intersections={of = c6 and row6,by=684}](684)circle(1pt);
\path[fill,name intersections={of = c6 and row7,by=784}](784)circle(1pt);
\path[fill,name intersections={of = c6 and row8,by=884}](884)circle(1pt);
\draw(484)--(884);
{\tiny
\coordinate[label=right:${\omega x^2 U}$] (I) at (484);
}
{\tiny
\coordinate[label=right:${h^4_0\omega x^2 U}$] (I) at (884);
}
\path[fill,name intersections={of = c6 and row3,by=384}](384)circle(1pt);
{\tiny
\coordinate[label=below:${c_0x^2 U}$] (I) at (1.2,2.5);
}
\path[fill,name intersections={of = c2 and row1,by=1112}](1112)circle(1pt);
{\tiny
\coordinate[label=left:${h_1x^2\mathrm{Sq}^1u U}$] (I) at (1112);
}
\path[fill,name intersections={of = c2 and row2,by=c2r2}](c2r2)circle(1pt);
{\tiny
\coordinate[label=left:${h_1^2ux^2 U}$] (I) at (c2r2);
}
\draw (3.8,0) -- (3.8,1.4);

\path[fill,name intersections={of = c1 and row1,by=139}](139)circle(1pt);
{\tiny
\coordinate[label=below:${h_2(v+ux)U}$] (I) at (1.1,0.7);
}
\path[fill,name intersections={of = c5 and row3,by=3310}](3310)circle(1pt);
{\tiny
\coordinate[label=above:${h_0^2h_2(\mathrm{Sq}^2\mathrm{Sq}^1u+x\mathrm{Sq}^1u)U}$] (I) at (4,2.1);
}
\draw(4,1.4)--(4,2.1);
\draw[dashed] (2tauxU)--(7112);
\end{tikzpicture}
\renewcommand{\figurename}{Fig}
 \caption{$\mm{Ext}^{s,t}_{\mathscr{A}_2}(H^\ast(\mm{M}\eta;\mb{Z}_2), \mathbb{Z}_2)$}
  \label{2Fig}
\end{figure}
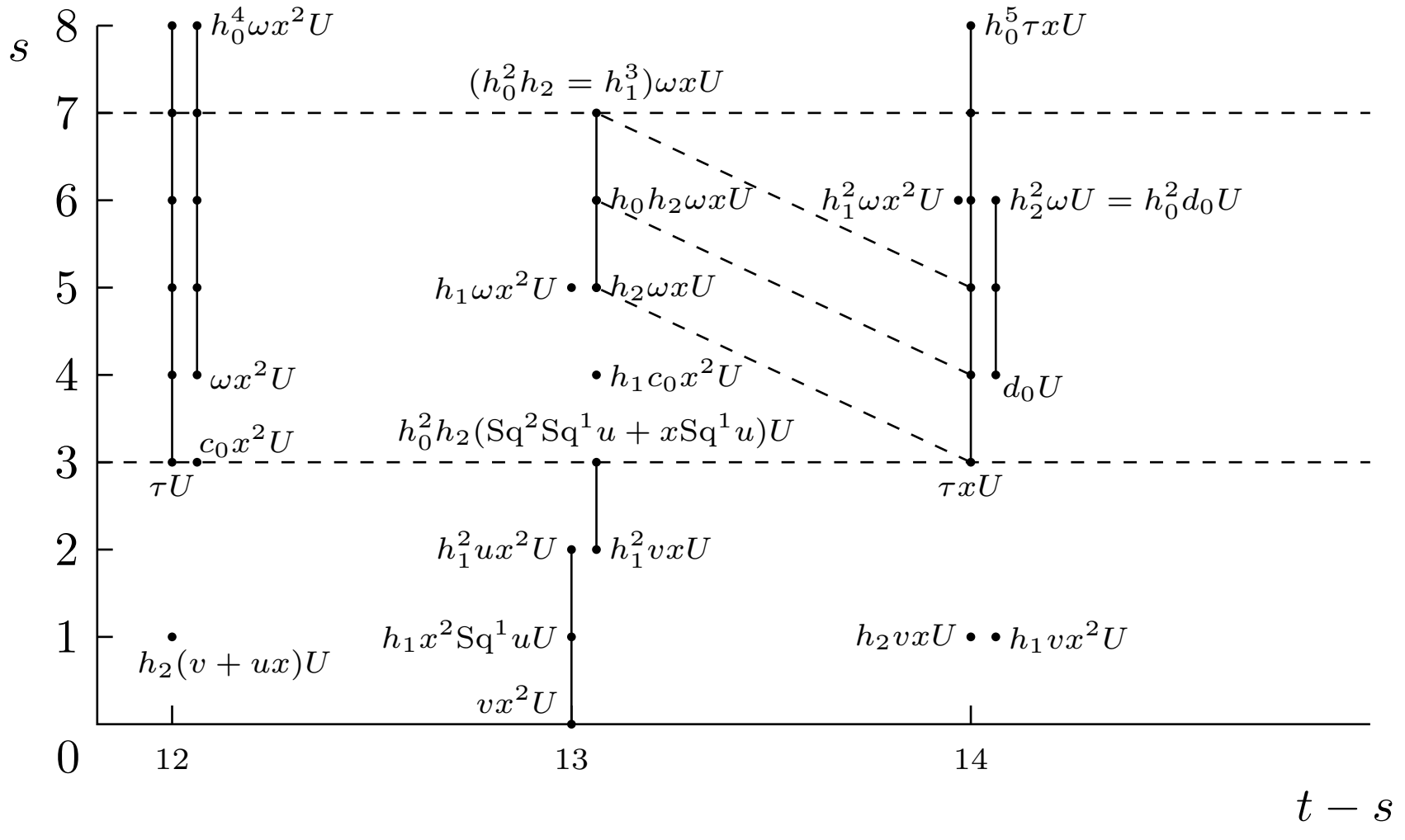

The details of the computations of the AAHSS are shown in Appendix \ref{AAHSSmod2}. Here, we give the final result in Fig.\ref{2Fig} where we use the elements of $E_1$-terms of the AAHSS to represent the generators of the $\mm{Ext}$-groups.

The vertical $h_0$-lines in Fig.~\ref{2Fig} are particularly important.  In the Adams spectral sequence, $h_0$ detects multiplication by $2$.  In other words, if $x,h_0x,$ $\ldots,h_0^rx$ survive to $E_\infty$ and $\theta$ is detected by $x$, then $2^j\theta$ is detected by $h_0^jx$ for $0\le j\le r$.  Thus surviving $h_0$-towers encode precisely the additive $2$-extensions.

The following relations computed by the bar complex control two families of nontrivial $2$-power extensions in degree $13$.

\begin{lemma}\label{h0module}
In $\mm{Ext}^{*,*}_{\mathscr A_2}(N,\mathbb Z_2)$ one has
\begin{equation}\label{h0relations}
\begin{aligned}
h_0vx^2U
&=h_1x^2\Sq^1uU,\\
h_0^2vx^2U
&=h_1^2ux^2U,\\
h_0h_1^2vxU
&=h_1^3(\Sq^2\Sq^1u+x\Sq^1u)U.
\end{aligned}
\end{equation}
\end{lemma}

\begin{proof}
By the $\mathscr A_2$-module structure of $N$,
\begin{equation}\label{vx2notin}
vx^2U\notin \Sq^iH^{13-i}(\mm M\eta;\mathbb Z_2),
\qquad i=1,2,4,
\end{equation}
so $vx^2U$ represents a nonzero class in $\mm{Ext}^{0,13}_{\mathscr A_2}(N,\mathbb Z_2)$.

The class $h_1x^2\Sq^1uU$ has the bar-cycle representative
\[
[\Sq^2]x^2\Sq^1uU+[\Sq^1]vx^2U.
\]
   The second summand is precisely the $h_0$-multiple of the class $vx^2U$, which yields the first equality in \eqref{h0relations}.  Likewise, the degree-two bar cycle
\[
\begin{aligned}
&[\Sq^2|\Sq^2]ux^2U+[\Sq^3|\Sq^1]ux^2U\\
&\qquad+[\Sq^1|\Sq^2]x^2\Sq^1uU+[\Sq^1|\Sq^1]vx^2U
\end{aligned}
\]
represents $h_1^2ux^2U$ and gives the second equality.

For the last equality, put
\[
\begin{aligned}
a&=(v+ux)U,\\
b&=(\Sq^2\Sq^1u+x\Sq^1u)U,\\
c&=vxU.
\end{aligned}
\]
The $\mathscr A_2$-module structure of $N$
and the Cartan formula give
\[
\Sq^1a=b,\quad
\Sq^2a=c.
\]

Consider the cobar complex for $N^\vee=H_\ast(\M\eta;\Z_2)$. Let $\delta$ denote the co-edge operation. Then
$$\delta(h_1^2c^\vee)=h_1^3a^\vee, \quad \delta(h_0h_2b^\vee)=h_0h_2h_0a^\vee.$$
Since $h_1^3=h_0^2h_2=h_0h_2h_0\in \Ext^{3,6}_{\mathscr A_2}(\Z_2\,\Z_2)$, there exists a co-chain $T$ so that
$$\delta(T)=h_1^3+h_0h_2h_0.$$
On the other hand, $a^\vee$ is a primitive comodule element. Hence,
$$\delta(T|a^\vee)=h_1^3a^\vee+h_0h_2h_0a^\vee.$$
And thus we obtain a co-cycle
$$h_1^2c^\vee+h_0h_2b^\vee+T|a^\vee.$$
By the leading term, this co-cycle  represents the class $h_1^2vxU\in \Ext^{2,15}_{\mathscr A_2}(N,\Z_2)$.

 After multiplication $h_0$, we obtain a new co-cycle
$$h_0h_1^2c^\vee+h^2_0h_2b^\vee+h_0T|a^\vee.$$
Since $h_0h_1^2$ is a co-boudary, this co-cycle has the leading term $h^2_0h_2b^\vee$, and represents the class $h_0^2h_2b\in \Ext^{3,16}_{\mathscr A_2}(N,\Z_2)$. This finishes the proof.
\end{proof}

\begin{lemma}\label{h0vx2survive}
The classes
\[
vx^2U,
\quad h_1x^2\Sq^1uU,
\quad h_1^2ux^2U,
\quad h_1^2vxU,
\quad h_1vx^2U
\]
survive to $E_\infty$ in the ASS for $\pi_*(\mm M_{\mathscr B})$.
\end{lemma}

\begin{proof}
By Proposition~\ref{image} and Lemma \ref{detectbyorientableclass}, there is an element
$
\Pi\in\pi_{13}(\mm M_{\mathscr B})$
satisfying that
\begin{enumerate}
	\item $\Pi$ represents a component of the normal $5$-smoothing of the $\mathcal B_\bold q$;
	\item The mod-$2$ Hurewicz image of $\Pi$ is $(vx^2U)^\vee$.
\end{enumerate}
Hence $vx^2U$ is a permanent cycle by Lemma~\ref{E20tASS}.  Since $h_0$ and $h_1$ are permanent cycles and the ASS is multiplicative, Lemma~\ref{h0module} implies that
all outgoing Adams differentials for these classes
\[
h_1x^2\Sq^1uU=h_0vx^2U,
\qquad
h_1^2ux^2U=h_0^2vx^2U,
\qquad
h_1vx^2U
\]
are zero. By degree reasons and $\mm{Ext}^{0,14}(N)=0$, they cannot be targets of any Adams differentials. Hence, they
survives to $E_\infty$.

  Finally, by degree reason and $\mm{Ext}^{0,14}(N)=0$, Fig.~\ref{2Fig} shows that all Adams differentials for $h_1^2vxU$ are zero.  Thus it also survives to $E_\infty$.
\end{proof}

\begin{remark}
As shown in Fig.~\ref{2Fig}, all possible outgoing Adams differentials from the class
\[
h_0^2h_2(\Sq^2\Sq^1u+x\Sq^1u)U
\]
vanish. However, we cannot directly rule out such an incoming differential
\[
d_2(h_2vxU)
=
h_0^2h_2(\Sq^2\Sq^1u+x\Sq^1u)U.
\]
We will show, using further arguments involving the AHSS and \(ko\)-homology, that this class survives to \(E_\infty\) in the ASS; see Lemma~\ref{h0h2_survives}.
\end{remark}

\subsection{Comparing the Adams and AH filtrations}
\quad

We now turn to the Atiyah--Hirzebruch spectral sequence (AHSS)
\begin{equation}\label{AHSSforMO8Meta}
E^2_{p,q}
=H_p\bigl(\mm M\eta;\pi_q(\mm{MO}\langle8\rangle)\bigr)
\Longrightarrow
\pi_{p+q}(\mm M_{\mathscr B}).
\end{equation}
The ASS and the AHSS converge to the same homotopy groups but carry different filtrations.  There is no direct morphism from one spectral sequence to the other used below.  Instead, we compare the two filtrations on the common abutment.

First, by considering the edge homomorphisms of the ASS and the AHSS, we obtain the following proposition.

\begin{proposition}\label{AHSS_Einfty_130}
In the AHSS \eqref{AHSSforMO8Meta}, 
$
(E^\infty_{13,0})_{(2)}\cong\mathbb Z_2.$
\end{proposition}

\begin{proof}
 Lemma~\ref{d2}, together with Tables~\ref{HBZ2} and~\ref{SqBtable}, gives
$
(E^3_{13,0})_{(2)}
\cong \mathbb Z_2\oplus\mathbb Z_2.$

Suppose, for contradiction, that
\[
(E^\infty_{13,0})_{(2)}
\cong\mathbb Z_2\oplus\mathbb Z_2.
\]
Then the edge homomorphism of the AHSS guarantees the existence of two elements
$
a,b\in\pi_{13}(\M_{\mathscr B})_{(2)}$
whose Hurewicz images in
$
H_{13}(\M\eta;\mathbb Z_{(2)})$
are nonzero and linearly independent. Since the mod-$2$ reduction map
\[
\rho_2:
H_{13}(\M\eta;\mathbb Z_{(2)})
\longrightarrow
H_{13}(\M\eta;\mathbb Z_2)
\]
is an isomorphism, their images under the mod-$2$ Hurewicz homomorphism are also nonzero and linearly independent.

Consequently, by the edge homomorphism of the ASS, the group \(E_\infty^{0,13}\) in the ASS must contain at least two \(\mathbb Z_2\)-summands. However, the \(E_2^{0,13}\)-term contains only one \(\mathbb Z_2\)-summand (cf. Fig.~\ref{2Fig}), which is a contradiction.
This completes the proof.
\end{proof}

We write $\{x\}\subset \bord_\ast(\eta)$ for the coset
corresponding to a class $x\in E_\infty^{s,t}$ in the ASS. By definition, every element of
$\{x\}$ maps to $x$ under the corresponding quotient map.

Let $\mm{AH}(w)$ denote the Atiyah--Hirzebruch (AH) filtration dimension of an element $w\in \bord_\ast(\eta)$. For a $2$-primary element $y\in \bord_\ast(\eta)_{(2)}$, let $\mm{Ad}(y)$ denote the Adams filtration dimension of $y$.

Since the Atiyah--Hirzebruch spectral sequence (AHSS) also
converges to $\bord_\ast(\eta)$, every element of $\{x\}$
has the AH filtration. We say that the coset $\{x\}$ has
AH filtration $p$ if every element of $\{x\}$ has AH
filtration exactly $p$. Equivalently, every element of
$\{x\}$ is provided by $E^2_{p,\ast}$ in the AHSS.

The next proposition is the basic passage from Adams filtration to AH filtration in degree $13$. 

\begin{proposition}\label{221}
$
\{h_1^2ux^2U\}$ and $\{h_1^2vxU\}$ have 
AH filtration at least $11$.
  Moreover, 
$
\{h_1vx^2U\}
$ has AH filtration at least $13$.
\end{proposition}

\begin{proof}
We first explain the argument for $\{h_1^2ux^2U\}$.

Suppose that there exists a $2$-primary element $\Pi$ detected by $h_1^2ux^2U$ such that $\mm{AH}(\Pi)=p\le10$.  

 By Proposition \ref{HBZ}, $H_{10}(\mm{M}\eta)=0$, thus we can take a subspectrum $\mm{M}\eta_{10}$ of $\M\eta$ such that the inclusion $i:\mm{M}\eta_{10}\to \mm{M}\eta$ gives that 
$$i_\ast:H_n(\mm{M}\eta_{10})\to H_n(\mm{M}\eta) $$
is an isomorphism for $n\le 10$, $H_n(\mm{M}\eta_{10})=0$ for $n>10$. Consider the following morphism
\[
\xymatrix@C=.3cm{
{E}^2_{p,q}(|_{10})=H_p(\mm{M}\eta_{10};\pi_q(\mm{MO}\langle 8\rangle))\ar@{=>}[rr]^-{}\ar[d]^-{}& & \pi_{p+q} ( \mm{MO}\langle 8\rangle \wedge \mm{M}\eta_{10})\ar[d]^-{}\\
{E}^2_{p,q}=H_p(\mm{M}\eta;\pi_q(\mm{MO}\langle 8\rangle))\ar@{=>}[rr]^-{}& & \pi_{p+q} ( \mm{MO}\langle 8\rangle\wedge \mm{M}\eta)
}
\]
By Lemma \ref{ahsspro}(1), for $p\le 10$, 
$${E}^\infty_{p,13-p}(|_{10})\to  {E}^\infty_{p,13-p}$$
is an epimorphism. By Proposition \ref{cosetspectralsequence}, there exists a $2$-primary element 
$$\bar{\Pi}\in \pi_{13}(\mm{MO}\langle 8\rangle \wedge \mm{M}\eta_{10})_{(2)}$$ provided by
${E}^2_{p,13-p}(|_{10})$ such that
$$i_\#:\pi_{13}(\mm{MO}\langle 8\rangle \wedge \mm{M}\eta_{10})\to \pi_{13}(\mm{MO}\langle 8\rangle \wedge \mm{M}\eta)$$
$$i_\# (\bar{\Pi})={\Pi}+\Pi^1$$
where $\Pi^1$ is a $2$-primary element provided by ${E}^2_{p_1,13-p_1}$,  $p_1<p$.

By the naturality of the ASS, $$\mm{Ad}(i_\#(\bar{\Pi}))\ge \mm{Ad}(\bar \Pi).$$

On the other hand, the AAHSS for $\mm M\eta_{10}$ gives
\[
\mm{Ext}_{\mathscr A_2}^{s,13+s}
\bigl(H^*(\mm M\eta_{10};\mathbb Z_2),\mathbb Z_2\bigr)=0,
\quad s=0,1,
\]
and its group in bidegree $(2,15)$ is generated by $h_2^2uU$, whose image in
$\mm{Ext}_{\mathscr A_2}^{2,15}(N,\mathbb Z_2)$ is zero. Thus 
$$\mm{Ad}({\Pi}+\Pi^1))=\mm{Ad}(i_\#(\bar{\Pi}))>2,\quad \mm{Ad}(\bar{\Pi})\ge 2,$$
 whereas $\mm{Ad}(\Pi)=2$. Therefore, we obtain a new element 
 $${\Pi}^1=i_\#(\bar{\Pi})-\Pi$$
 with $\mm{AH}(\Pi^1)=p_1<p$ and
  detected by $h_1^2ux^2U$ in the ASS.

 Repeating the argument decreases the AH filtration at each step and eventually gives an element $\Pi^*$ in negative filtration, whereas $\Pi^*\in \{h_1^2ux^2U\}$, a contradiction.  Hence $p\ge11$.

The same argument applied for $\{h^2_1vxU\}$ and $\{h_1vx^2U\}$
 shows the desired results. In particular, for the argument of $\{h_1vx^2U\}$, we take the subspectrum $\mm{M}\eta_{12}$ of $\M\eta$ such that the inclusion $i:\mm{M}\eta_{12}\to \mm{M}\eta$ gives that 
$$i_\ast:H_n(\mm{M}\eta_{12})\to H_n(\mm{M}\eta) $$
is an isomorphism for $n\le 12$, $H_n(\mm{M}\eta_{10})=0$ for $n>12$.  
\end{proof}

\begin{corollary}\label{E131infty}
In the AHSS \eqref{AHSSforMO8Meta},
$
E^\infty_{13,1}=E^3_{13,1}=\mathbb Z_2.$
\end{corollary}

\begin{proof}
Lemma~\ref{d2}, together with Tables~\ref{HBZ2}, \ref{SqBtable}, and~\ref{piMO8}, gives
$$
E^3_{13,1}=\mathbb Z_2,\quad E^3_{14,0}=0.$$
By Proposition~\ref{221} and $E^3_{14,0}=0$, an element of $\{h_1vx^2U\}$ has AH filtration $13$. Thus $E^\infty_{13,1}\ne 0$, and hence $E^\infty_{13,1}=E^3_{13,1}=\Z_2$.
\end{proof}

The first hidden extension is now forced by the $h_0$-tower in Lemma~\ref{h0module}.

\begin{corollary}\label{011}
In the AHSS \eqref{AHSSforMO8Meta}, the following assertions hold:
\begin{enumerate}
    \item $\{vx^2U\}$ is provided by $E^2_{13,0}$;
    \item $\{h_1x^2\Sq^1uU\}$ is provided by $E^2_{12,1}$;
    \item $\{h_1^2ux^2U\}$ and $\{h_1^2vxU\}$ are provided by $E^2_{11,2}$.
\end{enumerate}
\end{corollary}

\begin{proof}
We first record the relevant terms of the AHSS.
By Proposition~\ref{AHSS_Einfty_130} and
Lemma~\ref{d2}, together with Tables~\ref{HBZ2} and
\ref{SqBtable}, we have
\[
(E^\infty_{13,0})_{(2)}\cong\mathbb Z_2,\qquad
E^3_{12,1}\cong\mathbb Z_2,\qquad
E^3_{11,2}\cong\mathbb Z_2\oplus\mathbb Z_2.
\]
We work throughout with the $2$-primary components.

\medskip
\noindent
\textbf{(1)} Let $a\in\{vx^2U\}$.
By Lemma~\ref{h0module}, the elements $2a$ and $4a$
belong to the cosets
$\{h_1x^2\Sq^1uU\},$
$\{h_1^2ux^2U\},$
respectively.

By Proposition~\ref{221}, every element of
$\{h_1^2ux^2U\}$ has AH filtration at least $11$.
We claim that $\mm{AH}(a)=13$.
Indeed, if $\mm{AH}(a)\le 12$, then,
since the relevant $E^\infty$-terms in filtrations
$12$ and $11$ are annihilated by $2$, successive
multiplication by $2$ would give $
\operatorname{AH}(4a)\leq 10$.
This contradicts Proposition~\ref{221}, and thus completes the proof.

\medskip
\noindent
\textbf{(2)} Let $b\in\{h_1x^2\Sq^1uU\}$.
By Lemma~\ref{h0module}, we have
$
2b\in\{h_1^2ux^2U\}.$
Since every element of the latter coset has AH filtration
at least $11$, and $E^\infty_{11,2}$ is annihilated by
$2$, it follows that $b$ has AH filtration at least $12$.

Suppose that $b$ has AH filtration $13$.
Choose any $a\in\{vx^2U\}$.
By part~(1), both $a$ and $b$ have AH filtration $13$.
Since $E^\infty_{13,0}\cong\mathbb Z_2$, their difference
satisfies
\[
\operatorname{AH}(a-b)\leq 12.
\]
On the other hand, $\mm{Ad}(b)=1>\mm{Ad}(a)=0$, so 
\[
a-b\in\{vx^2U\}.
\]
This contradicts part~(1). Therefore every element of
$\{h_1x^2\Sq^1uU\}$ has AH filtration exactly $12$.
In particular,
$
E^\infty_{12,1}\cong\mathbb Z_2.$

\medskip
\noindent
\textbf{(3)} The proof is analogous to  part~(2).
\end{proof}

\subsection{Elements with high filtration}\quad

\begin{lemma}\label{ASS_cox2U}
	The class $c_0x^2U$ survives in the ASS for $\bord_\ast(\eta)$.
\end{lemma}
\begin{proof}
	By Fig.\ref{2Fig} and Lemma \ref{h0vx2survive}, there is no incoming Adams differential for $c_0x^2U$. Next we consider the outgoing Adams differentials. 

By Table \eqref{E4smn} and Lemma \ref{pert4}, for $s\ge 5$, the Adams $E_2^{s,11+s}$-page
 only contains the nonzero classes $h_2\omega U$ and $h_0h_2\omega U$. By Fig.\ref{1Fig}, we have
$$d_2(h_0^s\tau)=h_0^sh_2\omega$$
in the ASS for $\bord_\ast(\eta)$. Thus, we first have $d_2(c_0x^2U)=0$. Moreover, since the Adams term $E_3^{s,s+11}$ vanishes for $s\ge 5$, all outgoing Adams differentials for $c_0x^2U$ are zero. So,   
the class $c_0x^2U$ survives in the ASS.  
\end{proof}

\begin{proposition}\label{48}
  In the AHSS \eqref{AHSSforMO8Meta},
$\Z_2\subset E^2_{4,8}$  survives to $E^\infty$.
\end{proposition}

\begin{proof}
Let $\M\eta|_4$ denote the $4$-skeleton  of $\M\eta$. As the computation in the proof for Lemma \ref{ASS_cox2U}, we have
 that $c_0x^2U$ represents a nonzero class in $\Ext_{\mathscr A_2}^{3,15}(H^\ast(\M\eta|_4;\Z_2),\Z_2)$. Furthermore, 
  it survives in the ASS for $\bord_{\ast}(\eta|_4)$.
  $\{c_0x^2U\}$ has AH filtration at most $4$ since $H_i(\M\eta|_4;\Z)=0$ for $i\ge 5$.
  
 Thus, by Lemma \ref{ASS_cox2U}, together with naturality of the ASS and the AHSS, there exist elements in $\{c_0x^2U\}\subset \bord_{12}(\eta)$ with AH filtration at most $4$.

Suppose that there exists an element $\Pi\in \{c_0x^2U\}$ with $\mm{AH}(\Pi)\le 3$. By using the subspectrum $\M\eta|_2$ and the argument as Lemma \ref{221}, we can 
deduce a contradiction. 

Therefore, there exists an element of $\{c_0x^2U\}$ provided by $\Z_2\subset E_{4,8}^2$ in the AHSS. Thus $\Z_2\subset E_{4,8}^2$ survives to $E^\infty$. Indeed, the free part $\Z\subset E_{4,8}^2$ provides the elements detected by $h_0^swx^2U$ for $s\ge 0$.  
 \end{proof}

\begin{lemma}\label{Adams_to_AHSS49}
	If a $2$-primary class of $\bord_{13}(\eta)$ has Adams filtration at least $4$, then it has AH filtration most $4$.  
\end{lemma}
\begin{proof}
	We also consider the $4$-skeleton $\M\eta|_4$ of $\M\eta$. It is easy to check that 
	$$h_1\omega x^2U, h_1c_0x^2U, h_0^sh_2\omega xU\in \Ext_{\mathscr A_2}^{\ast,\ast+13}(H^\ast(\M\eta|_4;\Z_2),\Z_2).$$
For $s\ge 6,$ the groups $\Ext_{\mathscr A_2}^{s,s+12}(H^\ast(\M\eta|_4;\Z_2),\Z_2)$ have generators 
$$h_0^s\tau U\quad \text{and}\quad h_0^s\omega x^2U.$$
By the $h_0$-module structure, all outgoing differentials for $$h_1\omega x^2U, h_1c_0x^2U, h_0^sh_2\omega xU$$
are zero. Thus, the homomorphism between Adams $E_\infty$-terms
$$E_\infty^{s,s+13}(\M\eta|_4)\to E_\infty^{s,s+13}(\M\eta)$$
is surjective for $s\ge 4$. 
This implies that any $2$-primary element of $\bord_{13}(\eta)$ with Adams filtration at least $4$ has a lift to $\bord_{13}(\eta|_4)$. 

By the homology of $\M\eta|_4$, any element of $\bord_{13}(\eta|_4)$ has AH filtration at most $4$. This implies the desired result by naturality.
\end{proof}

\section{The $ko$-homology of $\mm{M}\eta$}\label{koMketa}

\subsection{The filtration comparison in degree $13$}\label{koMeta85_quot} \quad

Let
$
\alpha_8:
\mm{MO}\langle8\rangle
\longrightarrow
\mm{MSpin}
\xrightarrow{\alpha}
ko$
be the composite of the canonical map to spin bordism with the
Atiyah--Bott--Shapiro orientation \cite{AtiyahBottShapiro}.  Smashing
with $\mm{M}\eta$ induces
\[
(\alpha_8)_*:
\Omega_*^{\mm O\langle8\rangle}(\eta)
\longrightarrow
ko_*(\mm M\eta).
\]

For the AHSS of a connective spectrum $E$, we use the convention
\[
E^2_{p,q}(E,\mm M\eta)
=
H_p(\mm M\eta;\pi_q(E))
\Longrightarrow
E_{p+q}(\mm M\eta),
\]
with differential
$
d_r:E^r_{p,q}\longrightarrow E^r_{p-r,q+r-1}.$

For the two spectral sequences considered below, let
\[
F_{p,q}
\subset
\Omega_{p+q}^{\mm O\langle8\rangle}(\eta),
\qquad
\bar F_{p,q}
\subset
ko_{p+q}(\mm M\eta)
\]
denote the corresponding filtration subgroups (cf. Section \ref{sec:spectral-prelim}).  Thus
\[
F_{p,q}/F_{p-1,q+1}
\cong
E^\infty_{p,q}(\mm{MO}\langle8\rangle,\mm M\eta)
\]
and
\[
\bar F_{p,q}/\bar F_{p-1,q+1}
\cong
E^\infty_{p,q}(ko,\mm M\eta).
\]


The purpose of this subsection is to prove the following detection result.

\begin{theorem}\label{pFZ4}
 The following homomorphism induced by $\alpha_8$
\[
(\alpha_8)_*:
\frac{\Omega_{13}^{\mm O\langle8\rangle}(\eta)}
{F_{8,5}}
\longrightarrow
\frac{ko_{13}(\mm M\eta)}
{\bar F_{8,5}},
\]
is injective.
\end{theorem}

The only non-formal point in the proof is that the
class in filtration $10$ of the bordism AHSS has a nonzero image in filtration $9$ of the $ko$-AHSS.  We first isolate the relevant
$ko$-class.

\begin{lemma}\label{koketa94}
In the AHSS for $ko_*(\mm M\eta)$,
$
E^\infty_{9,4}(ko,\mm M\eta)
\cong
\mathbb Z_6\oplus\mathbb Z_s.$
\end{lemma}

\begin{proof}
By Tables~\ref{HBZ2} and~\ref{SqBtable}, the $d_3$ in the $ko$-AHSS, dual to $\Sq^2\Sq^1$ \cite{Shen2026Differentials}, gives
\[
E^4_{9,4}(ko,\mm M\eta)
\cong
\mathbb Z_6\oplus\mathbb Z_s.
\]
We show that no later differential affects this group.

It is easy to show $E^3_{14,0}=0$. So
the only possible incoming differential is
\[
d_4:
E^4_{13,1}(ko,\mm M\eta)
\longrightarrow
E^4_{9,4}(ko,\mm M\eta).
\]
Since $\alpha_8$ is an isomorphism on the coefficient groups in
degrees $0,1,2$, naturality identifies
$E^3_{13,1}(ko,\mm M\eta)$ with $E^3_{13,1}(\mm{MO}\langle8\rangle,\mm M\eta)$. The latter is $\Z_2$ and survives to
$E^\infty$ by Corollary~\ref{E131infty}; hence this $d_4$ is zero.

We next consider the possible outgoing differentials.  The differential
\[
d_5:
E^5_{9,4}(ko,\mm M\eta)
\longrightarrow
E^5_{4,8}(ko,\mm M\eta)
\]
is zero because the source is finite whereas
$E^5_{4,8}(ko,\mm M\eta)$ is torsion-free.  The $d_6$- and
$d_8$-targets vanish for degree reasons.

The only remaining possible nontrivial differential before $d_9$ is
\[
d_7:
E^7_{9,4}(ko,\mm M\eta)
\longrightarrow
E^7_{2,10}(ko,\mm M\eta).
\]
We claim that
\[
E^8_{2,10}(ko,\mm M\eta)
=
E^2_{2,10}(ko,\mm M\eta).
\]
This follows by comparison with the AHSS for $ko_*(\mm M\xi)$, where
$\xi$ is the bundle over $\CP^\infty$ whose pullback to $B$ is
$\eta$. The corresponding term satisfies
\[
E^8_{2,10}(ko,\mm M\xi)
=
E^2_{2,10}(ko,\mm M\xi),
\]
and naturality gives the asserted survival for $\mm M\eta$.
Consequently $d_7=0$.

Finally,
$
d_9:
E^9_{9,4}(ko,\mm M\eta)
\longrightarrow
E^9_{0,12}(ko,\mm M\eta)$
is zero, since the source is finite while
$
E^\infty_{0,12}(ko,\mm M\eta)
=
E^2_{0,12}(ko,\mm M\eta)
\cong\mathbb Z.$

Therefore, all differentials for $ 
E^r_{9,4}(ko,\mm M\eta)$ are zero where $r\ge 4$.
\end{proof}

Let $\bar\gamma\in ko_{13}(\mm M\eta)$ be an element
detected by a generator of the $\Z_6$-summand of
$E^\infty_{9,4}(ko,\mm M\eta)$, 
 whose order is of the form
$2^a3^b$ for some $a,b\ge 1$.
In the AHSS for $ko_\ast(\M\eta)$, since all $E^2_{p,13-p}$-terms are finite groups for $p\ge 9$, the filtration $\bar F_{9,4}$ is also a finite group. Hence, such an element $\bar \gamma$ exists.  

We next show that $\bar\gamma$ lies in the image of 
$\alpha_{8\ast}:\bord_{13}(\eta)\to ko_{13}(\M\eta). $

Consider the cofibration sequence
\begin{equation}\label{MO8koWcofib}
\mm{MO}\langle8\rangle
\xrightarrow{\alpha_8}
ko
\xrightarrow{h}
W
\xrightarrow{\partial}
\Sigma\mm{MO}\langle8\rangle.
\end{equation}
The associated long exact sequence shows that $W$ is $3$-connected,
\[
\pi_5(W)=\pi_6(W)=0,
\qquad
\pi_7(W)\cong\mathbb Z_2.
\]
We shall also use the following information.

\begin{lemma}\label{piW}
\begin{enumerate}[(i)]
\item
The homomorphism
\[
h_*:\pi_4(ko)\longrightarrow\pi_4(W)
\]
is multiplication by $24$.

\item
Under
\[
\partial_*:
\pi_n(W)
\longrightarrow
\pi_n(\Sigma\mm{MO}\langle8\rangle),
\]
the generators of
\[
\pi_9(W)\cong \Z_2\quad \text{and}\quad \pi_{10}(W)\cong\mathbb Z_2
\]
are mapped to $c_0$ and $h_1c_0$, respectively.  Moreover, a
generator of
\[
\pi_{11}(W)\cong\mathbb Z_3
\]
is mapped onto the subgroup of order $3$ in
\[
\pi_{11}(\Sigma\mm{MO}\langle8\rangle)
=
\pi_{10}(\mm{MO}\langle8\rangle)
\cong\mathbb Z_6.
\]

\item
There is an exact sequence
\begin{equation}\label{exactWMO8}
0
\longrightarrow
\pi_{13}(W)
\xrightarrow{\partial_*}
\pi_{12}(\mm{MO}\langle8\rangle)
\longrightarrow
\pi_{12}(ko)
\longrightarrow
\pi_{12}(W)
\longrightarrow0.
\end{equation}
\end{enumerate}
\end{lemma}

\begin{proof}
The assertions in (ii) and (iii) follow immediately from the long exact
homotopy sequence associated with \eqref{MO8koWcofib}, together with
the known homotopy groups of
$\mm{MO}\langle8\rangle$ and $ko$.

It remains to prove (i).  Since
$H_4(ko)\cong\mathbb Z,$
applying ordinary homology to \eqref{MO8koWcofib} gives
$
H_4(W)\cong\mathbb Z.$
Since $W$ is $3$-connected, the Hurewicz theorem gives
$
\pi_4(W)\cong H_4(W)\cong\mathbb Z.$
The long exact homotopy sequence therefore contains
\[
0
\longrightarrow
\pi_4(ko)\cong\mathbb Z
\xrightarrow{h_*}
\pi_4(W)\cong\mathbb Z
\longrightarrow
\pi_3(\mm{MO}\langle8\rangle)
\cong\mathbb Z_{24}
\longrightarrow0.
\]
Hence $h_*$ is multiplication by $24$.
\end{proof}

\begin{lemma}\label{kerkotoW}
The element $\bar\gamma$ lies in the image of
\[
(\alpha_8)_*:
\Omega_{13}^{\mm O\langle8\rangle}(\eta)
\longrightarrow
ko_{13}(\mm M\eta).
\]
\end{lemma}

\begin{proof}
Applying
$\pi_*(-\wedge\mm M\eta)$
to the cofibration sequence \eqref{MO8koWcofib} gives the exact
sequence
\[
\pi_{13}
\bigl(
\mm{MO}\langle8\rangle\wedge\mm M\eta
\bigr)
\xrightarrow{(\alpha_8)_*}
\pi_{13}(ko\wedge\mm M\eta)
\xrightarrow{h_*}
\pi_{13}(W\wedge\mm M\eta).
\]
Therefore it is enough to prove
$
h_*(\bar\gamma)=0.$

Let
$
\mathtt{k}_{(2)},$
$\mathtt{k}_{(3)}$ denote the $2,3$-primary components of $\bar\gamma$. 
It is therefore enough to prove
\[
h_*(\mathtt{k}_{(2)})=0,
\qquad
h_*(\mathtt{k}_{(3)})=0.
\]

We first consider $\mathtt{k}_{(2)}$.  Suppose, to the
contrary, that
\[
h_*(\mathtt{k}_{(2)})\ne0.
\]
Then $h_*(\mathtt{k}_{(2)})$ is an element of order $2^a$ in
$
\pi_{13}(W\wedge\mm M\eta).$
We determine which $E^\infty$-term can detect this element.

By Lemma~\ref{piW}(i),
\[
h_*:
\pi_4(ko)\longrightarrow\pi_4(W)
\]
is multiplication by $24$.  Hence the induced homomorphism
\[
h_*:
E^2_{9,4}(ko,\mm M\eta)
\longrightarrow
E^2_{9,4}(W,\mm M\eta)
\]
vanishes on the $\mathbb Z_6$-summand.  Consequently, if
$h_*(\mathtt{k}_{(2)})$ is nonzero, then its AH filtration is
strictly smaller than $9$.

By Lemma~\ref{piW}, among the terms
\[
E^\infty_{p,13-p}(W,\mm M\eta),
\qquad p<9,
\]
the only possible nontrivial $2$-primary term is
\[
E^\infty_{4,9}(W,\mm M\eta).
\]
There is one possible filtration-zero term to consider.  If
$\pi_{13}(W)\ne0$, then Lemma~\ref{piW}(iii) shows that
$\pi_{13}(W)$ is a subgroup of
$
\pi_{12}(\mm{MO}\langle8\rangle)\cong\mathbb Z.
$
Thus it is free.  Comparing the AHSSs induced by $\partial$, together
with Lemma~\ref{E_infty_012}, shows that
\[
E^\infty_{0,13}(W,\mm M\eta)
=
E^2_{0,13}(W,\mm M\eta),
\]
and this term is free.  It therefore cannot
detect the torsion element $h_*(\mathtt{k}_{(2)})$.

Consequently,
$
h_*(\mathtt{k}_{(2)})$
must be detected by
$
E^\infty_{4,9}(W,\mm M\eta).$ And we have
$$E^\infty_{4,9}(W,\mm M\eta)=E^2_{4,9}(W,\mm M\eta)=\Z_2.$$
Now consider the morphism of AHSSs induced by
$
\partial:
W\longrightarrow
\Sigma\mm{MO}\langle8\rangle.$
By Lemma~\ref{piW}(ii), the generator of
$\pi_9(W)\cong\mathbb Z_2$
is mapped to the nonzero class
\[
c_0
\in
\pi_9(\Sigma\mm{MO}\langle8\rangle)
=
\pi_8(\mm{MO}\langle8\rangle).
\]
Hence the induced homomorphism
\[
\partial_*:
E^2_{4,9}(W,\mm M\eta)
\longrightarrow
E^2_{4,9}
(\Sigma\mm{MO}\langle8\rangle,\mm M\eta)
\]
is nonzero.  By Proposition~\ref{48}, the corresponding
$\mathbb Z_2$-class on the target side survives to $E^\infty$.
Naturality therefore implies that
\[
\partial_*:
E^\infty_{4,9}(W,\mm M\eta)
\longrightarrow
E^\infty_{4,9}
(\Sigma\mm{MO}\langle8\rangle,\mm M\eta)
\]
is nonzero.  Then, we obtain
$\partial_*h_*(\mathtt{k}_{(2)})\ne0.$
This contradicts
\[
\partial_*\circ h_*=0
\]
in the long exact sequence associated with
\eqref{MO8koWcofib}.  Therefore
$
h_*(\mathtt{k}_{(2)})=0.$

We next consider the $3$-primary component.  Suppose that
\[
h_*(\mathtt{k}_{(3)})\ne0.
\]
As above, $\mm{AH}(h_*(\mathtt{k}_{(3)}))<9$.
By Lemma~\ref{piW}, the only possible nontrivial $3$-primary term in
total degree $13$ and filtration smaller than $9$ is
$
E^\infty_{2,11}(W,\mm M\eta).$
Hence $h_*(\mathtt{k}_{(3)})$ must be detected by this term which  satisfies
$$
E^\infty_{2,11}(W,\mm M\eta)=
E^2_{2,11}(W,\mm M\eta)=\Z_3.$$

By Lemma~\ref{piW}(ii), the generator of
$
\pi_{11}(W)\cong\mathbb Z_3
$
is mapped by $\partial_*$ onto the subgroup of order $3$ in
$
\pi_{11}(\Sigma\mm{MO}\langle8\rangle)
\cong\mathbb Z_6.
$
The corresponding class survives in the AHSS for
$\Sigma\mm{MO}\langle8\rangle\wedge\mm M\eta$ by Lemma \ref{b0x2}.  Thus naturality gives
$
\partial_*h_*(\mathtt{k}_{(3)})\ne0,$
again contradicting
$
\partial_*\circ h_*=0.$
Therefore
$
h_*(\mathtt{k}_{(3)})=0.
$

It follows that
$
h_*(\bar\gamma)=0.$
By exactness, $\bar\gamma$ lies in the image of $(\alpha_8)_*$.
\end{proof}

We can now prove Theorem~\ref{pFZ4}.

\begin{proof}[Proof of Theorem~\ref{pFZ4}]

  Since $\alpha_8$ induces isomorphisms on $\pi_i$ for $0\le i\le 2$,
   naturality gives injections
\[
E^\infty_{p,13-p}
(\mm{MO}\langle8\rangle,\mm M\eta)
\longrightarrow
E^\infty_{p,13-p}
(ko,\mm M\eta),
\qquad
11\le p\le13.
\]

To see how these associated-graded injections reduce the problem,
consider first the commutative diagram
\[
\xymatrix{
0\ar[r]
&
F_{12,1}/F_{8,5}
\ar[r]\ar[d]
&
F_{13,0}/F_{8,5}
\ar[r]\ar[d]^{(\alpha_8)_*}
&
F_{13,0}/F_{12,1}
\ar[r]\ar[d]
&
0
\\
0\ar[r]
&
\bar F_{12,1}/\bar F_{8,5}
\ar[r]
&
\bar F_{13,0}/\bar F_{8,5}
\ar[r]
&
\bar F_{13,0}/\bar F_{12,1}
\ar[r]
&
0.
}
\]
The right vertical homomorphism is the map on
$
E^\infty_{13,0},$
and is therefore injective.  Repeating the same argument for the
successive filtration quotients
\[
E^\infty_{12,1}
\qquad\text{and}\qquad
E^\infty_{11,2},
\]
we reduce the problem to proving that
\[
(\alpha_8)_*:
F_{10,3}/F_{8,5}
\longrightarrow
\bar F_{10,3}/\bar F_{8,5}
\]
is injective.  

We now treat the exceptional filtration-$10$ term.  On the $ko$-side, since
$
\pi_3(ko)=0,$
we have
$
\bar F_{10,3}=\bar F_{9,4}.$
By Lemma~\ref{koketa94},
\[
\bar F_{10,3}/\bar F_{8,5}
=
\bar F_{9,4}/\bar F_{8,5}
\cong
E^\infty_{9,4}(ko,\mm M\eta)
\cong
\mathbb Z_6\oplus\mathbb Z_s.
\]

By Lemma~\ref{kerkotoW}, choose
$
\gamma\in\Omega$
such that
$
(\alpha_8)_*(\gamma)=\bar\gamma.$
We claim that
\[
\gamma
\in
F_{10,3}\setminus F_{9,4}.
\]

Suppose first that
$
\gamma\in F_{9,4}.$
Since
$
E^2_{9,4}
(\mm{MO}\langle8\rangle,\mm M\eta)=0,$
$
F_{9,4}=F_{8,5}.$
Naturality would then imply
$
\bar\gamma\in\bar F_{8,5},$
contradicting the fact that $\bar\gamma$ represents a nonzero class in
$
\bar F_{9,4}/\bar F_{8,5}.$
Thus
$
\gamma\notin F_{9,4}.$

Suppose that
$
\gamma\notin F_{10,3}.$
Then 
$\mm{AH}(\gamma)=p\ge11$.  Since the map
\[
E^\infty_{p,13-p}
(\mm{MO}\langle8\rangle,\mm M\eta)
\longrightarrow
E^\infty_{p,13-p}
(ko,\mm M\eta)
\]
is injective for $p\ge11$, the image $\bar\gamma$ would satiafy $\mm{AH}(\bar \gamma)=p\ge 11$.  This contradicts $\mm{AH}(\bar \gamma)=9$.  Hence
$
\gamma\in F_{10,3}.$

Therefore, $\gamma$ represents a nonzero element of
\[
F_{10,3}/F_{9,4}\cong E^\infty_{10,3}
(\mm{MO}\langle8\rangle,\mm M\eta).
\]
 Since the image $\bar \gamma$ contains nontrivial $2$ and $3$ primary components, $\gamma$ also has such primary components. Combining this with Lemma \ref{E3103} gives 
 \begin{equation}\label{Einfty_103}
 	E^\infty_{10,3}
(\mm{MO}\langle8\rangle,\mm M\eta)=E^3_{10,3}
(\mm{MO}\langle8\rangle,\mm M\eta)=\Z_2\oplus \Z_3.
 \end{equation}
  Moreover, by the group structure of $\Z_6$, the induced map
\[
F_{10,3}/F_{8,5}
\longrightarrow
\bar F_{10,3}/\bar F_{8,5}
\]
 is injective. This finishes the proof.
\end{proof}

\begin{lemma}\label{h0h2_survives}
In the ASS for $\bord_\ast(\eta)$, the class
\[
h_0^2h_2(\Sq^2\Sq^1u+x\Sq^1u)U
\in
\Ext_{\mathscr A_2}^{3,16}
\bigl(H^\ast(\M\eta;\Z_2),\Z_2\bigr)
\]
survives to $E_\infty$ (cf. Fig.\ref{2Fig}). Moreover, there is a nontrivial extension
$$0\to E^\infty_{10,3}(\MO\langle8\rangle,\M\eta)\to F_{11,2}/F_{8,5}\to  E^\infty_{11,2}(\MO\langle8\rangle,\M\eta)\to 0.$$
\end{lemma}

\begin{proof}
By \eqref{Einfty_103}, there is a nonzero element $\alpha\in {\bord_{13}(\eta)}_{(2)}$ satisfying $$\mm{AH}(\alpha)=10.$$
In the ASS for $\bord_\ast(\eta)$, by Corollary~\ref{011}, the elements detected by
\[
vx^2U,\quad
h_1x^2\Sq^1uU,\quad
h_1^2ux^2U,\quad
h_1^2vxU
\]
have AH filtrations at least $11$. On the other hand,
Lemma~\ref{Adams_to_AHSS49} shows that an element of Adams
filtration at least $4$ has AH filtration at most $4$.
Consequently, none of these classes can detect $\alpha$.
Recall  Fig.\ref{2Fig}. The only remaining candidate is
\[
h_0^2h_2(\Sq^2\Sq^1u+x\Sq^1u)U.
\]
Since $\alpha$ is nonzero, this class must survive to
$E_\infty$ and detect $\alpha$.

 Lemma \ref{h0module} gives the $h_0$-action
$$h_0h_1^2vxU=h_0^2h_2(\Sq^2\Sq^1u+x\Sq^1u)U.$$
Let $\beta\in {\bord_{13}(\eta)}_{(2)}$ be detected by $h_1^2vxU$ in the ASS. Then $2\beta$ is detected by $h_0^2h_2(\Sq^2\Sq^1u+x\Sq^1u)U$. Hence $$2\beta=\alpha+\bar \alpha\in {\bord_{13}(\eta)}_{(2)}$$ where $\mm{Ad}(\bar \alpha)\ge 4$. By Lemma \ref{Adams_to_AHSS49}, $\bar \alpha\in F_{4,9}$. Thus, $\mm{AH}(2\beta)=10$.  
Recall 
$$F_{10,3}/F_{8,5}\cong E^\infty_{10,3}(\MO\langle8\rangle,\M\eta).$$
Combining this with the AH filtrations for $\beta$ (cf. Corollary \ref{011}) and $2\beta$, we obtain the desired extension.
\end{proof}

\begin{lemma}\label{E_infty_76}
In the AHSS for $\bord_\ast(\eta)$,
$
E^\infty_{7,6}=0.$
\end{lemma}

\begin{proof}
Recall Table \eqref{E2AHSSforMO8Meta}:
$
E^2_{7,6}\cong\Z_2.$
Hence, if $E^\infty_{7,6}\ne0$, there would exist a nonzero
element of $\bord_{13}(\eta)_{(2)}$ with AH filtration
exactly $7$.

 The preceding ASS calculations (cf. Fig.\ref{2Fig}),
together with Corollary~\ref{011} and
Lemma~\ref{h0h2_survives}, show that elements of Adams
filtration at most $3$ have AH filtration at least $10$.
Meanwhile, by Lemma~\ref{Adams_to_AHSS49}, elements of Adams
filtration at least $4$ have AH filtration at most $4$.

Thus no element of $\bord_{13}(\eta)_{(2)}$ can have AH
filtration exactly $7$. This contradiction proves that
$
E^\infty_{7,6}=0.$
\end{proof}

\subsection{The two- and three-primary parts of
\texorpdfstring{$ko_9(\mm M\eta)$}{ko9(M eta)}}\label{subsec:ko9}
\quad

In this subsection, we show that $ko_9(\mm M\eta)$ contains a summand $\Z_{24}$, which will be used to prove Proposition \ref{ttprime_differential}.
 We first use the
AHSS to determine the possible group orders and then the ASS to resolve the extensions.

We first recall the Adams spectral sequence
$$\Ext_{\mathscr A}^{s,t}(H^\ast(ko\wedge \M\eta;\Z_2),\Z_2)\Longrightarrow {_{2}\pi}_{t-s}(ko\wedge \M\eta)$$
 Let
$N=H^*(\mm M\eta;\mathbb Z_2)$.
Let $\mathscr A_1\subset\mathscr A$ be the subalgebra of the
mod-$2$ Steenrod algebra generated by $\Sq^1$ and $\Sq^2$.
The change-of-rings gives
\begin{equation}\label{eq:ko9-ASS}
\Ext_{\mathscr A}^{s,t}(H^\ast(ko\wedge \M\eta;\Z_2),\Z_2)\cong\Ext_{\mathscr A_1}^{s,t}(N,\mathbb Z_2).
\end{equation}

The relevant part of the $\mathscr A_1$-module $N$ is generated in
degrees $7$--$10$ by
\[
uU,\quad (\Sq^1u)U,\quad uxU,\quad vU,\quad
(x\Sq^1u)U,\quad (\Sq^2\Sq^1u)U.
\]
The Steenrod actions needed below are
\begin{equation}\label{eq:ko9-A1-actions}
\begin{array}{c|cc}
 &\Sq^1&\Sq^2\\ \hline
uU &(\Sq^1u)U&0\\
(\Sq^1u)U&0&(x\Sq^1u+\Sq^2\Sq^1u)U\\
uxU&(x\Sq^1u)U&\text{degree }11\\
vU&(\Sq^2\Sq^1u)U&\text{degree }11
\end{array}
\end{equation}
These identities follow from Table~\eqref{HBZ2} \eqref{SqBtable} and 
$\Sq^iU=w_i(\eta)U$.  

\begin{lemma}\label{ko9Meta2primary}
In stem $9$, the Adams $E_2$-term of \eqref{eq:ko9-ASS} is
\[
\Ext_{\mathscr A_1}^{s,s+9}(N,\mathbb Z_2)
\cong
\begin{cases}
\mathbb Z_2\langle\lambda_{ux},\lambda_v\rangle,&s=0,\\
\mathbb Z_2\langle h_0\lambda_v\rangle,&s=1,\\
\mathbb Z_2\langle h_0^2\lambda_v\rangle,&s=2,\\
0,&s\ge3,
\end{cases}
\]
where $\lambda_{ux}$ and $\lambda_v$ denote the filtration-zero
functionals dual to $uxU$ and $vU$, respectively.  All four classes survive in the ASS, and
\begin{equation}\label{eq:ko9-two-primary}
ko_9(\mm M\eta)_{(2)}\cong\mathbb Z_8\oplus\mathbb Z_2,
\qquad
ko_9(\mm M\eta)_{(3)}\cong\mathbb Z_3.
\end{equation}
In particular, every element detected by $\lambda_v$ has order $8$.
\end{lemma}

\begin{proof}
We begin with the AHSS
\[
E^2_{p,q}(ko,\mm M\eta)
=H_p(\mm M\eta;\pi_q(ko))
\Longrightarrow ko_{p+q}(\mm M\eta).
\]
By Proposition~\ref{HBZ}, the nonzero
$2$-primary terms in total degree $9$ on the $E^2$-page are
\[
\begin{aligned}
E^2_{9,0}(ko,\mm M\eta)_{(2)}&\cong(\mathbb Z_2)^2,\\
E^2_{8,1}(ko,\mm M\eta)&\cong\mathbb Z_2,\\
E^2_{7,2}(ko,\mm M\eta)&\cong\mathbb Z_2,\\
E^2_{0,9}(ko,\mm M\eta)&\cong\mathbb Z_2.
\end{aligned}
\]
The last term is killed by
\[
d_2:E^2_{2,8}(ko,\mm M\eta)
\longrightarrow E^2_{0,9}(ko,\mm M\eta),
\]
which is nonzero because it is dual to $\Sq^2U=xU$.
The $d_2$-differentials involving the remaining terms are zero:
in particular, the maps out of $E^2_{9,0}$ and into
$E^2_{7,2}$ are dual to $\Sq^2(uU)=0$.
The other possible differentials either have zero targets or map
torsion groups into free abelian groups.  Consequently, on the
$2$-primary part,
\[
E^\infty_{9,0}\cong(\mathbb Z_2)^2,
\qquad E^\infty_{8,1}\cong\mathbb Z_2,
\qquad E^\infty_{7,2}\cong\mathbb Z_2.
\]
In particular,
\begin{equation}\label{eq:ko9-ahss-order}
\bigl|ko_9(\mm M\eta)_{(2)}\bigr|=16.
\end{equation}

At the prime $3$, the only nonzero term in total degree $9$ is
\[
E^2_{9,0}(ko,\mm M\eta)_{(3)}\cong\mathbb Z_3.
\]
All outgoing differentials either have zero targets or map
torsion groups into free abelian groups.
  Hence
\begin{equation}\label{eq:ko9-three-primary}
ko_9(\mm M\eta)_{(3)}\cong E^2_{9,0}(ko,\mm M\eta)_{(3)}=\mathbb Z_3.
\end{equation}

We now determine the $2$-primary extensions. We mainly use the AAHSS 
$$E_1^{s,m,n}=\Ext_{\mathscr A_1}^{s,s+m}(\mathbb Z_2,\mathbb Z_2)\otimes N^n\Longrightarrow \Ext_{\mathscr A_1}^{s,s+m+n}(N,\mathbb Z_2),$$
 For $\Ext_{\mathscr A_1}^{s,s+m}(\mathbb Z_2,\mathbb Z_2)$ with $m\le 10$, one can refer to \cite[p.66]{Ravenel1986}.

  The actions
\eqref{eq:ko9-A1-actions} determine the differentials in the
AAHSS.  The resulting dimensions are
\[
\begin{array}{c|cccccc}
s&0&1&2&3,4&5&s\ge6\\ \hline
\dim E_{1,\mathrm{AAHSS}}\text{ in stem }9&2&3&3&2&3&2\\
\dim\Ext_{\mathscr A_1}^{s,s+9}(N,\mathbb Z_2)&2&1&1&0&0&0
\end{array}
\]
Thus the asserted dimensions of the Adams $E_2$-term follow from
the AAHSS.

It remains to identify the multiplication by $h_0$.  We use the
normalized bar complex only for this purpose.  By \eqref{eq:ko9-A1-actions} and the Adem relation
$\Sq^2\Sq^2=\Sq^1\Sq^2\Sq^1$, the expressions
\[
[\Sq^2\Sq^1]uU+[\Sq^1]uxU+[\Sq^1]vU
\]
and
\[
\begin{aligned}
&[\Sq^1|\Sq^2\Sq^1]uU+[\Sq^2|\Sq^2]uU\\
&\qquad+[\Sq^1|\Sq^1]uxU+[\Sq^1|\Sq^1]vU
\end{aligned}
\]
are bar cycles.  They pair nontrivially with both
$h_0\lambda_{ux}$ and $h_0\lambda_v$, and with both
$h_0^2\lambda_{ux}$ and $h_0^2\lambda_v$, respectively.
Thus we conclude that
\[
h_0\lambda_{ux}=h_0\lambda_v\ne0,
\qquad
h_0^2\lambda_{ux}=h_0^2\lambda_v\ne0.
\]

By \eqref{eq:ko9-ahss-order}, $ko_9(\mm M\eta)_{(2)}$ has order
$16$.  The Adams $E_2$-term just computed has four classes in stem
$9$, so strong convergence leaves no room for a nonzero Adams
differential.  Therefore, the $h_0$-multiple shows that any element detected
by $\lambda_v$ has order $8$.
It follows that
\[
ko_9(\mm M\eta)_{(2)}\cong\mathbb Z_8\oplus\mathbb Z_2.
\]
Together with \eqref{eq:ko9-three-primary}, this proves
\eqref{eq:ko9-two-primary}.  
\end{proof}

\section{Restricted Bazaikin bundles over the finite \(6\)-torsion model}
\label{sec:Y6-restriction}

The purpose of this section is to associate to each admissible
Bazaikin tuple \(\q\) a \(P\)-bundle over a fixed finite complex
\(Y\), and then to determine when two such bundles are
isomorphic.

We say a tuple $\q=(q_0,\ldots,q_5)$ a Bazaikin tuple if it satisfies $\sum_{i=0}^5q_i=0$ and (1.2) of \cite{FL}.

Let $s=\sigma_3(\q)/8$. By \cite{FL}, $\gcd(s,6)=1$. Set
\[ 
P=\SU(6)/\Sp(3),
\qquad
G=\SU(6)/\{\pm I\},
\qquad
H=S^1/\{\pm1\}\cong S^1.
\]
The group \(G\) acts naturally on \(P\).  For a Bazaikin
tuple 
$\q$,
the homomorphism
\[
\lambda_{\q}:H\longrightarrow G,
\qquad
\lambda_{\q}(\theta)
=
\left[
\diag\left(
e^{iq_0\theta/2},\ldots,e^{iq_5\theta/2}
\right)
\right]
\]
 induces an \(H\)-action on \(P\). Here
$\theta\in\mathbb R/2\pi\mathbb Z$ is the angular coordinate on the
effective circle $H$; changing $\theta$ by $2\pi$ multiplies the
diagonal matrix by $-I$, so the formula is well defined in $G$.  This is precisely the
free circle action whose quotient is the Bazaikin space
$
\mathcal B_{\q}=H\backslash P.$

Let
$
H\longrightarrow EH\longrightarrow BH$
be the universal principal \(H\)-bundle.  The above \(H\)-action $\lambda_\q $ on
\(P\) determines an associated \(P\)-bundle
\begin{equation}\label{eq:universal-Pq}
P
\longrightarrow
EH\times_{H}P
\longrightarrow
BH.
\end{equation}

We shall pull back this bundle to a fixed finite complex \(Y\).
The particular choice of \(Y\) will be made below; it is designed so
that its cohomology and complex \(K\)-theory retain the torsion
information needed to compare the bundles associated with different
tuples \(\q\).

\subsection{The finite model \(Y\) and the restricted \(P\)-bundle}
\quad

Let \(L\to\CP^5\) be the Hopf line bundle and write
\[
x=c_1(L)\in H^2(\CP^5;\Z).
\]
Consider the complex rank-four vector bundle
\[
V_{6,\mathbb C}
=
L^{\oplus3}\oplus L^{\otimes6},
\]
and let
\[
V_6=(V_{6,\mathbb C})_{\mathbb R}
\]
be its underlying oriented real rank-eight bundle.  Its Euler class and
second Stiefel--Whitney class are
\[
e(V_6)
=
c_4(V_{6,\mathbb C})
=
6x^4
\]
\[
w_2(V_6)
=
c_1(V_{6,\mathbb C})\bmod2
=
9\bar x
=
\bar x,
\]
where \(\bar x\) denotes the mod-\(2\) reduction of \(x\).

Let
$
\pi:S(V_6)\longrightarrow\CP^5$
be the unit sphere bundle.  Since the fibre is \(S^7\), the Gysin
sequence, together with \(e(V_6)=6x^4\), gives, through degree \(10\),
\[
H^{\le10}(S(V_6);\Z)
\cong
\left(\Z[x]/(6x^4,x^6)\right)^{\le10}.
\]
With \(\mathbb Z_2\)-coefficients the Euler class vanishes, so 
\[
H^*(S(V_6);\mathbb Z_2)
\cong
\mathbb Z_2[\bar x,\bar y_7]/(\bar x^6,\bar y_7^2),
\qquad
|\bar x|=2,\quad |\bar y_7|=7.
\]

Use the minimal CW model of $S(V_6)$ \cite{Wall1965}. This model has
one cell in each of the dimensions $7,9$ and $2j$, $0\le j\le5$.
Let
\[
Y=S(V_6)^{(10)}
\]
be the $10$-skeleton of this model. The inclusion induces an
isomorphism on integral cohomology in degrees $\le 10$. We use the same notation
$x$ and $\bar x$ for the restricted classes. Thus
\begin{equation}\label{eq:Y-cohomology}
H^*(Y;\Z)
\cong
\Z[x]/(6x^4,x^6),
\qquad |x|=2.
\end{equation}
The composite
\[
c_Y:
Y\hookrightarrow S(V_6)
\xrightarrow{\pi}
\CP^5
\hookrightarrow
\CP^\infty\simeq BH
\]
represents the generator $x\in H^2(Y;\Z)$.

We now consider the restriction of the $P$-bundle \eqref{eq:universal-Pq} to \(Y\): 
\begin{equation}\label{eq:Pq-restricted}
X_{\q}
=
c_Y^*(EH\times_H P)
\xrightarrow{\pi_Y}Y.
\end{equation}
Thus the problem is reduced to deciding when the bundles
\( X_{\q}\) and \( X_{\q'}\) over the fixed finite
complex \(Y\) are isomorphic.

We record two properties of \(Y\).  The first is a Steenrod operation
that will also be useful later in the paper.  The second is the
\(K\)-theory relation that is essential for the bundle comparison in
this section.

The tangent bundle of the sphere bundle satisfies
\[
TS(V_6)\oplus\epsilon^1
\cong
\pi^*\bigl(T\CP^5\oplus V_6\bigr).
\]
Since
$
w_2(T\CP^5)
=
0,$
we obtain
\[
w_2(TS(V_6))
=
\pi^*w_2(V_6)
=
\bar x.
\]
As \(S(V_6)\) is oriented, its second Wu class is
\[
v_2(S(V_6))=\bar x.
\]
Applying the Wu formula to
$
\bar x^4\bar y_7
\in H^{15}(S(V_6);\mathbb Z_2)$
gives
\[
\Sq^2(\bar x^4\bar y_7)=\bar x^5\bar y_7.
\]
On the other hand, the Cartan formula and
$
\Sq^2(\bar x^4)
=
0
$
give
\[
\bar x^4\Sq^2(\bar y_7)
=
\bar x^5\bar y_7.
\]
Since
$
H^9(S(V_6);\mathbb Z_2)
=
\mathbb Z_2\langle\bar x\bar y_7\rangle,$
we conclude that
\begin{equation}\label{eq:Y6-Sq2}
\Sq^2\bar y_7
=
\bar x\,\bar y_7.
\end{equation}

We next derive the \(K\)-theory relation used below.  For a complex
vector bundle \(E\), we use the $\lambda$-operation
\[
e_K(E)
=
\lambda_{-1}(E^*)
=
\Sigma_j(-1)^j[\Lambda^jE^*]
\]
where $E^\ast$ is the dual bundle of $E$. Note that $$\lambda_{-1}((L^{ n})^\ast)=1-(L^{n})^\ast=1-L^{-n}.$$
Since
$V_{6,\mathbb C}
=
L^{\oplus3}\oplus L^{\otimes 6},$
the multiplicativity of \(\lambda_{-1}\) gives
\[
e_K(V_{6,\mathbb C})
=
(1-L^{-1})^3(1-L^{-6}).
\]

This class vanishes after pullback to \(S(V_6)\).  Indeed,
\(\pi^*V_{6,\mathbb C}\) has the tautological nowhere-zero section
$
s(v)=(v,v).$
Hence
\[
\pi^*V_{6,\mathbb C}
\cong
\underline{\mathbb C}\oplus W
\]
for some complex rank-three bundle \(W\) and the trivial complex line bundle $\underline{\mb{C}}$.  Therefore
\[
\begin{aligned}
\pi^*e_K(V_{6,\mathbb C})
&=
\lambda_{-1}
\bigl((\pi^*V_{6,\mathbb C})^*\bigr)\\
&=
\lambda_{-1}(\underline{\mathbb C})\,
\lambda_{-1}(W^*)\\
&=0.
\end{aligned}
\]
In particular, the same class vanishes after restriction to \(Y\).

From now on, let \(L\) also denote the pullback of the Hopf line bundle
to \(Y\), and put
\[
t=[L]-1\in\widetilde K^0(Y).
\]
Since $K^0(\CP^5)=\Z[t]/t^6$ \cite{Sanderson1964}, we have
\[
t^6=0\in \widetilde K^0(Y).
\]
Modulo \(t^6\),
\[
1-L^{-1}
=
t-t^2+t^3-t^4+t^5
\]
and
\[
1-L^{-6}
=
6t-21t^2+56t^3-126t^4+252t^5.
\]
Thus, in $\widetilde K^0(Y)$,
\[
\begin{aligned}
0
&=
(1-L^{-1})^3(1-L^{-6})\\
&=
6t^4-39t^5
\pmod{t^6}.
\end{aligned}
\]
Multiplying by \(t\) gives \(6t^5=0\), and therefore
\[
12t^4
=
78t^5
=
13(6t^5)
=
0.
\]
We have proved
\begin{equation}\label{eq:Y6-t-orders}
\text{$12t^4=0$,
$6t^5=0$ in $\widetilde K^0(Y)$.}
\end{equation}
These two relations are the key \(K\)-theoretic input in the
isomorphism criterion below.

\subsection{Arithmetic input for the bundle comparison}
\quad

We now suppose that \(\q\) and \(\q'\) are Bazaikin tuples satisfying
\[
\sigma_2(\q)=\sigma_2(\q'),
\qquad
\sigma_3(\q)=\sigma_3(\q').
\]
Our goal in this subsection is to convert these equalities into the
divisibility conditions needed to annihilate the \(t^4\)- and
\(t^5\)-terms in the difference of the associated \(K\)-theory classes.

\begin{lemma}\label{lem:Y6-qmod4}
For every Bazaikin tuple \(\q\), five of the \(q_i\) lie in
one odd residue class modulo \(4\), while the remaining one lies in the
other.  

Let
$\varepsilon_{\q}\in\{1,-1\}$
be the residue class represented by the five entries.  If
$
\sigma_3(\q)=\sigma_3(\q'),$
then
$
\varepsilon_{\q}=\varepsilon_{\q'}.$
\end{lemma}

\begin{proof}
All entries are odd and $\sum_iq_i=0$ \cite[(1.2)]{FL}. If all six entries had
the same residue modulo $4$, their sum would be $2$ modulo $4$.
Thus both odd residue classes occur.

If there were at least two entries congruent to \(1\) and at least two
entries congruent to \(-1\) modulo \(4\), then two disjoint mixed pairs
would have sums divisible by \(4\).  This contradicts the admissibility
condition that the gcd of the corresponding two pair sums is \(2\) \cite[(1.2)]{FL}.

Put
$
\varepsilon=\varepsilon_{\q},$
$a_i=\frac{q_i-\varepsilon}{2}.$
Then exactly one of the \(a_i\) is odd and the other five are even.
We use $e_j(a)$ to denote the elementary symmetric functions for $a=(a_0,\cdots,a_5)$.

Since
$
q_i=2a_i+\varepsilon$ and $\sum_iq_i=0$, we have $e_1(a)=-3\varepsilon$. Moreover,
\[
\sigma_3(\q)
=
8e_3(a)+16\varepsilon e_2(a)-40\varepsilon.
\]
Since exactly one \(a_i\) is odd, \(e_2(a)\) is even and
\(e_3(a)\) is divisible by \(4\).  Hence
\[
\sigma_3(\q)
\equiv
-40\varepsilon
\pmod{32}.
\]
Therefore
\[
\varepsilon=1
\Longrightarrow
\sigma_3(\q)\equiv24\pmod{32},
\qquad
\varepsilon=-1
\Longrightarrow
\sigma_3(\q)\equiv8\pmod{32}.
\]
Thus \(\sigma_3(\q)\) determines \(\varepsilon_{\q}\).
\end{proof}

Hence, the two Bazaikin tuples have a common value
$
\varepsilon
=
\varepsilon_{\q}
=
\varepsilon_{\q'}.$
Set
\[
a_i=\frac{q_i-\varepsilon}{2},
\qquad
a_i'=\frac{q_i'-\varepsilon}{2},
\]
and write
\[
a=(a_0,\ldots,a_5),
\qquad
a'=(a_0',\ldots,a_5').
\]
Exactly one entry of each vector is odd.  Moreover,
$
e_1(a)
=
-3\varepsilon$.
The shift formulas give
\[
\sigma_2(\q)=4e_2(a)-15
\]
and
\[
\sigma_3(\q)
=
8e_3(a)+16\varepsilon e_2(a)-40\varepsilon.
\]
The same formulas hold for \(\q'\) and \(a'\). Since
\[
\sigma_2(\q)=\sigma_2(\q'),
\qquad
\sigma_3(\q)=\sigma_3(\q'),
\]
we obtain
\begin{equation}\label{eq:Y6-e123}
e_j(a)=e_j(a'),
\qquad
j=1,2,3.
\end{equation}

The next lemma supplies the \(3\)-primary part of the required
divisibility.

\begin{lemma}\label{lem:Y6-mod9}
Let \(\q,\q'\) be Bazaikin tuples satisfying
\[
\sigma_2(\q)=\sigma_2(\q'),
\qquad
\sigma_3(\q)=\sigma_3(\q').
\]
Then
\[
\sigma_4(\q)\equiv\sigma_4(\q')\pmod9,
\qquad
\sigma_5(\q)\equiv\sigma_5(\q')\pmod9.
\]
Consequently,
\[
e_4(a)\equiv e_4(a')\pmod9,
\qquad
e_5(a)\equiv e_5(a')\pmod9.
\]
\end{lemma}

\begin{proof}
By \cite[Corollary~3.2]{FL}, up to permutation of the entries and an
overall change of sign, the reduction of a Bazaikin tuple
modulo \(3\) is one of
\begin{equation}\label{eq:Y6-three-patterns}
1^6,\qquad
0^31^3,\qquad
0\,1^4(-1).
\end{equation}
Changing \(\q\) to \(-\q\) changes the elementary symmetric functions
according to
$
\sigma_k(-\q)=(-1)^k\sigma_k(\q).$
Thus it suffices to perform the mod-\(9\) calculation for the three
patterns in \eqref{eq:Y6-three-patterns}.

Choose
\[
r_i\in\{0,1,-1\},
\qquad
q_i=r_i+3u_i,
\quad
u_i\in\mathbb Z.
\]
Modulo \(9\), every term containing at least two factors \(3u_i\)
vanishes.  Hence
\begin{equation}\label{eq:Y6-mod9-expansion}
\sigma_k(\q)
\equiv
\sigma_k(r)
+
3\sum_{i=0}^5
u_i\,
\sigma_{k-1}(r_0,\ldots,\widehat{r_i},\ldots,r_5)
\pmod9.
\end{equation}
Then we have the following table:
\[
\begin{array}{c|c|c}
\q\bmod3
&
(\sigma_2,\sigma_3)\bmod9
&
(\sigma_4,\sigma_5)\bmod9
\\ \hline
1^6
&
(3,5)
&
(0,3)
\\[1mm]
0^31^3
&
(\alpha,\,4-\alpha),
\quad \alpha\in\{0,3,6\}
&
(\alpha+3,\,0)
\\[1mm]
0\,1^4(-1)
&
(\alpha,\beta),
\quad
\alpha\in\{2,5,8\},
\quad
\beta\in\{1,4,7\}
&
(3,\,\alpha+6).
\end{array}
\]
The three sign-reversed patterns are obtained by
\[
(\sigma_2,\sigma_3,\sigma_4,\sigma_5)
\longmapsto
(\sigma_2,-\sigma_3,\sigma_4,-\sigma_5).
\]
It follows that
\[
(\sigma_2,\sigma_3)\pmod9
\]
determines
\[
(\sigma_4,\sigma_5)\pmod9.
\]
Therefore
\[
\sigma_4(\q)\equiv\sigma_4(\q')\pmod9,
\qquad
\sigma_5(\q)\equiv\sigma_5(\q')\pmod9.
\]

We now translate these congruences for the elementary
symmetric functions of \(a\) and \(a'\).  Write
\[
\Delta\sigma_j
=
\sigma_j(\q)-\sigma_j(\q'),
\qquad
\Delta e_j
=
e_j(a)-e_j(a').
\]
For
$
q_i=2a_i+\varepsilon,$
the general shift formula is
\[
\sigma_k(\q)
=
\sum_{j=0}^k
2^j\varepsilon^{k-j}
\binom{6-j}{k-j}e_j(a).
\]
Subtracting the corresponding formula for \(\q'\), and using
\eqref{eq:Y6-e123}, gives
\[
\Delta\sigma_4
=
16\,\Delta e_4,
\]
\[
\Delta\sigma_5
=
32\varepsilon\,\Delta e_4
+
32\,\Delta e_5.
\]
Since
$
9\mid \Delta\sigma_4$,
$9\mid \Delta\sigma_5,$
and \(16,32\) are units modulo \(9\), we have
$
9\mid\Delta e_4$ and $9\mid\Delta e_5.$
This completes the proof.
\end{proof}

We next give some divisibilities matching the annihilation relations
\eqref{eq:Y6-t-orders}. 

\begin{lemma}\label{lem:Y6-e45}
With the notation above,
$
144\mid \Delta e_4,$
$144\mid\Delta e_5.$
Define
\[
B_r(a)=\sum_{i=0}^5\binom{a_i}{r}
\]
where $\binom{a_i}{r}$ is the generalized binomial coefficient. 
Then
\[
24\mid\bigl(B_4(a)-B_4(a')\bigr),
\qquad
6\mid\bigl(B_5(a)-B_5(a')\bigr).
\]
\end{lemma}

\begin{proof}
After reordering, write
\[
a_0\ \text{odd},
\qquad
a_j=2b_j,\quad1\le j\le5,
\]
and let
$
s_r=e_r(b_1,\ldots,b_5)$
 denote the \(r\)-th elementary symmetric function of
\(b_1,\ldots,b_5\).  We then have
\[
e_2(a)=2a_0s_1+4s_2,
\]
\[
e_3(a)=4a_0s_2+8s_3,
\]
\[
e_4(a)=8a_0s_3+16s_4.
\]
Since \(a_0\) is odd, reduction modulo \(2\) gives
\[
\frac{e_2(a)}2\equiv s_1,
\qquad
\frac{e_3(a)}4\equiv s_2,
\qquad
\frac{e_4(a)}8\equiv s_3
\pmod2.
\]

Let \(n\) be the number of odd entries among
\(b_1,\ldots,b_5\).  Modulo \(2\), \(b_j\) becomes \(1\)
or \(0\).  Hence a monomial in
\(s_r=e_r(b_1,\ldots,b_5)\) contributes \(1\) precisely when all its
\(r\) factors are chosen among the \(n\) odd entries.  Therefore
\[
s_r\equiv\binom{n}{r}\pmod2.
\]
Since
\[
\binom{n}{1}\binom{n}{2}
=
3\binom{n}{3}
+
2\binom{n}{2},
\]
we obtain
\[
s_1s_2
\equiv
s_3
\pmod2.
\]
Consequently,
\[
\frac{e_4(a)}8
\equiv
\frac{e_2(a)}2
\frac{e_3(a)}4
\pmod2,
\]
or equivalently,
\[
e_4(a)
\equiv
8
\left(\frac{e_2(a)}2\right)
\left(\frac{e_3(a)}4\right)
\pmod{16}.
\]
The same formula holds for \(a'\).  Since
$
e_2(a)=e_2(a')$ and $e_3(a)=e_3(a')$ (cf. \eqref{eq:Y6-e123}), it follows that
$
16\mid\Delta e_4.$

Moreover, every monomial occurring in \(e_5(a)\) contains at least
four of the five even entries \(a_1,\ldots,a_5\).  Hence
$
16\mid e_5(a).$
The same argument gives
$
16\mid e_5(a'),$
and therefore
$
16\mid\Delta e_5.$

Combining these divisibilities with Lemma~\ref{lem:Y6-mod9}, we obtain
\[
144\mid\Delta e_4,
\qquad
144\mid\Delta e_5.
\]

It remains to pass from the elementary symmetric functions
$e_j(a)$ to $B_4(a)$ and $B_5(a)$. Write
$\Delta B_r=B_r(a)-B_r(a')$. Put
\[
p_j(a)=\sum_{i=0}^5a_i^j,
\qquad
p_j(a')=\sum_{i=0}^5(a_i')^j,
\]
and write
\[
\Delta p_j=p_j(a)-p_j(a').
\]
Newton's identities give
\[
p_1(a)=e_1(a),
\]
\[
p_2(a)=e_1(a)p_1(a)-2e_2(a)
=
e_1^2(a)-2e_2(a),
\]
\[
p_3(a)=e_1(a)p_2(a)-e_2(a)p_1(a)+3e_3(a)
=
e_1^3(a)-3e_1(a)e_2(a)+3e_3(a).
\]
Since
$e_j(a)=e_j(a'),$ $j=1,2,3,$ (cf. \eqref{eq:Y6-e123}),
it follows that
\[
p_j(a)=p_j(a'),
\qquad j=1,2,3.
\]
Equivalently,
$
\Delta p_1=\Delta p_2=\Delta p_3=0.$

The fourth Newton identity is
\[
p_4(a)
=
e_1(a)p_3(a)-e_2(a)p_2(a)+e_3(a)p_1(a)-4e_4(a).
\]
Subtracting the identities for \(a\) and \(a'\), and using
\[
\Delta e_1=\Delta e_2=\Delta e_3
=
\Delta p_1=\Delta p_2=\Delta p_3=0,
\]
we obtain
\[
\Delta p_4=-4\Delta e_4.
\]

Similarly, the fifth Newton identity is
\[
p_5(a)
=
e_1(a)p_4(a)-e_2(a)p_3(a)+e_3(a)p_2(a)-e_4(a)p_1(a)+5e_5(a).
\]
Since \(e_1\) and \(p_1\) agree for \(a\) and \(a'\), and
\(p_1=e_1\), subtraction gives
\[
\begin{aligned}
\Delta p_5
&=
e_1\Delta p_4
-
p_1\Delta e_4
+
5\Delta e_5\\
&=
e_1(-4\Delta e_4)
-
e_1\Delta e_4
+
5\Delta e_5\\
&=
-5e_1\Delta e_4
+
5\Delta e_5.
\end{aligned}
\]

We next express \(B_4\) and \(B_5\) in terms of the power sums.  Since
\[
\binom{m}{4}
=
\frac{m^4-6m^3+11m^2-6m}{24},
\]
we have
\[
B_4(a)
=
\frac{
p_4(a)-6p_3(a)+11p_2(a)-6p_1(a)
}{24}.
\]
Thus
\[
\Delta B_4
=
\frac{
\Delta p_4-6\Delta p_3+11\Delta p_2-6\Delta p_1
}{24}
=
\frac{\Delta p_4}{24}
=
-\frac{\Delta e_4}{6}.
\]

Likewise,
\[
\binom{m}{5}
=
\frac{
m^5-10m^4+35m^3-50m^2+24m
}{120},
\]
so
\[
B_5(a)
=
\frac{
p_5(a)-10p_4(a)+35p_3(a)-50p_2(a)+24p_1(a)
}{120}.
\]
Therefore
\[
\Delta B_5
=
\frac{
\Delta p_5-10\Delta p_4
}{120}
=
\frac{
\Delta e_5+(8-e_1)\Delta e_4
}{24}.
\]

Finally,
$
144\mid\Delta e_4$
implies
$
24\mid\Delta B_4.$
 $
144\mid\Delta e_4$ and
$144\mid\Delta e_5$
imply 
\[
144\mid \Delta e_5+(8-e_1)\Delta e_4,
\]
 and hence
$
6\mid \Delta B_5
$.
\end{proof}

\subsection{The isomorphism criterion for the restricted bundles}
\quad

We can now compare the bundles defined in
\eqref{eq:Pq-restricted}.

\begin{proposition}\label{prop:Y-restriction-isomorphic}
Let \(\q,\q'\) be Bazaikin tuples satisfying
\[
\sigma_2(\q)=\sigma_2(\q'),
\qquad
\sigma_3(\q)=\sigma_3(\q').
\]
Then the corresponding restricted \(P\)-bundles over \(Y\) are
isomorphic:
\[
 X_{\q}\cong X_{\q'}.
\]
\end{proposition}

\begin{proof}
We first reduce the comparison of the \(P\)-bundles to a comparison of
rank-six complex vector bundles.  Set
\[
\widehat G
=
\frac{\SU(6)\times S^1}{\langle(-I,-1)\rangle}.
\]
There are exact sequences
\[
1\longrightarrow S^1
\longrightarrow\widehat G
\longrightarrow G
\longrightarrow1
\]
and
\[
1\longrightarrow\Z/3
\longrightarrow\widehat G
\stackrel{\rho}{\longrightarrow}U(6)
\longrightarrow1,
\qquad
\rho[A,z]=zA.
\]
The homomorphism \(\lambda_{\q}\) lifts to
\[
\widehat\lambda_{\q}(\theta)
=
\left[
\diag\left(
e^{iq_0\theta/2},\ldots,e^{iq_5\theta/2}
\right),
e^{-i\varepsilon\theta/2}
\right]\in \widehat G,
\]
and under \(\rho\) the associated rank-six complex vector bundle over
\(Y\) is
\[
E_{\q}
=
\bigoplus_{i=0}^5L^{a_i}.
\]
Similarly,
\[
E_{\q'}
=
\bigoplus_{i=0}^5L^{a_i'}.
\]

We next compare their classes in \(K^0(Y)\).  Recall that
\[
t=[L]-1,
\qquad
t^6=0.
\]
For any integer \(a_i\), the generalized binomial expansion gives
\[
(1+t)^{a_i}
=
\sum_{r=0}^5\binom{a_i}{r}t^r
\]
in \(K^0(Y)\).  Hence
\[
\begin{aligned}{}
[E_{\q}]
&=
\sum_{i=0}^5[L^{a_i}]
=
\sum_{i=0}^5(1+t)^{a_i}\\
&=
\sum_{r=0}^5B_r(a)t^r.
\end{aligned}
\]
The same formula holds for \(E_{\q'}\).  For \(r\le3\), the coefficient
\(B_r(a)\) is determined by \(e_1(a),e_2(a),e_3(a)\).  Explicitly,
\[
B_1(a)=e_1,
\]
\[
B_2(a)
=
\frac12\bigl(e_1^2-e_1-2e_2\bigr),
\]
and
\[
B_3(a)
=
\frac16\Bigl(
e_1^3-3e_1^2+2e_1
-3e_1e_2+6e_2+3e_3
\Bigr).
\]
Thus \eqref{eq:Y6-e123} gives
\[
B_r(a)=B_r(a'),
\qquad
r=1,2,3.
\]
Since \(B_0(a)=B_0(a')=6\), it follows that
\[
[E_{\q}]-[E_{\q'}]
=
\Delta B_4\,t^4+\Delta B_5\,t^5.
\]
By Lemma~\ref{lem:Y6-e45},
$
24\mid\Delta B_4,$
$6\mid\Delta B_5,$
whereas \eqref{eq:Y6-t-orders} gives
$
12t^4=0,$
$6t^5=0.$
Consequently,
\[
\Delta B_4\,t^4=0,
\qquad
\Delta B_5\,t^5=0,
\]
and therefore
$
[E_{\q}]=[E_{\q'}]$
in $K^0(Y).$

Both bundles have complex rank \(6\), while
$
\dim Y=10<12=2\cdot6.$
 Therefore,
$
E_{\q}\cong E_{\q'}.$
Thus the corresponding principal \(U(6)\)-bundles are isomorphic. Both \(U(6)\)-bundles admit \(\widehat G\)-lifts.  Since
$
H^1(Y;\Z_3)=0,$
such a lift through
\[
B\Z_3\longrightarrow B\widehat G\longrightarrow BU(6)
\]
is unique up to isomorphism.  Hence the principal \(\widehat G\)-bundles,
and therefore the induced principal \(G\)-bundles, are isomorphic.
Passing to the associated \(G\)-space \(P\) gives
$
 X_{\q}\cong X_{\q'}.$
\end{proof}


\subsection{The common source and target}
\label{subsec:common-H9-Thom}\quad

We next describe the cohomology of $X_\q$ and the maps needed to
compare the degree-nine classes of two Bazaikin spaces.
The pullback construction gives the homotopy-commutative diagram
\begin{equation}\label{eq:H9-circle-square}
\begin{tikzcd}[column sep=large,row sep=large]
P\arrow[r,equal]\arrow[d]
  &P\arrow[r,"\bar f"]\arrow[d,"\pi"']
  &\K(\Z,5)\arrow[d,"i"]\\
X_\q\arrow[r,"j"]\arrow[d,"\pi_Y"']
  &EH\times_H P\simeq\mathcal B_\q\arrow[r,"f"]\arrow[d]
  &B\arrow[d,"p"]\\
Y\arrow[r,"c_Y"]
  &\CP^\infty\arrow[r,equal]
  &\CP^\infty.
\end{tikzcd}
\end{equation}
Here the middle column is a fibration after replacing the Borel
construction by $\mathcal B_\q$. The equivalence with
$\mathcal B_\q$ uses the freeness of the effective $H$-action.

Choose $\alpha_9\in H^9(P)$ and the compatible generator
$\beta_9\in H^9(\mathcal B_\q)$ as in
\eqref{beta9top9}, so that
\begin{equation}\label{pi_inverse9}
\pi^*\beta_9=s\alpha_9.
\end{equation}
We continue to write $x$ for the pullback of the degree-two generator
to $X_\q$.

\subsubsection{Integral cohomology and Steenrod operations}

\begin{proposition}\label{prop:restricted-P-ring}
There are classes $y,z,w$ of degrees $7,9,11$, respectively, such
that the nonzero integral cohomology groups through degree $14$ are
\begingroup
\small
\setlength{\arraycolsep}{3pt}
\[
\begin{array}{c|ccccccccc}
k&0&2&4&6&7&9&11&13&14\\\hline
H^k(X_\q)
&\Z&\Z\langle x\rangle&\Z\langle x^2\rangle
&\Z_s\langle x^3\rangle&\Z\langle y\rangle
&\Z\langle xy,z\rangle&\Z\langle xz,w\rangle
&\begin{gathered}\Z\langle x^2z\rangle\\
\oplus\,\Z_6\langle xw\rangle\end{gathered}&0
\end{array}
\]
\endgroup
As a graded ring truncated above degree $14$,
\begin{equation}\label{cohomo_Xq}
	H^{\le14}(X_\q)
\cong
\left(
\Z[x,y,z,w]/(sx^3,\ x^4,\ x^2y-6w,\ 6xw)
\right)^{\le13}.
\end{equation}
Moreover, $z$ restricts to $\alpha_9$ on the fibre $P$.
\end{proposition}

\begin{proof}
Consider the integral cohomological Serre spectral sequence
\[
E_2^{p,q}=H^p(Y;H^q(P))
\Longrightarrow H^{p+q}(X_\q).
\]
 Recall that $H^*(P)$ is the exterior algebra on
$\alpha_5$ and $\alpha_9$.

For the fibration for $\mathcal B_\q $ over $\CP^\infty$, the transgression is
$d_6(\alpha_5)=sx^3$; see \cite[Proposition~2.1]{FL}.
Naturality and multiplicativity give
\[
d_6(x^j\alpha_5)=sx^{j+3}
\]
for the pullback fibration. Thus
\[
E_7^{0,5}=0,\qquad E_7^{6,0}\cong\Z_s.
\]
For $j=1,2$, the target is $\Z_6$, and multiplication by $s$ is
surjective because $\gcd(s,6)=1$. Hence
\[
E_7^{2,5}=6\Z\langle x\alpha_5\rangle,\qquad
E_7^{4,5}=6\Z\langle x^2\alpha_5\rangle,
\qquad E_7^{8,0}=E_7^{10,0}=0.
\]
The remaining terms in the $q=5$ row relevant here are
\[
E_\infty^{6,5}=\Z\langle x^3\alpha_5\rangle,
\qquad
E_\infty^{8,5}=\Z_6\langle x^4\alpha_5\rangle.
\]
They have no possible nonzero incoming or outgoing differential.

The only potentially nonzero differential from the $q=9$ row is
\[
d_{10}:E_{10}^{0,9}\longrightarrow E_{10}^{10,0}.
\]
Its target is already zero. Thus $\alpha_9,x\alpha_9,x^2\alpha_9$
survive, giving
\[
E_\infty^{0,9}\cong\Z,\qquad
E_\infty^{2,9}\cong\Z,\qquad
E_\infty^{4,9}\cong\Z.
\]
These are all the contributions in total degrees at most $13$.
The additive extensions split because their quotients are free.

Choose $y$ detected by $6x\alpha_5$. Then $xy$ generates the
filtration-four subgroup of $H^9(X_\q)$. Choose $w$ detected by $x^3\alpha_5$ in degree $11$. Note that $x^2y$ is detected by
$6x^3\alpha_5$. Since there is no nonzero filtration above $6$ in
total degree $11$, this gives the exact relation
\[
x^2y=6w.
\]
Likewise, $xw$ is detected by $x^4\alpha_5$  generates the filtration-eight subgroup in degree
$13$; there is no higher integral filtration in this degree, so
$xw$ has order exactly $6$.

Finally, choose $z$ restricting to $\alpha_9$. Its multiples $xz$
and $x^2z$ lift generators of the free filtration quotients in
degrees $11$ and $13$. This yields the displayed bases and all the
ring relations through degree $13$.

 The only possible term in total degree
$14$ is killed by the injective transgression
$d_6(\alpha_5\alpha_9)=s x^3\alpha_9$. Hence $H^{14}(X_\q)=0$.
\end{proof}

The universal coefficient theorem now gives
\begin{equation}\label{eq:Y6-X-homology-low}
\begin{gathered}
H_4(X_\q)\cong\Z,\quad H_5(X_\q)\cong\Z_s,
\quad H_7(X_\q)\cong\Z,\\
H_8(X_\q)=H_{10}(X_\q)=0,\quad
H_9(X_\q)\cong\Z^2,\\
H_{11}(X_\q)\cong\Z^2,\quad
H_{12}(X_\q)\cong\Z_6,\quad H_{13}(X_\q)\cong\Z.
\end{gathered}
\end{equation}

For an integral cohomology class $a$, write $\bar a$ for its
reduction modulo $2$.

\begin{proposition}\label{prop:restricted-P-squares}
The $\Sq^2$-operations with target degree at most $13$ are
\begin{equation}\label{eq:restricted-P-Sq2}
\begin{aligned}
\Sq^2\bar x&=\bar x^2,&
\Sq^2\bar y&=\bar x\bar y,\\
\Sq^2(\bar x\bar z)&=\bar x^2\bar z,&
\Sq^2\bar w&=\epsilon\,\bar x\bar w,
\end{aligned}
\end{equation}
for some $\epsilon\in\Z_2$; all the other $\Sq^2$-operations in
this range vanish. Moreover,
\begin{equation}\label{eq:restricted-P-Sq4}
\Sq^4:H^k(X_\q;\Z_2)\longrightarrow H^{k+4}(X_\q;\Z_2)
\end{equation}
is zero for $k+4\le13$.
\end{proposition}

\begin{proof}
The mod-$2$ Serre spectral sequence shows that
\[
\pi_Y^*(\bar y_7)=\bar y.
\]
 By \eqref{eq:Y6-Sq2} and naturality,
\[
\Sq^2\bar y=\bar x\bar y.
\]
The Cartan formula gives $\Sq^2(\bar x\bar y)=0$.
Since $\dim Y=10$, both $\bar y$ and $\bar x\bar y$ have zero
fourth square, being pulled back from $Y$.

The kernel of restriction $H^9(X_\q)\to H^9(P)$ is
generated by $xy$. Hence \eqref{pi_inverse9} gives an integer $t$
such that
\begin{equation}\label{Pback_Prim_degree9}
j^*\beta_9=sz+txy.
\end{equation}
Since $s$ is odd,
$
\bar z-j^*\bar\beta_9\in\Z_2\langle\bar x\bar y\rangle.$
By \eqref{SqBq}, $\Sq^2\bar\beta_9=\Sq^4\bar\beta_9=0$. By $\Sq^2(\bar x\bar y)=\Sq^4(\bar x\bar y)=0$, 
it follows that $\Sq^2\bar z=\Sq^4\bar z=0$, and hence
\[
\Sq^2(\bar x\bar z)=\bar x^2\bar z.
\]

Finally, the reduction $\bar w$ is provided by $E^{6,5}_\infty$. By \cite[formulas~(2.25) and~(2.75)]{Si}, $\Sq^2\bar w$ admits Serre filtration $\ge 6$. Hence $\Sq^2\bar w$ is provided by $E^{p,13-p}_\infty$ for some $p\ge 6$. 
 The only possible such term is
$E_\infty^{8,5}=\Z_2\langle \bar x^4\bar \alpha_5\rangle$, represented by
$\bar x\bar w$. Therefore,
$\Sq^2\bar w=\epsilon\bar x\bar w,$ $\epsilon\in\Z_2.$

The remaining assertions follow from the Cartan formula and \eqref{cohomo_Xq}.
\end{proof}

\begin{proposition}\label{prop:restricted-P-P1}
Let $a_3$ denote reduction modulo $3$. The first reduced power,
with target degree at most $13$, satisfies
\[
\begin{aligned}
\mathcal P^1(x_3)&=0,& \mathcal P^1((x^2)_3)&=0,\\
\mathcal P^1(y_3)&=0,& \mathcal P^1((xy)_3)&=0,\\
\mathcal P^1(z_3)&=k_3(x^2z)_3,
\end{aligned}
\]
where $p_1(\mathcal B_\q)=kx^2$ and $k_3=k\bmod3$.
\end{proposition}

\begin{proof}
Since $\gcd(s,3)=1$, the relation $sx^3=0$ implies $x_3^3=0$.
Thus $\mathcal P^1(x_3)=x_3^3=0$, and the Cartan formula gives
$\mathcal P^1((x^2)_3)=0$.

The mod-$3$ Serre spectral sequence gives an isomorphism
\[
\pi_Y^*:H^7(Y;\Z_3)\xrightarrow{\cong}H^7(X_\q;\Z_3).
\]
Consequently $y_3$ comes from the base, and
$\mathcal P^1(y_3)=0$ because $\dim Y=10$. The Cartan formula
then gives $\mathcal P^1((xy)_3)=0$.

For the remaining class, we first recall 
$
\mathcal P^1((\beta_9)_3)=k_3(x^2\beta_9)_3$ (cf. \eqref{P1Bq}).
Using \eqref{Pback_Prim_degree9}, naturality now gives
\[
\begin{aligned}
s\mathcal P^1(z_3)
&=\mathcal P^1\bigl(j^*(\beta_9)_3\bigr)\\
&=j^*\mathcal P^1((\beta_9)_3)\\
&=k_3x_3^2\bigl(sz_3+t(xy)_3\bigr)\\
&=sk_3(x^2z)_3,
\end{aligned}
\]
since $x_3^3=0$. Cancelling the unit $s\bmod3$ proves the claim.
\end{proof}

\subsubsection{The common maps and the integral basis}

From now on, assume
\[
\sigma_2(\q)=\sigma_2(\q'),\qquad
\sigma_3(\q)=\sigma_3(\q')>0.
\]
The two Bazaikin spaces have the common space $B$ and bundle $\eta$
from Section~\ref{Normalbg}. Fix the isomorphism
$X_\q\cong X_{\q'}$ from
Proposition~\ref{prop:Y-restriction-isomorphic}, and denote the
common total space by $X=X_\q $.
Let
\[
j:X\longrightarrow\mathcal B_\q,
\qquad j':X\longrightarrow\mathcal B_{\q'}
\]
be the induced maps, and write $f,f'$ for the maps to $B$.

\begin{lemma}\label{fjfjprime}
The maps $f\circ j$ and $f'\circ j'$ are homotopic.
\end{lemma}

\begin{proof}
The bundle isomorphism covers the identity of $Y$, so
\[
p\circ f\circ j\simeq p\circ f'\circ j'.
\]
Since the homotopy fiber of $p:B\to \CP^\infty$ is $\K(\Z,5)$, $H^5(X)=0$ implies that the lifts are homotopic.
\end{proof}

By the choices of $z$ (cf. Proposition \ref{prop:restricted-P-ring}) and $\beta_9$ (cf. \eqref{pi_inverse9}), we have
\begin{equation}\label{cohomology_z_xy_beta}
j^*\beta_9=sz+txy,
\qquad
(j')^*\beta_9'=sz+t'xy.
\end{equation}
Let $z^\vee,(xy)^\vee$ be the integral homology basis
Kronecker-dual to $z,xy$. Then
\begin{equation}\label{homology_z_xy_beta}
\begin{aligned}
j_*z^\vee&=s\beta_9^\vee,&
j_*(xy)^\vee&=t\beta_9^\vee,\\
j'_*z^\vee&=s(\beta_9')^\vee,&
j'_*(xy)^\vee&=t'(\beta_9')^\vee.
\end{aligned}
\end{equation}

\begin{remark}\label{rem:Y6-primitive-beta9}
We verify that $j^*\beta_9$ is primitive. This also proves
$\gcd(s,t)=1$.
Let $Q_Y=c_Y^*EH\to Y$ be the principal circle bundle with Euler
class $x$. Its Gysin sequence and \eqref{eq:Y-cohomology} give
\[
\begin{gathered}
H^i(Q_Y)=0\quad(1\le i\le6),\qquad
H^7(Q_Y)\cong\Z,\\
H^8(Q_Y)=H^9(Q_Y)=0.
\end{gathered}
\]

Since the action of $H$ on $P$ is free, the description
$X_\q  =Q_Y\times_H P$ identifies $j$ with a fibration
\[
Q_Y\longrightarrow X_\q \xrightarrow{j}\mathcal B_\q.
\]
Its Serre spectral sequence has exactly two nonzero terms in total
degree $9$, namely $E_2^{9,0}\cong\Z$ and
$E_2^{2,7}\cong\Z$, which collapse at $E_2$. Thus
\[
0\longrightarrow H^9(\mathcal B_\q;\Z)
\xrightarrow{j^*}H^9(\mathcal P;\Z)
\longrightarrow\Z\longrightarrow0.
\]
So $j^*\beta_9$ is primitive. The same proof
gives $\gcd(s,t')=1$.
\end{remark}

Our remaining goal is the following congruence.
\begin{proposition}\label{ttprime_differential}
With these choices of integral generators \eqref{homology_z_xy_beta},
\[
t\equiv t'\pmod{24}.
\]
\end{proposition}

\subsubsection{The degree-nine edge homomorphisms}

By Lemma~\ref{fjfjprime}, the two pullback bundles
$(f\circ j)^*\eta$ and $(f'\circ j')^*\eta$ are isomorphic.
Set
\[
\varphi=f^*\eta,\qquad \varphi'=(f')^*\eta,
\qquad \omega=j^*\varphi\cong(j')^*\varphi'.
\]
The characteristic classes in \eqref{eq:eta-char} give
\begin{equation}\label{char_varphi_omega}
\begin{aligned}
w_2(\varphi)&=\bar x,&w_4(\varphi)&=0,\\
w_2(\omega)&=\bar x,&w_4(\omega)&=0.
\end{aligned}
\end{equation}
By
Lemma~\ref{fjfjprime}, the relavent normalized Thom spectra then fit into
a homotopy-commutative square
\begin{equation}\label{corr:Thom-square}
\begin{tikzcd}[column sep=large,row sep=large]
\mm M\omega\arrow[r,"j"]\arrow[d,"j'"']
  &\mm M\varphi\arrow[d,"f"]\\
\mm M\varphi'\arrow[r,"f'"']&\mm M\eta.
\end{tikzcd}
\end{equation}
We use the same symbols for these Thom maps and their base maps.

Let $\chi:ko\to H\Z$ be the zeroth Postnikov section. 

\begin{proposition}\label{prop:ko9-Thom-groups}
The induced homomorphisms
\[
\chi_\omega:ko_9(\mm M\omega)\longrightarrow H_9(\mm M\omega),
\qquad
\chi_\varphi:ko_9(\mm M\varphi)\longrightarrow H_9(\mm M\varphi)
\]
 are surjective. Moreover, there is a split exact sequence
\[
0\longrightarrow\Z_s\longrightarrow ko_9(\mm M\varphi)
\xrightarrow{\chi_\varphi}\Z\longrightarrow0.
\]
\end{proposition}

\begin{proof}
We use the AHSS
$$E^2_{p,q}(ko,X)=H_p(X;\pi_q(ko))\Longrightarrow ko_{p+q}(X).$$
By Lemma \ref{chi}, the homomorphism induced by $\chi$ is identified with the edge homomorphism of the AHSS.

\medskip\noindent
\textbf{Step 1: the edge term for $\mm M\omega$.}

The Thom isomorphism and \eqref{eq:Y6-X-homology-low} give
$E^2_{9,0}\cong\Z^2$. By
\eqref{eq:restricted-P-Sq2} and \eqref{char_varphi_omega},
\[
\Sq^2(\bar yU)
=(\Sq^2\bar y+\bar x\bar y)U=0.
\]
Thus the
$d_2$ from $(9,0)$ is zero.
The $d_3$-target vanishes by $H_6(X;\Z_2)=0$;
the targets for $r=4,6,7,8$ have zero coefficient groups.

The remaining possibilities are
\[
d_5:E^5_{9,0}\longrightarrow E^5_{4,4},\qquad
d_9:E^9_{9,0}\longrightarrow E^9_{0,8}.
\]
Here $E^5_{4,4}\cong E^9_{0,8}\cong\Z$: their earlier possible incoming sources are
zero or finite, and their earlier outgoing targets vanish. 

The rational splitting
\[
ko_{\mathbb Q}\simeq\bigvee_{j\ge0}\Sigma^{4j}H\mathbb Q
\]
implies that the rational AHSS degenerates. Hence the rational differentials $d_5,$ $d_9$ vanish. Thus both vanish
integrally by $E^5_{4,4}\cong E^9_{0,8}\cong\Z$. Therefore
\begin{equation}\label{eq:ko9-omega-E90}
E^\infty_{9,0}(ko,\mm M\omega)
=E^2_{9,0}(ko,\mm M\omega)\cong\Z^2,
\end{equation}
which proves surjectivity of $\chi_\omega$.

\medskip\noindent
\textbf{Step 2: the group $ko_9(\mm M\varphi)$.}

The only nonzero terms in total degree $9$ are
\[
\begin{array}{c|ccc}
(p,q)&(0,9)&(5,4)&(9,0)\\\hline
E^2_{p,q}(ko,\mm M\varphi)&\Z_2&\Z_s&\Z.
\end{array}
\]
Since $\Sq^2U=\bar\zeta U$, the differential
\[
d_2:E^2_{2,8}\longrightarrow E^2_{0,9}
\]
is surjective. Thus $E^\infty_{0,9}=0$.

For $(5,4)$, all incoming sources are zero except for
the $2$-primary $d_4$-source; its map to $\Z_s$ is zero. The only possible nontrivial outgoing differential is the
$d_5$ whose target is the free group $E^5_{0,8}=\Z$.
Consequently $E^\infty_{5,4}\cong\Z_s$.

For $(9,0)$, the targets of $d_2,d_3,d_4,d_6,d_7,d_8$ vanish.
The remaining $d_5,d_9$ vanish by the same torsion-free-target
and rational-degeneration argument as in Step~1. Hence
$
E^\infty_{9,0}
=E^2_{9,0}\cong\Z.$
This proves surjectivity of $\chi_\varphi$. 

The kernel of $\chi_\varphi$ is the
filtration-five subgroup $\Z_s$. The resulting extension splits
because the quotient is free.
\end{proof}

\subsubsection{Comparison of the two coefficients}

We now complete the proof of Proposition~\ref{ttprime_differential}.
The calculation used here is the one in \eqref{eq:ko9-two-primary}:
\[
ko_9(\mm M\eta)_{(2)}\cong\Z_8\oplus\Z_2,
\qquad
ko_9(\mm M\eta)_{(3)}\cong\Z_3.
\]
We also use its detection statement: every element detected by the
Adams filtration-zero functional $\lambda_v$ dual to $vU$ has
order $8$ on the $2$-primary part. This statement is what detects the
order-eight component. A generator of $ko_9(\M\eta)_{(3)}$ is detected by a class in $H_9(\M\eta;\Z)$ whose mod-3 reduction is dual to $z_9U\in H^9(\M\eta;\Z_3)$.

\begin{proof}[Proof of Proposition~\ref{ttprime_differential}]
Fix a projection
\[
\operatorname{pr}_{24}:ko_9(\mm M\eta)\longrightarrow\Z_{24}
\]
by identifying
$\Z_8\oplus\Z_3\cong\Z_{24}$, and killing all other summands. 
The projection is fixed once and used for both Bazaikin spaces.

Applying $ko_9$ to \eqref{corr:Thom-square} gives the commutative
diagram
\begin{equation}\label{ko_commuta}
\begin{tikzcd}[column sep=large,row sep=large]
ko_9(\mm M\omega)\arrow[r,"j_*"]\arrow[d,"j'_*"']
 &ko_9(\mm M\varphi)\arrow[d,"f_*"]\\
ko_9(\mm M\varphi')\arrow[r,"f'_*"']
 &ko_9(\mm M\eta)\arrow[r,"\operatorname{pr}_{24}"]&\Z_{24}.
\end{tikzcd}
\end{equation}
By
Proposition~\ref{prop:ko9-Thom-groups}, there exist edge lifts
\[
A_z,A_{xy}\in ko_9(\mm M\omega),\qquad
\Gamma\in ko_9(\mm M\varphi),\qquad
\Gamma'\in ko_9(\mm M\varphi')
\]
with edge images $z^\vee U,(xy)^\vee U,\beta_9^\vee U$ and
$(\beta_9')^\vee U$, respectively. 

We first check that
\[
C=\operatorname{pr}_{24}f_*(\Gamma)
\]
is a generator of $\Z_{24}$. Modulo $2$, we have
$f^*(u)=0$ and $f^*(v)=\bar\beta_9$ by \eqref{fastvbeta9}
and $H^7(\mathcal B_\q;\Z_2)=0$.
Naturality of the Adams filtration-zero identification therefore
shows that $f_*(\Gamma)_{(2)}$ is detected by $\lambda_v$.
It has order $8$ by Lemma~\ref{ko9Meta2primary}.

Modulo $3$, we have $f^\ast(z_9)= \mathfrak m(\beta_9)_3$ for some $0\ne \mathfrak m\in \Z_3$ by \eqref{fastz_9_beta_9}.
 Thus $f_*(\Gamma)$ has nonzero $3$-primary
edge image, and consequently nonzero $3$-primary component.
This proves that $C$ is a generator of $\Z_{24}\subset ko_9(\M\eta)$.

By naturality of the edge maps and \eqref{homology_z_xy_beta},
\[
j_*(A_z)=s\Gamma+\tau_z,
\qquad
j'_*(A_z)=s\Gamma'+\tau'_z,
\]
where $s\tau_z=s\tau'_z=0$. These torsion corrections vanish
after projection to $\Z_{24}$, since $\gcd(s,24)=1$.
Applying \eqref{ko_commuta} to $A_z$ therefore gives
\[
s\operatorname{pr}_{24}f_*(\Gamma)
=s\operatorname{pr}_{24}f'_*(\Gamma').
\]
Multiplication by $s$ is invertible in $\Z_{24}$, so
\[
\operatorname{pr}_{24}f_*(\Gamma)
=\operatorname{pr}_{24}f'_*(\Gamma')=C.
\]

Apply the same argument to $A_{xy}$. We have
\[
j_*(A_{xy})=t\Gamma+\tau_{xy},
\qquad
j'_*(A_{xy})=t'\Gamma'+\tau'_{xy},
\]
with $s\tau_{xy}=s\tau'_{xy}=0$. Commutativity and projection
now give
\[
tC=t'C\quad\text{in }\Z_{24}.
\]
Since $C$ has order $24$, this is exactly
$t\equiv t'\pmod{24}$.
\end{proof}


\section{Detecting the normal $5$-smoothing}\label{Detect}

The normal $5$-smoothing \eqref{eq:normal-smoothing} of a Bazaikin
space $\mathcal B_{\q}$ determines a bordism class
\[
[\mathcal B_{\q},(g,f)]
\in\Omega_{13}^{\mm O\langle8\rangle}(\eta).
\]
We first represent this class in stable homotopy and then compare
its image in $ko$-homology. As before, $s>0$ and $\gcd(s,6)=1$.

Throughout this section, all Thom spectra are normalized to virtual
rank zero. We write $U$ for a cohomological Thom class and $aU_*$
for the image of a homology class $a$ under the Thom isomorphism.
The same symbols $f,j,f',j'$ denote the maps of base spaces and the
induced Thom spectrum maps. A map on a smash product is understood
to be the identity on the coefficient spectrum unless otherwise stated.

By the Thom--Pontryagin construction, the bordism class corresponds
to the stable homotopy class of the composite
\begin{equation}\label{homotopy_repre}
S^{13}\xrightarrow{\ell}\mm M_{\mathcal B_{\q}}
\longrightarrow\mm M_{\mathscr B}
=\mm{MO}\langle8\rangle\wedge\mm M\eta,
\end{equation}
where $\mm M_{\mathcal B_{\q}}$ is the Thom spectrum of the stable
normal bundle of $\mathcal B_{\q}$. The ordinary Hurewicz image of
$\ell$ is $[\mathcal B_{\q}]U_*$.

The normal map factors as
\[
\mathcal B_{\q}\xrightarrow{\Delta}
\mathcal B_{\q}\times\mathcal B_{\q}
\xrightarrow{g\times f}\mm{BO}\langle8\rangle\times B=\mathscr B.
\]
Put $\varphi=f^*\eta$. The corresponding factorization of Thom maps is
\[
S^{13}\xrightarrow{\ell}\mm M_{\mathcal B_{\q}}
\longrightarrow\mm{MO}\langle8\rangle\wedge\mm M\varphi
\xrightarrow{1\wedge f}\mm{MO}\langle8\rangle\wedge\mm M\eta.
\]
Let $\varsigma_{\q}$ be the class of the composite through
$\mm{MO}\langle8\rangle\wedge\mm M\varphi$. Thus
\[
\varsigma_{\q}\in\pi_{13}
(\mm{MO}\langle8\rangle\wedge\mm M\varphi),
\qquad
[\mathcal B_{\q},(g,f)]\longleftrightarrow f_*(\varsigma_{\q}).
\]
Applying $\alpha_8:\mm{MO}\langle8\rangle\to ko$, we define
\[
\rho_{\q}=(\alpha_8)_*(\varsigma_{\q})
\in ko_{13}(\mm M\varphi).
\]
Naturality gives the commutative diagram
\[
\begin{tikzcd}[column sep=large,row sep=large]
\pi_{13}(\mm{MO}\langle8\rangle\wedge\mm M\varphi)
\arrow[r,"f_*"]\arrow[d,"(\alpha_8)_*"']
&\pi_{13}(\mm{MO}\langle8\rangle\wedge\mm M\eta)
\arrow[d,"(\alpha_8)_*"]\\
ko_{13}(\mm M\varphi)\arrow[r,"f_*"]
&ko_{13}(\mm M\eta).
\end{tikzcd}
\]

We retain $F_{p,q}$ and $\bar F_{p,q}$ for the AH filtrations of
$\Omega_{p+q}^{\mm O\langle8\rangle}(\eta)$ and
$ko_{p+q}(\mm M\eta)$, respectively. By Theorem~\ref{pFZ4},
\[
(\alpha_8)_*:
\Omega_{13}^{\mm O\langle8\rangle}(\eta)/F_{8,5}
\longrightarrow ko_{13}(\mm M\eta)/\bar F_{8,5}
\]
is injective. We will use the following immediate consequence.

\begin{lemma}\label{common_image_varsigma}
If
$f_*(\rho_{\q})=f'_*(\rho_{\q'})
\in ko_{13}(\mm M\eta)/\bar F_{8,5},$
then
\[
f_*(\varsigma_{\q})=f'_*(\varsigma_{\q'})
\in \Omega_{13}^{\mm O\langle8\rangle}(\eta)/F_{8,5}.
\]
The same implication holds after localization at any prime.
\end{lemma}

\begin{proof}
Apply the injectivity above to
$f_*(\varsigma_{\q})-f'_*(\varsigma_{\q'})$ and use the commutative
diagram. Localization preserves injectivity.
\end{proof}

We therefore compare $f_*(\rho_{\q})$ and $f'_*(\rho_{\q'})$ when
$\sigma_2(\q)=\sigma_2(\q')$ and
$\sigma_3(\q)=\sigma_3(\q')>0$.
The comparison at primes not dividing $s$ will be carried out below;
the primes dividing $s$ are handled by the linking invariant.

\subsection{The $ko$-homology of $\mm M\omega$ and $\mm M\varphi$}
\label{subsec:ko13-Momega-Mphi}\quad

For a connective coefficient spectrum $E$, we use the convention
\[
E^2_{p,q}(E,X)=H_p(X;\pi_q(E))\Longrightarrow E_{p+q}(X),
\qquad d_r:E^r_{p,q}\longrightarrow E^r_{p-r,q+r-1}.
\]
Write $\mm F_{p,q}$ and $\mb F_{p,q}$ for the AH filtrations of
$ko_{p+q}(\mm M\varphi)$ and $ko_{p+q}(\mm M\omega)$, respectively;
$\mm F_{p,q}'$ denotes the AH filtration for $\varphi'$.

For $\mm M\varphi$, the nonzero $E^2$-terms in total degree $13$ are
\[
\begin{array}{c|ccccc}
(p,q)&(13,0)&(11,2)&(9,4)&(5,8)&(4,9)\\\hline
E^2_{p,q}&\Z&\Z_2&\Z&\Z_s&\Z_2.
\end{array}
\]

\begin{lemma}\label{E130koMvarphi}
In the AHSS for $ko_*(\mm M\varphi)$,
\[
E^\infty_{13,0}(ko,\mm M\varphi)
=E^2_{13,0}(ko,\mm M\varphi)\cong\Z.
\]
\end{lemma}

\begin{proof}
Let $\chi:ko\to H\Z$ be the zeroth Postnikov section.
The composite $\chi\circ\alpha_8$ is the ordinary integral Thom
orientation. The Thom--Pontryagin construction therefore gives
\[
\chi_*(\rho_{\q})
=[\mathcal B_{\q}]U_*\in H_{13}(\mm M\varphi).
\]
This is a generator, so the edge homomorphism is surjective.
Thus,
the asserted equality follows.
\end{proof}

\begin{lemma}\label{koM4}
In the AHSS for $ko_*(\mm M\varphi)$,
\[
E^\infty_{11,2}(ko,\mm M\varphi)
=
E^2_{11,2}(ko,\mm M\varphi)
\cong\Z_2,
\]
and
\[
E^\infty_{9,4}(ko,\mm M\varphi)
=
E^2_{9,4}(ko,\mm M\varphi)
\cong\Z .
\]
Moreover, the extension
\[
0\longrightarrow
E^\infty_{9,4}(ko,\mm M\varphi)
\longrightarrow
\mm{ F}_{11,2}/\mm{F}_{8,5}
\longrightarrow
E^\infty_{11,2}(ko,\mm M\varphi)
\longrightarrow0
\]
is nontrivial. Consequently,
$
\mm{ F}_{11,2}/\mm{F}_{8,5}\cong\Z .$
\end{lemma}

\begin{proof}
For the term
$E^2_{9,4}\cong\Z,$
the only possible differentials are the outgoing differentials $d_5,d_7,d_9$. After
rationalization,
\[
ko_{\mathbb Q}
\simeq
\bigvee_{i\ge0}\Sigma^{4i}H\mathbb Q,
\]
so the rational AHSS collapses. Hence $d_5,d_9$  with
target a free abelian group vanish integrally.

It remains to exclude
$
d_7:E^7_{9,4}(ko,\mm M\varphi)
\longrightarrow E^7_{2,10}(ko,\mm M\varphi).$

Recall that $\varphi= c^*\xi$, where
$\xi=-kL$ is the bundle over $\CP^\infty$, the map $c:\mathcal B_\q\to \CP^\infty$ classifies the principal
circle bundle (cf. Section \ref{oriented}).
The induced Thom map
$
\mm M\varphi\longrightarrow\mm M\xi$
gives a morphism of AHSSs.
In particular, the induced map
\[
E^2_{2,10}(ko,\mm M\varphi)
\longrightarrow E^2_{2,10}(ko,\mm M\xi)
\]
is an isomorphism.

For $E^r_{2,10}(ko,\mm M\xi)$, all outgoing differentials are zero. The possible incoming $d_2$ vanishes
because
\[
\Sq^2(\bar\zeta U)
=(\bar\zeta^2+\bar\zeta^2)U=0,
\]
and the incoming $d_3$ vanishes since its source
has odd order.  All other incoming differentials
before page $10$ have zero source. Since 
$$\Sq^2(\bar x^6 U)
=\bar x^6\, w_2(\xi)U=\bar x^7U$$ 
we obtain $E^3_{12,1}(ko,\M\xi)=0$. This implies that
$$d_{10}:E^{10}_{12,1}(ko,\M\xi)\to E^{10}_{2,10}(ko,\M\xi)$$
is zero. By $E^{2}_{13,0}(ko,\M\xi)=0$, we conclude that 
\[
E^\infty_{2,10}(ko,\mm M\xi)
=
E^2_{2,10}(ko,\mm M\xi).
\]
Thus, naturality implies
\begin{equation}\label{Einfty_210_Mxi}
	E^{\infty}_{2,10}(ko,\M\varphi)=E^{2}_{2,10}(ko,\M\varphi).
\end{equation} 
Then $d_7:E^7_{9,4}(ko,\mm M\varphi)\to E^7_{2,10}(ko,\mm M\varphi)$ is zero. Therefore, 
\[
E^\infty_{9,4}(ko,\mm M\varphi)
=
E^2_{9,4}(ko,\mm M\varphi)
\cong\Z.
\]

By Propositions~\ref{Mvarphih12b9zeta} and
\ref{E2112h12beta9zeta}, the coset
$
\{h_1^2\beta_9 x U\}
\subset
ko_{13}(\mm M\varphi)$
 has AH filtration $11$.   
Under the natural map
$f:\mm M\varphi\longrightarrow\mm M\eta,$ we have
$$\{h_1^2\beta_9x U\}\longmapsto \{h_1^2vxU\}.$$
The latter has AH filtration $11$ (cf. Corollary \ref{AHSS_Meta_112}), and its double is provided by
$
E^\infty_{9,4}(ko,\mm M\eta)$ (cf. Lemma \ref{h0module}, \ref{h0h2_survives} and Theorem \ref{pFZ4}).
By naturality of the AHSS, $
2 \{h_1^2\beta_9x U\} $
has AH filtration $9$.
Thus, the desired nontrivial extension follows.
\end{proof}

Since $E^2_{12,1}=0$, the filtration gives
\[
0\longrightarrow\mm F_{11,2}/\mm F_{8,5}
\longrightarrow ko_{13}(\mm M\varphi)/\mm F_{8,5}
\longrightarrow E^\infty_{13,0}\longrightarrow0.
\]
Both outer groups are $\Z$, so this sequence splits. Thus
\begin{equation}\label{eq:detect-varphi-free}
ko_{13}(\mm M\varphi)/\mm F_{8,5}
=\Z\langle\rho_{\q},\bar\rho_{\q}\rangle.
\end{equation}
In particular, every torsion element belongs to $\mm F_{8,5}$.
The analogous notation and statements will be used for $\varphi'$.

We now consider $\omega=j^*\varphi$ over the common
$P$-bundle $X$. Recall \eqref{eq:Y6-X-homology-low}
\[
H_9(\mm M\omega;\Z)=\Z\langle z^\vee U_*,(xy)^\vee U_*\rangle,
\qquad
H_{11}(\mm M\omega;\Z)
=\Z\langle(xz)^\vee U_*,w^\vee U_*\rangle.
\]
$$H_{13}(\M\omega ;\Z)\cong\Z\langle (x^2z)^\vee U_\ast \rangle, \qquad H_{14}(\M\omega ;\mathbb Q)=0.$$

\begin{lemma}\label{koXq}
In the AHSS for $ko_*(\mm M\omega)$,
\[
E^\infty_{13,0}(ko,\M\omega)=E^2_{13,0}(ko,\M\omega)\cong \Z,
\quad E^\infty_{9,4}(ko,\M\omega)=E^2_{9,4}(ko,\M\omega)\cong\Z^2.
\]
Let $\theta\in ko_{13}(\M\omega)$ be the element detected by $(x^2z)^\vee U_\ast$.

Moreover,   $(xz)^\vee U_*\in E^2_{11,2}$ survives to $E^\infty$ and detects $\bar \theta\in ko_{13}(\M\omega)$. 
Let $\psi_z,\psi_{xy}\in ko_{13}(\M\omega)$ be the elements detected by $z^\vee U_\ast,  (xy)^\vee U_*\in E^2_{9,4}$. Then, $\bar \theta$ is torsion free, and there exist $a,b\in \Z$ such that 
\begin{equation}\label{eq:detect-general-extension}
	sa+tb\equiv 1\pmod 2\quad \text{and}\quad 2\bar \theta=a\psi_z+b \psi_{xy}\pmod{\mb{F}^{8,5}}
\end{equation}
 where $t$ is defined in \eqref{cohomology_z_xy_beta}, $\mb{F}^{p,q}$ denotes the AHSS filtration for $ko_*(\mm M\omega)$.
\end{lemma}

\begin{proof}
The generator of $E^2_{13,0}$ is dual to $zx^2U$.
By $w_2(\omega)\ne 0$ and \eqref{eq:restricted-P-Sq2},
\[
\Sq^2(xzU)=0,
\qquad
\Sq^2(wU)\in\Z_2\langle xwU\rangle.
\]
 Thus the outgoing $d_2$ on this generator is zero.
The other targets before $d_9$ vanish by the coefficient groups of
$ko$ and the homology groups of $\M\omega$.

The remaining possible differentials are $d_9,d_{11},d_{13}$.
At the relevant pages, the targets of $d_9$ and $d_{13}$ are
subgroups of the corresponding free integral homology groups. By the rational argument, the $d_9$ and $d_{13}$ are zero. 
 The target of $d_{11}$ is instead $\Z_2$;
its survival as in \eqref{Einfty_210_Mxi} gives $d_{11}=0$.
This implies $E^\infty_{13,0}=E^2_{13,0}$.

For $E^2_{9,4}\cong\Z^2$, the outgoing $d_5$ and $d_9$ again have
free targets and vanish rationally, while $d_7=0$ by
\eqref{Einfty_210_Mxi}. Incoming differentials have finite sources. For the possible source in bidegree $(14,0)$,
$H_{14}(\M\omega ;\mathbb Q)=0$ gives the finiteness.
Thus $E^\infty_{9,4}=E^2_{9,4}$.

For $E^r_{11,2}$ there are no outgoing differentials: the possible
$d_7$ and $d_{11}$ have free targets, while $d_9=0$ by \eqref{Einfty_210_Mxi}. 
By \eqref{cohomology_z_xy_beta}, we have
\begin{equation}\label{j_beta11}
	j^*(\beta_{11})=xj^*(\beta_9)
=sxz+t x^2y=sxz+6t w.
\end{equation}
Modulo $2$, this is $xz$, since $s$ is odd. Hence $j_*$ maps
$(xz)^\vee U_*$ to the surviving nonzero class of
$E^\infty_{11,2}(ko,\mm M\varphi)$, and maps $w^\vee U_*$ to zero.
Naturality prevents any incoming differential from killing
$(xz)^\vee U_*$.

Finally, the map $j_*$ gives a diagram of extensions
\[
\begin{tikzcd}[column sep=small,row sep=large]
0\arrow[r]& E^\infty_{9,4}(ko,\M\omega)= \Z^2\arrow[r]\arrow[d,"{(s,t)}"']
&\mathcal G\arrow[r]\arrow[d]
&\Z_2\subset E^\infty_{11,2}(ko,\M\omega)\arrow[r]\arrow[d,"\cong"]&0\\
0\arrow[r]&E^\infty_{9,4}(ko,\M\varphi)=\Z\arrow[r]
&\mm F_{11,2}/\mm F_{8,5}\arrow[r]
&E^\infty_{11,2}(ko,\M\varphi)=\Z_2\arrow[r]&0
\end{tikzcd}
\]
where we omit $w^\vee U_\ast\in E^\infty_{11,2}(ko,\M\omega)$. Since
the bottom sequence is non-split  (cf. Lemma~\ref{koM4}), the top sequence is
non-split.
 This proves 
$$2\bar \theta=a\psi_z+b \psi_{xy}\pmod{\mb{F}_{8,5}}$$
where $\bar \theta$ is detected by $(xz)^\vee U_\ast$.

Now we show $sa+tb\equiv 1\pmod 2$.
Recall that $\bar\rho_\q$ is provided by $E^2_{11,2}(ko,\M \varphi)$ and $2\bar\rho_\q$ is provided by $E^2_{9,4}(ko,\M \varphi)$ (cf. \eqref{eq:detect-varphi-free}). By \eqref{j_beta11},
\[
j_*:
E^\infty_{11,2}(ko,\mm M\omega)
\longrightarrow
E^\infty_{11,2}(ko,\mm M\varphi)
\]
maps $(xz)^\vee U_\ast$ to $\beta_{11}^\vee U_\ast$.
Thus
$$j_\ast(\bar \theta)=m \bar\rho_\q \pmod{\mm{ F}_{8,5}},\quad m\equiv 1\pmod 2.$$
 \eqref{homology_z_xy_beta} and the naturality of Hurewicz homomorphism gives 
\[
j_*(\psi_z)=2s\bar \rho_\q,\qquad
j_*(\psi_{xy})=2t\bar \rho_\q
\pmod{\mm{ F}_{8,5}}.
\]
Therefore, we obtain
\[
2(sa+tb)\bar \rho_\q
=
2j_*(\bar\theta)
=
2m\bar \rho_\q \pmod{\mm{ F}_{8,5}}.
\]
Since $2\bar \rho_\q$ generates an infinite cyclic group,
$
sa+tb=m\equiv1\pmod2.$
\end{proof}

\subsection{Integral homology of $ko$-Thom spectra and Hurewicz images}
\label{subsec:Hurewicz-ko13}\quad

We first record the coefficient information needed for the integral
homology calculation. Write $\iota:S^0\to H\Z$ for the unit.
For a spectrum $X$, the ordinary Hurewicz homomorphism is
\[
\iota_*:ko_{13}(X)=\pi_{13}(ko\wedge X)
\longrightarrow H_{13}(ko\wedge X).
\]
It is distinct from the edge homomorphism
$\chi_*:ko_{13}(X)\to H_{13}(X)$.

\begin{proposition}\label{hko}
The integral homology of $ko$ satisfies:
\begin{enumerate}
\item $H_n(ko)$ has a free abelian summand only if $n\equiv0\pmod4$.
\item $H_n(ko)=0$ for $1\le n\le3$ and for $n=5$.
\item $H_6(ko)$ and $H_7(ko)$ are finite $2$-primary groups.
\item $H_4(ko)\cong\Z$, and, after compatible choices of generators,
\[
\iota_*:\pi_4(ko)\longrightarrow H_4(ko)
\]
is multiplication by $24$.
\end{enumerate}
\end{proposition}

\begin{proof}
The rational splitting
\[
ko_{\mathbb Q}\simeq\bigvee_{i\ge0}\Sigma^{4i}H\mathbb Q
\]
shows that $H_n(ko;\mathbb Q)$ is zero unless $n\equiv0\pmod4$.
Since the integral homology groups are finitely generated in each
degree, this proves (1).

The Atiyah--Bott--Shapiro orientation
$\alpha:\mm{MSpin}\to ko$ is an $8$-equivalence
\cite{HebeJoach}. The integral Thom isomorphism therefore gives
\[
H_n(ko)\cong H_n(\mm{MSpin})
\cong H_n(\mm{BSpin}),\qquad n\le7.
\]
The low-dimensional calculation in \cite{Giambalvo1969}
gives (2), (3), and $H_4(ko)\cong\Z$.

For the coefficient Hurewicz homomorphism, consider
\[
E^2_{p,q}=H_p(ko;\pi_q(S^0))\Longrightarrow\pi_{p+q}(ko).
\]
By (2), the only possible differential into $(0,3)$ is
\[
d_4:E^4_{4,0}=\Z\longrightarrow E^4_{0,3}=\Z_{24}.
\]
Since $\pi_3(ko)=0$, this map is surjective.
The outgoing $d_2$ and $d_3$ from $(4,0)$ have zero targets. Hence
\[
E^\infty_{4,0}=\ker d_4=24\Z\subset H_4(ko).
\]
This is the image of the Hurewicz homomorphism. As $\pi_4(ko)\cong\Z$,
the homomorphism is multiplication by $24$, up to the choice of sign.
\end{proof}

We next use the AHSS
\[
E^2_{p,q}(H\wedge ko,X)=H_p(X;H_q(ko))
\Longrightarrow H_{p+q}(ko\wedge X)
\]
to locate the free homology classes that occur in the Hurewicz images.

\begin{lemma}\label{H1394}
For $X=\mm M\varphi$ or $\mm M\omega$ and $q=0,4$,
\[
E^\infty_{13-q,q}(H\wedge ko,X)
=E^2_{13-q,q}(H\wedge ko,X).
\]
\end{lemma}

\begin{proof}
We first give an argument for $\M\varphi$.

$
E^2_{13,0}(H\wedge ko,\mm M\varphi)
=
H_{13}(\mm M\varphi)
\cong\Z$
is generated by
$
[\mathcal B_{\q}]U,$
which is the Hurewicz image of the Thom--Pontryagin class
(cf. Lemma~\ref{E130koMvarphi}). Hence it survives to
\(E^\infty\) by the naturality of AHSS.

Now consider
$E^2_{9,4}(H\wedge ko,\mm M_\varphi).$
The possible outgoing differentials are
\[
d_5,d_7,d_9.
\]
We compare with the Thom spectrum
$
\mm M\xi ,
$
where \(\xi\) is defined in Section \ref{Normalbg}.

The Thom map
$
c:\mm M\varphi\longrightarrow \mm M\xi
$
induces a morphism of AHSSs. Since
\[
H_p(\mm M\varphi;\Z)
\cong
H_p(\mm M\xi;\Z),
\qquad p\le4,
\]
the relevant target terms agree.
For \(\mm M\xi\), the K{\"u}nneth product
\[
H_*(ko)\otimes H_*(\mm M\xi)
\longrightarrow
H_*(ko\wedge\mm M\xi)
\]
is an isomorphism, and thus the AHSS collapses.
Therefore the corresponding differentials $d_5,d_7,d_9$ for
\(\mm M\varphi\) vanish by naturality.

The proof for \(\mm M\omega\) is identical, using $E^\infty_{13,0}(ko,\M\omega)=E^2_{13,0}(ko,\M\omega)$ (cf. Lemma \ref{koXq}) and the Thom map
$
 j_\ast:\mm M\omega\longrightarrow\mm M\varphi .
$
\end{proof}

Fix a generator $\mathtt{k}_4\in H_4(ko)$.  Lemma~\ref{H1394} now gives
the free generators below. The symbol $\wedge $ denotes the wedge product
$$\wedge: H_p(X;\Z)\otimes H_q(Y;\Z)\to H_{p+q}(X\wedge Y;\Z).$$
 The symbol $\otimes$ denotes a lift from
an associated-graded $E^\infty$-group.

\begin{corollary}\label{homologykovarphiomega}
There are integral bases
\[
H_{13}(ko\wedge\mm M\varphi)/\operatorname{Tor}
=\Z\langle
1\wedge[\mathcal B_{\q}]U_*,\
\mathtt{k}_4\otimes\beta_9^\vee U_*
\rangle
\]
and
\[
H_{13}(ko\wedge\mm M\omega)/\operatorname{Tor}=\Z\langle
1\wedge(zx^2)^\vee U_*,\
\mathtt{k}_4\otimes z^\vee U_*,\
\mathtt{k}_4\otimes(xy)^\vee U_*
\rangle.
\]
\end{corollary}

\begin{proof}
The only nonzero rational K\"unneth contributions in total degree
$13$ occur in bidegrees $(0,13)$ and $(4,9)$. Thus these free cross products form bases after
quotienting by torsion.
\end{proof}

We retain the classes $\rho_{\q},\bar\rho_{\q}\in ko_{13}(\M\varphi)$.
Recall the classes $\theta,\bar \theta, \psi_z,\psi_{xy} \in ko_{13}(\mm M\omega)$ (cf. Lemma~\ref{koXq}). We now compute their Hurewicz
images.

\begin{lemma}\label{Huretheta}
For some integers $t_1,t_2$ with $6\mid t_2$,
\[
\begin{aligned}
\iota_*(\theta)
={}&1\wedge(zx^2)^\vee U_*
+t_1\mathtt{k}_4\otimes z^\vee U_*\\
&+t_2\mathtt{k}_4\otimes(xy)^\vee U_*\in H_{13}(ko\wedge\mm M\omega;\Z)/\operatorname{Tor}
\end{aligned}
\]
\end{lemma}

\begin{proof}
The Hurewicz map induces the homomorphism
$$E^\infty_{13,0}(ko,\M\omega)\to E^\infty_{13,0}(H\Z\wedge ko,\M\omega)$$
where the coefficient
map is the identity on $\Z$. The edge image of
$\theta$ is $(zx^2)^\vee U_*$. Hence its coefficient on the
filtration-$13$ generator is $1$. Corollary~\ref{homologykovarphiomega}
gives the stated expression with $t_1,t_2\in\Z$.

It remains to determine the divisibility of the coefficient
$t_2$.  
Consider the following commutative diagram:
\[
\begin{tikzcd}[column sep=large]
\pi_{13}(ko\wedge\mm M\omega)
\arrow[r,"\chi\wedge1"]
\arrow[d,"\iota_\ast=(\iota\wedge 1\wedge 1)_\ast"']
&
\pi_{13}(H\Z\wedge\mm M\omega)\arrow[d,"\iota_{2\ast}=(\iota\wedge 1\wedge 1)"]
\\
\pi_{13}(H\wedge ko\wedge\mm M\omega)
\arrow[r,"\iota_{1\ast}=(1\wedge \iota\wedge 1)_\ast"]& \pi_{13}(H\Z\wedge H\Z\wedge\mm M\omega)
\end{tikzcd}
\]

The right-hand homomorphism  has a mod prime $p$ reduction:
$$\iota_{2\ast} :{{H}\mathbb{Z}_p}_\ast(\mm{M}\omega)\to  {H\mathbb{Z}_p}_\ast({H}\mathbb{Z}_p\wedge \mm{M}\omega)\cong \mathscr{A}_{p\ast}\otimes {H\mathbb{Z}_p}_\ast(\mm{M}\omega)$$
\begin{equation}
	(zx^2)_p^\vee U\to 1\otimes (zx^2)_p^\vee U+\cdots \in \mathscr{A}_{p\ast}\otimes H_{\ast}(\mm{M}\omega;\mb{Z}_p) \label{comoduleoperation}
\end{equation}
 which indeed is the $\mathscr{A}_{p\ast}$-comodule operation \cite{Sw}, where $\mathscr{A}_{p\ast}$ is the dual of the mod $p$ Steerod algebra.

 The mod $p$ reduction of $\iota_{1\ast}$ is
$$H_\ast(ko;\Z_p)\otimes H_\ast(\M\omega;\Z_p)\xrightarrow{\chi_p\otimes 1} H_\ast(H\Z_p;\Z_p)\otimes H_\ast(\M\omega;\Z_p)$$
For $p=2$, $$\chi_2(\mathtt{k}_4)={\Sq^4}^\vee\in H_4(H\Z_2;\Z_2).$$ 
Recall $w_2(\omega)\ne 0$ and $w_4(\omega)=0$. By Proposition \eqref{prop:restricted-P-squares},
$$\Sq^4(\bar x\bar y) U)=\Sq^2(\bar x\bar y)w_2(\omega)U+(\bar x\bar y)w_4(\omega)U=0.$$ 
Then, in \eqref{comoduleoperation}, the image of $(zx^2)^\vee U$ does not contain the term
$${\Sq^4}^\vee\otimes (\bar x\bar y)^\vee U.$$
Therefore, the above commutative diagram gives 
 $2|t_2.$

For $p=3$, $$\chi_3(\mathtt{k}_4)={\mathcal P^1}^\vee\in H_4(H\Z_3;\Z_3).$$ By Proposition \ref{prop:restricted-P-P1},
$$\mathcal P^1(( xy)_3 U)=\mathcal P^1(xy)_3 U+( x y)_3\,p_1(\omega)U=k(x^3y)_3=0$$
where $p_1(w)=kx^2\in H^4( X_\q)$. Similarly, we have $3|t_2$.
\end{proof}

\begin{lemma}\label{Hure-rho130-phi}
For the class $\rho_{\q}=(\alpha_8)_*(\varsigma_{\q})$,
\[
\iota_*(\rho_{\q})=1\wedge[\mathcal B_{\q}]U_*
\quad\text{in }H_{13}(ko\wedge\mm M\varphi)/\operatorname{Tor}.
\]
\end{lemma}

\begin{proof}
Before applying $\alpha_8$, the Hurewicz image of $\varsigma_{\q}$ is obtained from
$$(g,1)_*[\mathcal B_{\q}]U_\ast\in H_{13}(\mm{MO}\langle8\rangle\wedge\M\varphi)$$ by the Thom isomorphism.
There is only one nonzero rational K\"unneth summand in total degree
$13$ for $\mm{MO}\langle8\rangle\wedge\M\varphi$. And the free summand is generated by $1\wedge [\mathcal B_{\q}]U_\ast$. 
Thus the Hurewicz image of $\varsigma_{\q}$ is
$1\wedge[\mathcal B_{\q}]U_*$ modulo torsion.
Applying $\alpha_8$ and using naturality of $\wedge$ gives the result.
\end{proof}

\begin{lemma}\label{Hure-rho94}
With the notations above,
\[
\iota_*(\bar\rho_{\q})
=12\mathtt{k}_4\otimes\beta_9^\vee U_*
\quad\text{in }H_{13}(ko\wedge\mm M\varphi)/\operatorname{Tor},
\]
\[
\iota_*(\bar\theta)
=12a\mathtt{k}_4\otimes z^\vee U_*
+12b\mathtt{k}_4\otimes(xy)^\vee U_*
\quad\text{in }H_{13}(ko\wedge\mm M\omega)/\operatorname{Tor}
\]
where $a,b\in \Z$ are defined in Lemma \ref{koXq}.
\end{lemma}

\begin{proof}
The coefficient Hurewicz map in degree $4$ is multiplication by
$24$, by Proposition~\ref{hko}. Hence
\[
\begin{aligned}
\iota_*(2\bar\rho_{\q})&=24\mathtt{k}_4\otimes\beta_9^\vee U_*,\\
\iota_*(\psi_z)&=24\mathtt{k}_4\otimes z^\vee U_*,\\
\iota_*(\psi_{xy})&=24\mathtt{k}_4\otimes(xy)^\vee U_*
\end{aligned}
\]
modulo torsion. Then we have the first formula.
Applying Hurewicz homomorphism to $2\bar \theta=a\psi_z+b\psi_{xy}$ gives the second formula.
\end{proof}

\subsection{Final detection}

\begin{lemma}\label{common_image_rhoq}
Suppose 
$\sigma_2(\q)=\sigma_2(\q'),$
$\sigma_3(\q)=\sigma_3(\q')>0.$
Then, for every prime $p\nmid s$,
\[
f_*(\rho_{\q})=f'_*(\rho_{\q'})
\quad\text{in}\quad
\bigl(ko_{13}(\mm M\eta)/\bar F_{8,5}\bigr)_{(p)},
\]
where $\bar F_{p,q}$ denotes the AHSS filtration
for $ko_{p+q}(\mm M\eta)$.
\end{lemma}

\begin{proof}
By the restricted-bundle comparison, we have a common
space $X$, a common bundle $\omega$, and the
homotopy-commutative diagram
\begin{equation}\label{commuta_Thomspec}
\begin{tikzcd}[column sep=large,row sep=large]
\mm M\omega
\arrow[r,"j"]
\arrow[d,"j'"']
&
\mm M\varphi
\arrow[d,"f"]
\\
\mm M\varphi'
\arrow[r,"f'"']
&
\mm M\eta .
\end{tikzcd}
\end{equation}

For convenience, put
\[
Q=ko_{13}(\mm M\eta)/\bar F_{8,5},
\qquad
A=\bar F_{11,2}/\bar F_{8,5}\subset Q.
\]
By Theorem~\ref{pFZ4} and Lemma~\ref{koketa94},
there is an exact sequence
\[
0\longrightarrow
\Z_6\oplus\Z_s
\longrightarrow A
\longrightarrow(\Z_2)^2
\longrightarrow0.
\]
Since $\gcd(s,6)=1$, it follows that
\begin{equation}\label{eq:detect-target-exponent}
12A_{(p)}=0
\qquad\text{for every prime }p\nmid s.
\end{equation}

We prove the lemma by comparing the images of the
free generators of $ko_{13}(\mm M\varphi)$ and
$ko_{13}(\mm M\varphi')$.

\medskip
\noindent
\textbf{Step 1. Comparison of the filtration-eleven classes.}

We first use the Hurewicz homomorphism to determine
the maps induced by $j$ and $j'$ on these free quotients.

Recall
$
ko_{13}(\mm M\varphi)/\mm F_{8,5}
\cong
\Z\langle\rho_{\q}\rangle
\oplus
\Z\langle\bar\rho_{\q}\rangle$.
Recall
$\psi_z,\psi_{xy}\in\mb F_{9,4}\subset ko_{13}(\M\omega)$
lifting $z^\vee U_*$ and $(xy)^\vee U_*$,
respectively. Applying Hurewicz homomorphism, we have
\[
\begin{aligned}
\iota_*(\psi_z)
&=
24\mathtt{k}_4\otimes z^\vee U_*,
\\
\iota_*(\psi_{xy})
&=
24\mathtt{k}_4\otimes(xy)^\vee U_*
\end{aligned}
\]
modulo torsion.

On the other hand, by Lemma~\ref{Hure-rho94},
\[
\iota_*(\bar\rho_{\q})
=
12\mathtt{k}_4\otimes\beta_9^\vee U_*.
\]

Now consider the commutative diagram
\[
\begin{tikzcd}[column sep=large,row sep=large]
ko_{13}(\mm M\omega)
\arrow[r,"j_*"]
\arrow[d,"\iota_*"']
&
ko_{13}(\mm M\varphi)
\arrow[d,"\iota_*"]
\\
H_{13}(ko\wedge\mm M\omega;\Z)
\arrow[r,"(1\wedge j)_*"']
&
H_{13}(ko\wedge\mm M\varphi;\Z).
\end{tikzcd}
\]
Then we have
\begin{equation}\label{eq:detect-j-basis}
\begin{aligned}
j_*(\psi_z)&=2s\bar\rho_{\q},
&
j_*(\psi_{xy})&=2t\bar\rho_{\q},
\\
j'_*(\psi_z)&=2s\bar\rho_{\q'},
&
j'_*(\psi_{xy})&=2t'\bar\rho_{\q'}
\end{aligned}
\end{equation}
in the respective quotients modulo filtration $8$.

Applying \eqref{commuta_Thomspec} to $\psi_z$ gives
\[
2s\bigl(
f_*(\bar\rho_{\q})
-
f'_*(\bar\rho_{\q'})
\bigr)=0
\quad\text{in }A.
\]
For every odd prime $p\nmid s$,
multiplication by $2s$ is invertible after
$p$-localization. Therefore,
$
f_*(\bar\rho_{\q})
=
f'_*(\bar\rho_{\q'})$ in $A_{(p)}.$

It remains to establish the same equality at $p=2$.

Recall $\bar\theta\in\mb F_{11,2}$ lifting
$(xz)^\vee U_*$. By Lemma \ref{Hure-rho94},
\[
\iota_*(\bar\theta)
=
12a\mathtt{k}_4\otimes z^\vee U_*
+
12b\mathtt{k}_4\otimes(xy)^\vee U_*
\]
modulo torsion.
Applying $j_*$ and $j'_*$, respectively, we obtain
\begin{equation}\label{j_bar_theta}
\begin{aligned}
j_*(\bar\theta)
&=(sa+tb)\bar\rho_{\q},
\\
j'_*(\bar\theta)
&=(sa+t'b)\bar\rho_{\q'}
\end{aligned}
\end{equation}
modulo filtration $8$.

By the commutativity of \eqref{commuta_Thomspec},
equation \eqref{j_bar_theta} yields
\begin{equation}\label{compute_bartheta}
(sa+tb)f_*(\bar\rho_{\q})
=
(sa+t'b)f'_*(\bar\rho_{\q'})
\quad\text{in }A.
\end{equation}
By Proposition~\ref{ttprime_differential},
$
24\mid(t-t').$
Hence
\[
(sa+t'b)-(sa+tb)=(t'-t)b
\]
is divisible by $24$.
In view of \eqref{eq:detect-target-exponent},
the right-hand side of \eqref{compute_bartheta}
can therefore be rewritten as
\[
(sa+tb)f'_*(\bar\rho_{\q'}).
\]
Take the $2$-localization, thus
\[
(sa+tb)
\bigl(
f_*(\bar\rho_{\q})
-
f'_*(\bar\rho_{\q'})
\bigr)=0
\quad\text{in }A_{(2)}.
\]
Since $sa+tb$ is odd (cf. \eqref{eq:detect-general-extension}), multiplication by this
integer is invertible on the finite $2$-primary group.
It follows that
$
f_*(\bar\rho_{\q})
=
f'_*(\bar\rho_{\q'})$
in $A_{(2)}.$

Combining the odd-primary and $2$-primary cases,
we conclude that
\begin{equation}\label{compar_bar_rho}
f_*(\bar\rho_{\q})
=
f'_*(\bar\rho_{\q'})
\quad\text{in }A_{(p)}
\qquad(p\nmid s).
\end{equation}

\medskip
\noindent
\textbf{Step 2. Comparison of the filtration-thirteen classes.}

Recall $\theta\in ko_{13}(\mm M\omega)$
detected by
\[
(zx^2)^\vee U_*
\in E^\infty_{13,0}(ko,\mm M\omega).
\]
We now determine $j_*(\theta)$ and
$j'_*(\theta)$ by their Hurewicz images.

By Lemma~\ref{Huretheta},
\begin{equation}\label{eq:detect-hure-theta}
\begin{aligned}
\iota_*(\theta)
={}&
1\wedge(zx^2)^\vee U_*
\\
&+
t_1\mathtt{k}_4\otimes z^\vee U_*
+
t_2\mathtt{k}_4\otimes(xy)^\vee U_*
\end{aligned}
\end{equation}
modulo torsion, where
$
t_1,t_2\in\Z,$ $6\mid t_2.$

Recall that
$
j_*((zx^2)^\vee U_*)
=
s[\mathcal B_{\q}]U_*$.
Together with \eqref{homology_z_xy_beta},
naturality of the Hurewicz homomorphism gives
\[
\begin{aligned}
\iota_*j_*(\theta)
={}&
s\bigl(1\wedge[\mathcal B_{\q}]U_*\bigr)
\\
&+
(st_1+tt_2)
\mathtt{k}_4\otimes\beta_9^\vee U_*
\end{aligned}
\]
modulo torsion.

On the other hand, by Lemmas
\ref{Hure-rho130-phi} and \ref{Hure-rho94},
\[
\iota_*(\rho_{\q})
=
1\wedge[\mathcal B_{\q}]U_*,
\quad
\iota_*(\bar\rho_{\q})
=
12\mathtt{k}_4\otimes\beta_9^\vee U_*.
\]
The analogous formulas hold for $\q'$.

Since
$
ko_{13}(\mm M\varphi)/\mm F_{8,5}
\cong
\Z\langle\rho_{\q}\rangle
\oplus
\Z\langle\bar\rho_{\q}\rangle,$
the Hurewicz calculation determines the
coefficients of $j_*(\theta)$ uniquely.
In particular, there exist integers $n,n'$ such that
\begin{equation}\label{eq:detect-j-theta}
\begin{aligned}
j_*(\theta)
&=
s\rho_{\q}+n\bar\rho_{\q},
\\
j'_*(\theta)
&=
s\rho_{\q'}+n'\bar\rho_{\q'}
\end{aligned}
\end{equation}
in the respective quotients modulo filtration $8$, where
\begin{equation}\label{eq:detect-n-coefficients}
12n=st_1+tt_2,
\qquad
12n'=st_1+t't_2.
\end{equation}

Subtracting these identities gives
\[
12(n-n')=(t-t')t_2.
\]
Since
$
24\mid(t-t')$ and $6\mid t_2,$
we obtain
\begin{equation}\label{diff_n_nprime}
12\mid(n-n').
\end{equation}

Now fix a prime $p\nmid s$.
By \eqref{compar_bar_rho},
\[
f_*(\bar\rho_{\q})
=
f'_*(\bar\rho_{\q'})
\quad\text{in }A_{(p)}.
\]
Consequently,
\[
f_*(n\bar\rho_{\q})
-
f'_*(n'\bar\rho_{\q'})
=
(n-n')f_*(\bar\rho_{\q}).
\]
By \eqref{diff_n_nprime} and
\eqref{eq:detect-target-exponent},
the right-hand side vanishes.
Therefore,
\begin{equation}\label{compar_final_barrho}
f_*(n\bar\rho_{\q})
=
f'_*(n'\bar\rho_{\q'})
\quad\text{in }A_{(p)}\subset  Q_{(p)}.
\end{equation}

Finally, the homotopy-commutative square
\eqref{commuta_Thomspec} implies
\begin{equation}\label{compute_theta}
f_*j_*(\theta)
=
f'_*j'_*(\theta)
\quad\text{in }Q.
\end{equation}
 Using \eqref{eq:detect-j-theta}, \eqref{compar_final_barrho} and \eqref{compute_theta}, we obtain
\[
s\bigl(
f_*(\rho_{\q})
-
f'_*(\rho_{\q'})
\bigr)=0
\quad\text{in }Q_{(p)}.
\]
Since $p\nmid s$, multiplication by $s$
is invertible after $p$-localization.
Hence
\[
f_*(\rho_{\q})
=
f'_*(\rho_{\q'})
\quad\text{in }Q_{(p)}.
\]
This completes the proof.
\end{proof}

We conclude by combining this comparison with the linking invariant.

\begin{proposition}\label{Bordism_compare_lowF_85}
Suppose $\gcd(s,5)=1$,
\[
\sigma_i(\q)=\sigma_i(\q')\quad(i=2,3),
\qquad\sigma_5(\q)\equiv\sigma_5(\q')\pmod s.
\]
Then the corresponding normal $5$-smoothings satisfy
\[
[\mathcal B_{\q},(g,f)]
=[\mathcal B_{\q'},(g',f')]
\quad\text{in }
\Omega_{13}^{\mm O\langle8\rangle}(\eta)/F_{8,5}.
\]
\end{proposition}

\begin{proof}
By Lemma~\ref{130Zs}, the $\Z_s$-summand of this quotient is
specified by the oriented linking invariant. The assumed congruence
for $\sigma_5$ therefore identifies these components.

For $p=2,3,5$, we have $p\nmid s$. Lemmas
\ref{common_image_rhoq} and \ref{common_image_varsigma} give the desired equality
after localization at each of these primes. The previously computed
group decomposition is
\[
\Omega_{13}^{\mm O\langle8\rangle}(\eta)/F_{8,5}
\cong\Z_s\oplus
\bigoplus_{p=2,3,5}
\bigl(\Omega_{13}^{\mm O\langle8\rangle}(\eta)/F_{8,5}\bigr)_{(p)}.
\]
Equality of all the displayed primary components proves the assertion.
\end{proof}

\section{Bordism in the PL category}\label{PL}

The natural map $\mm{BO}\to\mm{BPL}$ lifts to a map
$
\mm{BO}\langle8\rangle
\longrightarrow
\mm{BPL}\langle8\rangle$
between the $7$-connected covers and induces a Thom spectrum map
\[
T:\mm{MO}\langle8\rangle
\longrightarrow
\mm{MPL}\langle8\rangle.
\]
Regarding $\eta$ as a stable PL bundle, we obtain the forgetful map
\[
T_*:
\Omega_{13}^{\mm O\langle8\rangle}(\eta)
\longrightarrow
\Omega_{13}^{\mm{PL}\langle8\rangle}(\eta).
\]
Throughout this section, $s>0$ and $\gcd(s,6)=1$.

Write
\[
\begin{aligned}
E^2_{p,q}
&=
H_p(\mm M\eta;\Omega_q^{\mm O\langle8\rangle}),\\
\widetilde{ E}^2_{p,q}
&=
H_p(\mm M\eta;\Omega_q^{\mm{PL}\langle8\rangle})
\end{aligned}
\]
for the $E^2$-pages of the two AHSSs.
Their filtrations are denoted by $F_{p,q}$ and
${\widetilde F}^{p,q}$, respectively. 

We shall prove that $T_*$ vanishes on the  filtration subgroup
$ F_{4,9}$ and acts by multiplication by $7$ on the
$\Z_s$-contribution in filtration $5$. We then explain how this
changes the comparison condition involving $\sigma_4$.

\begin{theorem}\label{top13}
$T_*$ vanishes on the  filtration subgroup
$ F_{4,9}$.
\end{theorem}

\begin{proof}
By $\Omega_{11}^{\mm{PL}\langle8\rangle}
=
\Omega_{13}^{\mm{PL}\langle8\rangle}
=0$ and the low-dimensional homology of $\M\eta$, we have
$\widetilde{E}^2_{p,13-p}=0$ for $p\le 3$. Thus,
\begin{equation}\label{PL49}
	\widetilde{F}^{4,9}\cong \widetilde{E}^\infty_{4,9}.
\end{equation}
Then it is enough to consider 
$T_\ast:{E}^2_{4,9}\to \widetilde{E}^2_{4,9}.$
This homomorphism is zero since the coefficient homomorphism
$T_\ast:\Omega_{9}^{\mm{O}\langle8\rangle}\to \Omega_{9}^{\mm{PL}\langle8\rangle}$
is zero. Hence, the desired result follows.
\end{proof}

\subsection{$9$-dimensional PL bordism}

\begin{lemma}\label{PL9primary}
 For odd prime $p$, $\Omega_9^{\mm{PL}\langle8\rangle}$ has no $p$-prime component.
\end{lemma}

\begin{proof}
It is well known that
there is a short exact sequence
\[
0\longrightarrow
\pi_9(\mm{BO})
\longrightarrow
\pi_9(\mm{BPL})
\longrightarrow
\pi_8(\mm{PL/O})
\longrightarrow0.
\]
By Bott periodicity and the calculation
of Kervaire--Milnor \cite{KervaMilnor}, we have
\[
\pi_9(\mm{BO})\cong\Z_2,
\qquad
\pi_8(\mm{PL/O})
\cong\Theta_8\cong\Z_2.
\]
Let $p$ be an odd prime. Then, $\pi_9(\mm{BPL})_{(p)}=0$, and thus
$$\pi_9(\mm{BPL}\langle 8\rangle)_{(p)}=0.$$
Since $\mm{BPL}\langle8\rangle$ is $7$-connected, 
Hurewicz theorem gives a surjection 
\[
\pi_9(\mm{BPL}\langle8\rangle)
\to H_9(\mm{BPL}\langle8\rangle ).
\]
Consequently, 
${H}_9(\mm{BPL}\langle8\rangle)_{(p)}=0$. 
By the Thom isomorphism,
$${H}_9(\mm{MPL}\langle8\rangle)_{(p)}=0.$$

Now consider the AHSS
\[
E^2_{p,q}(S^0)
=
H_p(\mm{MPL}\langle8\rangle;
\pi_q(S^0))
\Longrightarrow
\pi_{p+q}(\mm{MPL}\langle8\rangle).
\]
The only possibly
nonzero terms in total degree $9$ are
\[
E^2_{0,9}(S^0),\qquad
E^2_{8,1}(S^0),\qquad
E^2_{9,0}(S^0).
\]
All three have no $p$-prime component.
 The same holds for
the AHSS extensions, proving
the assertion.
\end{proof}

\subsection{The coefficient map in degree eight}\quad

We next compute the coefficient map responsible for the factor $7$.

\begin{lemma}\label{BO8toBPL8times7}
There is a short exact sequence
\[
0
\longrightarrow H_8(\mm{BO}\langle8\rangle;\Z)
\xrightarrow{T_\ast} H_8(\mm{BPL}\langle8\rangle;\Z)
\longrightarrow\Z_{28}
\longrightarrow0,
\]
where the first two groups are $\Z$ and $\Z\oplus\Z_4$.
After compatible choices of generators, the induced map $T_\ast$ is multiplication by $7$.
\end{lemma}

\begin{proof}
The standard smoothing calculation gives
\[
0
\longrightarrow\pi_8(\mm{BO})\cong\Z
\longrightarrow\pi_8(\mm{BPL})\cong\Z\oplus\Z_4
\longrightarrow\Z_{28}
\longrightarrow0;
\]
see \cite[Lemma~5.1]{Kramer}.
Passing to the $7$-connected covers does not change $\pi_8$, and
Hurewicz identifies $\pi_8$ with $H_8$ for both covers.
The torsion subgroup $\Z_4$ injects into the $\Z_{28}$.
Thus, by taking suitable generators of  $H_8(\mm{BPL}\langle 8\rangle)$ and using the exactness, we identify $T_\ast$ with the multiplication by $7$.
\end{proof}

\begin{lemma}\label{PL8Zs}
The coefficient groups in degrees at most seven are
\[
\begin{array}{c|cccccccc}
n&0&1&2&3&4&5&6&7\\ \hline
\Omega_n^{\mm{PL}\langle8\rangle}
&\Z&\Z_2&\Z_2&\Z_{24}&0&0&\Z_2&0.
\end{array}
\]
For every $s>0$ with $\gcd(s,6)=1$, compatible identifications give
\begin{equation}\label{eq:PL-degree8-times7}
T_*\otimes1:
\Omega_8^{\mm O\langle8\rangle}\otimes\Z_s
\longrightarrow
\Omega_8^{\mm{PL}\langle8\rangle}\otimes\Z_s,
\qquad
\Z_s\xrightarrow{\ \times7\ }\Z_s.
\end{equation}
\end{lemma}

\begin{proof}
The Thom map $T$ is a $7$-equivalence. Hence it induces
isomorphisms on coefficient groups through degree six and a
surjection in degree seven. The table follows from the smooth
coefficient calculation and
$\Omega_7^{\mm O\langle8\rangle}=0$.

To compute the degree-eight map, we first work over
$
R=\Z[1/6].$

By Lemma~\ref{BO8toBPL8times7} and the Thom isomorphisms,
\[
H_8(\mm{MO}\langle8\rangle;R)
\longrightarrow
H_8(\mm{MPL}\langle8\rangle;R)
\]
is multiplication by $7$ between two copies of $R$.

For the Thom spectrum $A=\mm{MO}\langle8\rangle$ or $\mm{MPL}\langle8\rangle$, consider the AHSS
\[
H_p(A;\pi_q(S^0)\otimes R)
\Longrightarrow
\pi_{p+q}(A)\otimes R.
\]
Since $H_i(A;R)=0$ for $0<i<8$, the only differential entering
$(0,7)$ is $d_8$. Moreover,
\[
\pi_7(S^0)\otimes R\cong\Z_5,
\qquad
\pi_8(S^0)\otimes R=0.
\]
Since $\pi_7(A)=0$, $d_8$ is surjective. Thus naturality gives
\[
\begin{tikzcd}[column sep=small,row sep=large]
0\arrow[r]
&
\Omega_8^{\mm O\langle8\rangle}\otimes R
\arrow[r]\arrow[d,"T_*"]
&
R\arrow[r,"d_8"]\arrow[d,"\times7"]
&
\Z_5\arrow[r]\arrow[d,equal]
&
0
\\
0\arrow[r]
&
\Omega_8^{\mm{PL}\langle8\rangle}\otimes R
\arrow[r]
&
R\arrow[r,"\tilde d_8"]
&
\Z_5\arrow[r]
&
0.
\end{tikzcd}
\]
Both kernels of $d_8$ and $\tilde d_8$ are $5R\subset R$.
Therefore the left vertical map, expressed in the corresponding
generators, is also multiplication by $7$.

Finally, since $\gcd(s,6)=1$,
$
R\otimes\Z_s\cong\Z_s
$.
Tensoring the homomorphism
$$T_\ast:\Omega_8^{\mm O\langle8\rangle}\otimes R\to \Omega_8^{\mm {PL}\langle8\rangle}\otimes R$$ 
with $\Z_s$ proves
\eqref{eq:PL-degree8-times7}.
\end{proof}

\subsection{The filtration-five contribution}
\quad

  We now consider the  kernel on a certain subgroup. 

\begin{theorem}\label{7ordersubgroupinkerT}
Let $C_s\subset F_{5,8}$ be the subgroup of elements whose orders
are coprime to $6$. Then $C_s\cong\Z_s$, and
\[
\ker(T_*|_{C_s})\cong\Z_{\gcd(s,7)},
\qquad
T_*(C_s)\cong\Z_{s/\gcd(s,7)}.
\]

More generally,
\begin{equation}\label{eq:PL-kernel-through-five}
\ker(T_*|_{ F_{5,8}})
=
 F_{4,9}\oplus\ker(T_*|_{C_s}).
\end{equation}
\end{theorem}

\begin{proof}
The calculation in Section~\ref{pprimary} gives
\[
E^\infty_{5,8}
=
E^2_{5,8}
\cong\Z_s.
\]
The lower subgroup $ F_{4,9}$ is finite with only $2$- and
$3$-primary torsion. So,
\begin{equation}\label{F58_split}
	F_{5,8}\cong 
 F_{4,9}\oplus C_s.
\end{equation}

By Lemma~\ref{PL8Zs} and the homology of $\M\eta$, we have
\[
\widetilde E^2_{5,8}
=
H_5(\mm M\eta;\Omega_8^{\mm{PL}\langle8\rangle})
\cong\Z_s,
\]
and the map
$
E^2_{5,8}\longrightarrow\widetilde E^2_{5,8}$
is multiplication by $7$. 

By Lemma~\ref{PL8Zs} and the homology of $\mm M\eta$, the source of every possible incoming differential of $\widetilde E^r_{5,8}$ is zero or finite with only
$2$- and $3$-primary torsion.
Therefore, all incoming differentials for $\widetilde E^r_{5,8}$ are zero. 

Since every class of $E^2_{5,8}$ is a permanent cycle,
naturality implies that its image under multiplication
by $7$ supports no outgoing differential.
Therefore, the subgroup
$
7\Z_s\subset\widetilde E^2_{5,8}$
survives to $\widetilde E^\infty_{5,8}$.

By Lemma \ref{PL9primary} and \eqref{PL49},  
$
\widetilde { F}_{4,9}$ only has $2$-prime component. Thus,
$$\widetilde { F}_{5,8}\cong \widetilde { F}_{4,9}\oplus \widetilde {E}^\infty_{5,8}$$
Combing this with 
\eqref{F58_split} gives
$$\mm{Ker}(T_\ast)\cong \mm{Ker}(T_\ast|_{F_{4,9}})\oplus \mm{Ker}(T_\ast|_{C_s}).$$
Note that $T_*$ vanishes on $F_{4,9}$;
$T_*|_{C_s}$ has exactly the kernel of multiplication by $7$
on $\Z_s$. This finishes the proof.
\end{proof}

\begin{remark}\label{rem:PL-seven-local}
The multiplication-by-$7$ formula holds for every $s$ coprime to $6$,
not only when $7\mid s$.
It is injective when $7\nmid s$; when $7\mid s$, it kills exactly
one subgroup of order $7$, even if a higher power of $7$ divides $s$.

The proof only requires the survival of $7\Z_s$ on the PL side;
it does not assume that every class of $\widetilde E^2_{5,8}$ survives.
\end{remark}

\subsection{The comparison condition involving
\texorpdfstring{$7\sigma_4$}{7 sigma4}}
\quad

We now relate the filtration-five calculation to the characteristic
invariants in Section~\ref{Zs}.
For this interpretation, retain the assumption $\gcd(s,5)=1$
used there to identify the smooth $\Z_s$-coordinate by the integral
Hurewicz homomorphism.
This extra assumption is not needed for the preceding statements
about $T_*$.

\begin{corollary}\label{cor:PL-seven-sigma4}
Let $\q,\q'$ be Bazaikin tuples satisfying
\[
\sigma_2(\q)=\sigma_2(\q'),
\qquad
\sigma_3(\q)=\sigma_3(\q')=8s>0,
\qquad
\gcd(s,5)=1,
\]
and $\sigma_5(\q)\equiv \sigma_5(\q')\pmod s$.

Let $b_s(\q)$ and $b_s(\q')$ denote the sums of the primary
components, at primes dividing $s$, of the normal $5$-smoothings of $\mathcal B_\q $ and $\mathcal B_{\q'}$.
Then
\begin{equation}\label{eq:PL-seven-sigma4}
T_*b_s(\q)=T_*b_s(\q')
\quad\Longleftrightarrow\quad
7\sigma_4(\q)\equiv7\sigma_4(\q')\pmod s.
\end{equation}
\end{corollary}

\begin{proof}
The case $s=1$ is immediate, so assume $s>1$.

If $\gcd(s,5)=1$, by Lemma  \ref{E_58130_odd} and Proposition \ref{ZsHure}, we have
$$\bord_{13}(\eta)\cong \bord_{13}(\eta)_{(2)}\oplus \bord_{13}(\eta)_{(3)}\oplus \Z_5\oplus \Z_s\oplus \Z_s.$$ 

By Lemma~\ref{130Zs}, the filtration-$13$ $\Z_s$-coordinate is
determined by the linking invariant.
Since these invariants agree,
\[
b_s(\q)-b_s(\q')\in C_s\cong\Z_s.
\]
Fix a generator $g_s$ of $C_s$ and write
\[
b_s(\q)-b_s(\q')
=
d\,g_s,
\qquad d\in\Z_s.
\]
Theorem~\ref{7ordersubgroupinkerT} gives
\begin{equation}\label{eq:PL-coordinate-seven}
T_*b_s(\q)=T_*b_s(\q')
\quad\Longleftrightarrow\quad
7d=0\quad\text{in }\Z_s.
\end{equation}

Write
$
p_2(\mathcal B_{\q})
=
p_{2,\q}\,x^4,$ $p_{2,\q}\in\Z_s,$
and use the analogous notation for $\q'$. 
The characteristic-number calculation \eqref{p2deterZs58} gives
$$b_s(\q)-b_s(\q')= (p_{2,\q}\ell k(\mathcal B_\q )-p_{2,\q'}\ell k(\mathcal B_{\q'} ))\mathfrak u g_s$$
where $\mathfrak u$ is a unit of $\Z_s$. Then 
$$d=(p_{2,\q}\ell k(\mathcal B_\q )-p_{2,\q'}\ell k(\mathcal B_{\q'} ))\mathfrak u \in \Z_s.$$
Since $\sigma_5(\q)\equiv \sigma_5(\q')\pmod s$ is a unit of $\Z_s$ (cf. \cite[Theorem~2.4]{FL}), $$\ell k(\mathcal B_\q )\equiv\ell k(\mathcal B_{\q'})\pmod s $$
is also a unit of $\Z_s$. Hence, $7d=0\in \Z_s$ is equivalent to
$$7p_{2,\q}-7p_{2,\q'}=0\in \Z_s$$ 

By \cite[Theorem~3.1]{FL},
$
8p_{2,\q}
=
3\sigma_2(\q)^2-\sigma_4(\q)
\in \Z_s.$
Since $\sigma_2(\q)=\sigma_2(\q')$ and $8$ is a unit of $\Z_s$, $7d=0\in \Z_s$ is equivalent to
\[
7\sigma_4(\q)\equiv 7\sigma_4(\q')\pmod s.
\]
This finishes the proof.
\end{proof}


\section{Proofs of main Theorems}\label{mainproof}

\begin{proof}[Proof of Theorem \ref{classify}]
Recall the AHSS for $\Omega_{\ast}^{\mm{O}\langle 8\rangle}(\eta)$
$$E^2_{p,q}=H_p(\mm{M}\eta;\pi_q(\mm{MO}\langle 8\rangle))\Longrightarrow \Omega_{p+q}^{\mm{O}\langle 8\rangle}(\eta)$$
Let ${F}_{p,\ast}$ be its associated filtration.

Recall $\gcd(s,6)=1$.
Suppose $\gcd(s,5)=1$ and
\[
\sigma_i(\q)=\sigma_i(\q')\quad(i=2,3),
\qquad\sigma_5(\q)\equiv\sigma_5(\q')\pmod s.
\]
By Proposition \ref{Bordism_compare_lowF_85}, the differential of the normal $5$-smoothings satisfies
\[
[\mathcal B_{\q},(g,f)]
-[\mathcal B_{\q'},(g',f')]
\quad\text{in }
F_{8,5}\subset \Omega_{13}^{\mm{O}\langle 8\rangle}(\eta).
\]
Since $E^\infty_{p,13-p}=0$ for $p=8,7,6$ (cf. Table \ref{E2AHSSforMO8Meta} and Lemma \ref{E_infty_76}), we have
$$F_{8,5}=F_{5,8}.$$
 By $F_{5,8}\cong F_{4,9}\oplus E^\infty_{5,8}$ (cf. \eqref{F58_split}), the differential has the representative
 \[
[\mathcal B_{\q},(g,f)]
-[\mathcal B_{\q'},(g',f')]=(\alpha_{4,9},\alpha_{5,8})\in F_{4,9}\oplus E^\infty_{5,8}
=
F_{5,8}.
\]

By Corollary \ref{cor:PL-seven-sigma4}, if $7\sigma_4(\q)\equiv 7\sigma_4(\q')\pmod s$, the forgetful map
\[
T_*:
\Omega_{13}^{\mm O\langle8\rangle}(\eta)
\longrightarrow
\Omega_{13}^{\mm{PL}\langle8\rangle}(\eta)
\]
maps $\alpha_{5,8}$ to zero. By Theorem \ref{top13}, $T_\ast(\alpha_{4,9})=0$. Therefore,
$$[\mathcal B_{\q},(T\circ g,f)]=[\mathcal B_{\q'},(T\circ g',f')]\in \Omega_{13}^{\mm{PL}\langle8\rangle}(\eta).$$
By Proposition \ref{PLhcobor}, $\mathcal B_{\q}$ is PL homeomorphic to $\mathcal B_{\q'}$.

As is well known, the cohomology ring and the link form are homotopy invariants. Furthermore, the rational Pontryagin classes are homeomorphic invariants. Thus, the first Pontryagin class is a homeomorphic invariant by $H^4(\mathcal{B}_{\bold{q}})=\mathbb{Z}$. By \cite{Sharma}, $7p_2$ is a homeomorphic invariant. By the relations between these invariants and the $\sigma_i$ (i=2,3,4,5), the necessary holds. 
\end{proof}

%
%

\begin{proof}[Proof of Theorem \ref{diff}] 
If $\mathcal B_\q $ is homeomorphic to $\mathcal B_{\q'}$, then they are PL homeomorphic. Thus their normal $5$-smoothings have the same image in $\Omega_{13}^{\mm{PL}\langle 8\rangle}(\eta)$.  
 It is easy to check that
$$\bord_{13}(\eta)/F_{8,5}\to \Omega_{13}^{\mm{PL}\langle 8\rangle}(\eta)/\widetilde{F}_{8,5}$$
is injective. Hence, 
\[
[\mathcal B_{\q},(g,f)]
-[\mathcal B_{\q'},(g',f')]
\quad\text{in }
F_{8,5}\subset \Omega_{13}^{\mm{O}\langle 8\rangle}(\eta).
\]

Note that $|F_{4,9}|\le 12$.
From
Theorem \ref{top13}, \ref{7ordersubgroupinkerT}, there is a subgroup of order $12\cdot 7^n$ at most of $\Omega_{13}^{\mm{O}\langle 8\rangle}(\eta)$ which vanishes under the homomorphism $T_\ast:\Omega_{13}^{\mm{O}\langle 8\rangle}(\eta)\to \Omega_{13}^{\mm{PL}\langle 8\rangle}(\eta)$ where $n=0$ or 1.
 So Theorem \ref{diff} follows by Theorem \ref{hcobordism}.
\end{proof}

\appendix

\section{Surgery obstruction}\label{Surgeryobstruction}

In this section, we apply Kreck's modified surgery theory
\cite{K} to show that bordant normal
$5$-smoothings of Bazaikin spaces give rise to an
$h$-cobordism.

\begin{theorem}\label{hcobordism}
Suppose that the smooth (resp. PL) normal
$5$-smoothings of $\mathcal B_{\bold q}$ and
$\mathcal B_{\bold q'}$ are bordant in 
$\Omega_{13}^{\mm O\langle8\rangle}(\eta)$
(resp.\ $\Omega_{13}^{\mm{PL}\langle8\rangle}(\eta)$),
then they are smoothly (resp.\ PL) $h$-cobordant.
\end{theorem}

\begin{proof}
We first consider the smooth category.
Write
\[
M_0=\mathcal B_{\bold q},
\qquad
M_1=\mathcal B_{\bold q'},
\qquad
\mathscr B_{s,k}=\mm{BO}\langle8\rangle\times B.
\]
Let $(W,\phi)$ be a normal $\mathscr B$-bordism
between the smooth normal
$5$-smoothings of
$M_0$ and $M_1$.
By surgery below the middle dimension
\cite[Proposition~4]{K}, we may assume that
\[
\phi:W\longrightarrow\mathscr B
\]
is a $7$-equivalence. Consequently,
\[
H_i(W,M_j;\mathbb Z)=0,
\qquad i\le6,\quad j=0,1.
\]

The cohomology of the boundary components,
together with the long exact sequences of the
pairs and Poincare--Lefschetz duality, gives
\[
H_7(W;\mathbb Z)
\cong
\mathbb Z^{2r}\oplus\mathbb Z_s
\]
for some $r\ge0$, where
$s=|\sigma_3(\bold q)|/8$.
Moreover,
\[
H_7(W,M_j;\mathbb Z)\cong\mathbb Z^{2r},
\qquad j=0,1.
\]
The middle-dimensional intersection pairing
induces a unimodular skew-symmetric form
\[
\lambda:A\times A\longrightarrow\mathbb Z,
\qquad
A=H_7(W;\mathbb Z)/\operatorname{Tor}.
\]
Thus $(A,\lambda)$ is a symplectic lattice.

Consider the homomorphism
\[
\phi_*:
H_7(W;\mathbb Z)\longrightarrow
H_7(\mathscr B;\mathbb Z).
\]
By Proposition~\ref{HBZ} and the Hurewicz theorem,
\[
H_7(\mathscr B;\mathbb Z)
\cong\mathbb Z_2\oplus\mathbb Z_s.
\]
The map $\phi_\ast$ restricts to an isomorphism
on the $\mathbb Z_s$-summands. Indeed, this follows
from the naturality of the universal coefficient
theorem and the surjection
\[
H^8(B;\mathbb Z)
\longrightarrow H^8(M_0;\mathbb Z)
\cong\mathbb Z_s,
\]
which sends $x^4$ to a generator.

After choosing the free summand of $H_7(W)$
appropriately, the restriction of $\phi_\ast$ to
$A$ is therefore a homomorphism
\[
\psi:A\longrightarrow\mathbb Z_2.
\]

By the following morphism between exact sequences
 \[
\xymatrix@C=.3cm{
\pi_8(\mm{BO}\langle 8\rangle) =\pi_8(\mathscr{ B})\ar[r]^-{} \ar[d]^-{\cong}& \pi_8(\mathscr{ B},W) \ar[r]^-{} \ar[d]^-{\cong}& \pi_7(W)  \ar[r]^-{\phi_\#} \ar[d]^-{h_\ast}&\pi_7(\mathscr{ B})=0 \ar[d]^-{}&\\
 H_8(\mm{BO}\langle 8\rangle)= H_8(\mathscr{ B})\ar[r]^-{} &H_8(\mathscr{ B},W) \ar[r]^-{}& H_7(W)\ar[r]^-{\phi_\ast} &H_7(\mathscr{ B})\ar[r]^-{}& 0
}
\]
we identify the
middle-dimensional surgery kernel $\mm{Ker}\phi_\#$ with
\[
\ker\psi\subset A.
\]

If $\psi=0$, we may choose any
symplectic basis
of $A$.
If $\psi\ne0$, the symplectic linear algebra
over $\mathbb F_2$ allows us to choose  a symplectic
basis
$$\{a_1,b_1,\cdots,a_r,b_r\},\quad a_i\cup b_j=\delta_{ij}$$
 satisfying
\[
\psi(a_i)=0,\qquad 1\le i\le r.
\]
In either case, the subgroup
\[
U=\langle a_1,\ldots,a_r\rangle
\subset\ker\psi
\]
satisfies
$
\lambda|_{U\times U}=0.$
And the intersection pairing induces an
isomorphism
\[
U\xrightarrow{\ \cong\ }(A/U)^*.
\]

It remains to consider the quadratic part of
Kreck's surgery obstruction.
Since
$
\pi_1(\mathscr B)=0$
and the universal class $w_8$ evaluates
nontrivially on
\[
\pi_8(\mm{BO}\langle8\rangle)
\subset\pi_8(\mathscr B),
\]
we are in the exceptional case $q=7$ of
\cite[Section~5]{K}.
In this case the modified quadratic refinement
$\widetilde\mu$ takes values in
\[
\mathbb Z/\langle2,1\rangle=0.
\]
Hence $\widetilde\mu$ vanishes identically.

It follows that $U$ satisfies the elementary
conditions for Kreck's modified surgery
obstruction
\[
\theta(W,\phi)\in\widetilde l_{14}(0).
\]
By \cite[Theorem~3]{K}, $(W,\phi)$ is normally
$\mathscr B$-bordant relative to its boundary
to a smooth $s$-cobordism. In particular,
$M_0$ and $M_1$ are smoothly $h$-cobordant.

For the PL category, replace
$\mm{BO}\langle8\rangle$ by
$\mm{BPL}\langle8\rangle$ and use the
corresponding PL normal $5$-smoothings.
The same homological argument applies.
The nontrivial evaluation of $w_8$ persists
under the natural map
\[
\mm{BO}\langle8\rangle
\longrightarrow\mm{BPL}\langle8\rangle.
\]
Thus the corresponding PL modified surgery
obstruction is again elementary.
The PL version of modified surgery
 yields a
PL $h$-cobordism.
\end{proof}

\section{Algebraic Atiyah-Hirzebruch spectral sequence}\label{AAHSSmod2}

We use the algebraic Atiyah--Hirzebruch spectral sequence
(AAHSS) to compute the relevant
$\Ext_{\mathscr A_2}$-groups of
$H^*(\mm M\eta;\Z_2)$.
Throughout this section, all tensor products
are taken over $\Z_2$ unless otherwise indicated.

We first recall some representatives in the bar complex.
Let $\mathscr A_*$ denote the dual Steenrod algebra.
By \cite[Theorem~2.4.1]{Adams1960}, the classes
\[
h_i\in
\Ext_{\mathscr A_2}^{1,2^i}(\Z_2,\Z_2),
\qquad i=0,1,2,
\]
are represented in the cobar complex by
$\xi_1^{2^i}\in\mathscr A_*$,
dual to $\Sq^{2^i}\in\mathscr A$.

Using the graded duality between the bar and cobar complexes,
we may choose the following bar representatives for
the corresponding classes.

\begin{center}
\renewcommand{\arraystretch}{1.3}
\begin{tabular}{|c|l|}
\hline
Class & Representative in the bar complex\\
\hline
$h_0^s$, $s\ge1$
&
$[\Sq^1|\cdots|\Sq^1]$
\\
\hline
$h_1^2$
&
$[\Sq^2|\Sq^2]+[\Sq^3|\Sq^1]$
\\
\hline
$h_0h_2$
&
$[\Sq^1|\Sq^4]+[\Sq^2|\Sq^3]+[\Sq^4|\Sq^1]$
\\
\hline
$h_2^2$
&
$[\Sq^4|\Sq^4]+[\Sq^7|\Sq^1]+[\Sq^6|\Sq^2]$
\\
\hline
\end{tabular}
\end{center}

Similarly, a representative for $h_0^2h_2$ is
\[
\begin{aligned}
&
[\Sq^4|\Sq^1|\Sq^1]
+[\Sq^3|\Sq^1|\Sq^2]
+[\Sq^2|\Sq^2|\Sq^2]
\\
&\quad
+[\Sq^1|\Sq^1|\Sq^4]
+[\Sq^1|\Sq^2|\Sq^3]
+[\Sq^1|\Sq^4|\Sq^1]
\\
&\quad
+[\Sq^2|\Sq^3|\Sq^1].
\end{aligned}
\]

Put
$
N=H^*(\mm M\eta;\Z_2),$
$N^n=H^n(\mm M\eta;\Z_2),$
and let
\[
N_n=\bigoplus_{i\ge n}N^i
\]
be the filtration by cohomological degree.

The associated AAHSS is
\[
E_1^{s,m,n}
=
\Tor_{\mathscr A_2}^{s,s+m}(\Z_2,\Z_2)
\otimes N^n
\Longrightarrow
\Tor_{\mathscr A_2}^{s,s+m+n}(\Z_2,N),
\]
with differentials
\begin{equation}\label{eq:AAHSS-differential}
d_r:
E_r^{s,m,n}
\longrightarrow
E_r^{s-1,m+1-r,n+r}.
\end{equation}

Here $s$ denotes the homological degree,
$m$ the stem of the coefficient class,
and $n$ the cohomological filtration.
In particular, the total stem is $m+n$.


We restrict attention to the range
\[
11\le m+n\le14.
\]

\subsubsection*{The first two pages}

The nonzero $E_1$-terms in this range are listed
in Table~\ref{E1smn}.
The differential $d_1$ is determined by the
$\Sq^1$-action on $N$.

\begin{table}[htbp]
\centering
\small
\renewcommand{\arraystretch}{1.2}
\caption{Nonzero $E_1^{s,m,n}$-terms for $12\le m+n\le14$}
\label{E1smn}

\begin{tabular}{|c|l|}
\hline
$(s,m,n)$ & Generators\\
\hline

$(s,0,11)$
&
$h_0^s ux^2U,\quad
h_0^svxU,\quad
h_0^s\Sq^4uU$
\\
\hline

$(s,0,12)$
&
$h_0^sx^2\Sq^1uU,\quad
h_0^sx\Sq^2\Sq^1uU,\quad
h_0^s\Sq^5uU$
\\
\hline

$(s,0,13)$
&
$h_0^svx^2U,\quad
h_0^sx\Sq^4uU,\quad
h_0^s\Sq^6uU$
\\
\hline

$(s,0,14)$
&
$\begin{aligned}
&h_0^su^2U,\quad
h_0^sx^2\Sq^2\Sq^1uU,\\
&h_0^sx\Sq^5uU,\quad
h_0^s\Sq^6\Sq^1uU
\end{aligned}$
\\
\hline

$(1,1,11)$
&
$h_1ux^2U,\quad h_1vxU,\quad h_1\Sq^4uU$
\\
\hline

$(1,1,12)$
&
$h_1x^2\Sq^1uU,\quad
h_1x\Sq^2\Sq^1uU,\quad
h_1\Sq^5uU$
\\
\hline

$(1,1,13)$
&
$h_1vx^2U,\quad
h_1x\Sq^4uU,\quad
h_1\Sq^6uU$
\\
\hline

$(2,2,10)$
&
$h_1^2\Sq^2\Sq^1uU,\quad
h_1^2x\Sq^1uU$
\\
\hline

$(2,2,11)$
&
$h_1^2ux^2U,\quad
h_1^2vxU,\quad
h_1^2\Sq^4uU$
\\
\hline

$(2,2,12)$
&
$h_1^2x^2\Sq^1uU,\quad
h_1^2x\Sq^2\Sq^1uU,\quad
h_1^2\Sq^5uU$
\\
\hline

$(s+1,3,9)$
&
$h_0^sh_2uxU,\quad h_0^sh_2vU$
\\
\hline

$(s+1,3,10)$
&
$h_0^sh_2\Sq^2\Sq^1uU,\quad
h_0^sh_2x\Sq^1uU$
\\
\hline

$(s+1,3,11)$
&
$h_0^sh_2ux^2U,\quad
h_0^sh_2vxU,\quad
h_0^sh_2\Sq^4uU$
\\
\hline

$(2,6,7)$ & $h_2^2uU$\\
\hline
$(2,6,8)$ & $h_2^2\Sq^1uU$\\
\hline
$(3,8,4)$ & $c_0x^2U$\\
\hline

$(s+4,8,4)$
&
$h_0^s\omega x^2U$
\\
\hline

$(4,9,4)$ & $h_1c_0x^2U$\\
\hline
$(5,9,2)$ & $h_1\omega xU$\\
\hline
$(5,9,4)$ & $h_1\omega x^2U$\\
\hline

$(6,10,2)$ & $h_1^2\omega xU$\\
\hline
$(6,10,4)$ & $h_1^2\omega x^2U$\\
\hline
$(s+5,11,0)$
&
$h_0^sh_2\omega U$
\\
\hline

$(s+5,11,2)$
&
$h_0^sh_2\omega xU$
\\
\hline

$(s+3,12,0)$ & $h_0^s\tau U$\\
\hline
$(s+3,12,2)$ & $h_0^s\tau xU$\\
\hline

$(s+4,14,0)$ & $h_0^sd_0U$\\
\hline
\end{tabular}

\medskip
\parbox{0.94\textwidth}{\small
In the rows involving $h_0^s$, we take $s\ge0$
unless otherwise specified.
For the $h_2$ and $h_2\omega$-families, $0\le s\le2$. We only list the $E_1^{s,m,11-m}$-terms with $s\ge 5$.
}
\end{table}

Before passing to the $E_2$-page, we record a useful
observation about classes supported in cohomological
degrees $2$ and $4$.

\begin{lemma}\label{pert}
Every element of
$
E_1^{s,m,2}
$ or
$E_1^{s,m,4}$
supports no nonzero outgoing differential in the AAHSS.
\end{lemma}

\begin{proof}
By the Steenrod operations on $N$,
\[
\Sq^I(xU)=0,
\qquad
\Sq^I(x^2U)=0
\]
for every nonempty admissible sequence $I$
of positive total degree.

Consider a representative
\[
\sum[a_1|\cdots|a_s]x^2U
\in E_1^{s,m,4}.
\]
Since all positive-degree Steenrod operations
annihilate $x^2U$, the representative lifts
to a cycle in the filtered bar complex $B(N)$.
Consequently, it supports no outgoing AAHSS
differential.

The same argument applies to
$
\sum[a_1|\cdots|a_s]xU
\in E_1^{s,m,2}.$
\end{proof}

\medskip


The differential $d_1$ is computed from the
$\Sq^1$-action on $N$, using the bar representatives
listed above. The resulting nonzero $E_2$-terms
are given in Table~\ref{E2smn}.

\begin{table}[htbp]
\centering
\small
\renewcommand{\arraystretch}{1.2}
\caption{Nonzero $E_2^{s,m,n}$-terms for $11\le m+n\le14$}
\label{E2smn}

\begin{tabular}{|c|l|}
\hline
$(s,m,n)$ & Generators\\
\hline

$(0,0,13)$
&
$vx^2U,\quad x\Sq^4uU,\quad\Sq^6uU$
\\
\hline

$(0,0,14)$ & $\Sq^6\Sq^1uU$\\
\hline

$(1,1,11)$
&
$h_1ux^2U,\quad h_1vxU,\quad h_1\Sq^4uU$
\\
\hline

$(1,1,12)$
&
$h_1x^2\Sq^1uU,\quad
h_1x\Sq^2\Sq^1uU,\quad
h_1\Sq^5uU$
\\
\hline

$(1,1,13)$
&
$h_1vx^2U,\quad
h_1x\Sq^4uU,\quad
h_1\Sq^6uU$
\\
\hline

$(1,3,9)$ & $h_2uxU,\quad h_2vU$\\
\hline

$(1,3,11)$
&
$h_2ux^2U,\quad h_2vxU,\quad h_2\Sq^4uU$
\\
\hline

$(2,2,10)$
&
$h_1^2\Sq^2\Sq^1uU,\quad h_1^2x\Sq^1uU$
\\
\hline

$(2,2,11)$
&
$h_1^2ux^2U,\quad
h_1^2vxU,\quad
h_1^2\Sq^4uU$
\\
\hline

$(2,2,12)$
&
$h_1^2x^2\Sq^1uU,\quad
h_1^2x\Sq^2\Sq^1uU,\quad
h_1^2\Sq^5uU$
\\
\hline

$(2,6,7)$ & $h_2^2uU$\\
\hline
$(2,6,8)$ & $h_2^2\Sq^1uU$\\
\hline

$(3,3,10)$
&
$h_0^2h_2\Sq^2\Sq^1uU,\quad
h_0^2h_2x\Sq^1uU$
\\
\hline

$(3,8,4)$ & $c_0x^2U$\\
\hline

$(s+4,8,4)$ & $h_0^s\omega x^2U$\\
\hline

$(4,9,4)$ & $h_1c_0x^2U$\\
\hline
$(5,9,4)$ & $h_1\omega x^2U$\\
\hline

$(6,10,2)$ & $h_1^2\omega xU$\\
\hline
$(6,10,4)$ & $h_1^2\omega x^2U$\\
\hline

$(s+5,11,0)$ & $h_0^sh_2\omega U$\\
\hline

$(s+5,11,2)$ & $h_0^sh_2\omega xU$\\
\hline

$(s+3,12,0)$ & $h_0^s\tau U$\\
\hline
$(s+3,12,2)$ & $h_0^s\tau xU$\\
\hline

$(s+4,14,0)$ & $h_0^sd_0U$\\
\hline
\end{tabular}

\medskip
\parbox{0.94\textwidth}{\small
The ranges of $s$ are the same as those in
Table~\ref{E1smn}.
}
\end{table}

\subsubsection*{Classes surviving in the AAHSS}
We use the bar and cobar complexes to establish the survival of several
families.

Recall that
$
N_*=H_*(\mm M\eta;\Z_2)$
is dual to $N$ in each degree.
Let $B(N)$ and $cB(N_*)$ denote the bar and
cobar complexes, with differentials
$\partial$ and $\delta$, respectively.


\begin{lemma}\label{pert2}
In the AAHSS,
\[
E_\infty^{s+5,11,2}
=
E_1^{s+5,11,2},
\qquad 0\le s\le2,
\]
\[
E_\infty^{s+3,12,2}
=
E_1^{s+3,12,2},
\qquad s\ge0.
\]
\end{lemma}

\begin{proof}
We first consider
$
h_0^sh_2\omega xU
\in E_1^{s+5,11,2}.$

In the cobar complex $cB(N_*)$,
the cobar differential has the form
\[
\delta(h_0^sh_2\omega xU)
=
h_0^sh_2\omega h_1U.
\]

By the relations in
$\Ext_{\mathscr A_2}^{*,*}(\Z_2,\Z_2)$
(cf. Fig.~\ref{1Fig}), the class
$
h_0^sh_2\omega h_1$
is zero. Hence there exists a cochain
$\alpha\in cB(\Z_2)$ such that
$
\delta\alpha
=
h_0^sh_2\omega h_1.$

Consequently,
$
h_0^sh_2\omega xU+\alpha U$
is a cocycle in $cB(N_*)$. Its dual in the  bar complex cannot be
killed by an incoming differential. Moreover, its dual has the leading term $h_0^sh_2\omega xU$ in the filtration $2$ of the AAHSS. 
This implies that the class $h_0^sh_2\omega xU\in E_r^{s+5,11,2}$ cannot be
killed by an incoming differential in the AAHSS.

On the other hand, Lemma~\ref{pert}
excludes all outgoing differentials.
Therefore
$
E_\infty^{s+5,11,2}
=
E_1^{s+5,11,2}.$

The same cobar lifting argument applies to
$
h_0^s\tau xU
\in E_1^{s+3,12,2}.$ Then we obtain that $h_0^s\tau xU$ cannot be killed by an incoming  differential in the AAHSS.
 Together with Lemma~\ref{pert},
this proves
$
E_\infty^{s+3,12,2}
=
E_1^{s+3,12,2}.$
\end{proof}

\begin{lemma}\label{pert3}
The following classes survive to $E_\infty$:
\[
E_\infty^{3,8,4}
=
E_1^{3,8,4},
\]
\[
E_\infty^{s+4,8,4}
=
E_1^{s+4,8,4},
\qquad s\ge0,
\]
\[
E_\infty^{s,9,4}
=
E_1^{s,9,4},
\qquad s=4,5,
\]
\[
E_\infty^{6,10,4}
=
E_1^{6,10,4}.
\]
\end{lemma}

\begin{proof}
We first consider
$
c_0x^2U\in E_1^{3,8,4}.$
By the Steenrod action on $N$,
\[
\Sq^4U=0,
\qquad
\Sq^2(xU)=0.
\]
The corresponding cobar calculation gives
\[
\delta(c_0x^2U)=0.
\]
Thus $c_0x^2U$ represents a cocycle
in $cB(N_*)$ which pairs nontrivially
with the corresponding bar representative.
This excludes incoming differentials in the AAHSS.
By Lemma~\ref{pert}, there are no
outgoing differentials.
Therefore
\[
E_\infty^{3,8,4}
=
E_1^{3,8,4}.
\]

The same argument applies to
$
h_0^s\omega x^2U,
h_1c_0x^2U,
h_1\omega x^2U,
h_1^2\omega x^2U.$
\end{proof}

\begin{lemma}\label{pert4}
The following equalities hold:
\[
E_\infty^{s+5,11,0}
=
E_1^{s+5,11,0},
\quad 0\le s\le 1,\quad E_3^{7,11,0}
=
0,
\]
\[
E_\infty^{s+3,12,0}
=
E_1^{s+3,12,0},
\qquad s\ge0,
\]
\[
E_\infty^{s+4,14,0}
=
E_1^{s+4,14,0},
\qquad 0\le s\le2.
\]
\end{lemma}

\begin{proof}
By Lemma~\ref{pert2} and the bidegrees
of the differentials, the classes
\[
h_2\omega U,\quad h_0h_2\omega U,\quad h_0^s\tau U,
\quad
h_0^sd_0U
\]
survive at least to the $E_3$-page. Since $h_0^2h_2\omega U=h_1^3\omega U$ and $\Sq^2U=x$,
$$E_3^{7,11,0}
=
0.$$

Recall that
$\Sq^I U=0$
for every nonempty admissible sequence
$I$ of total degree at least $3$.
Using this vanishing in the bar complex
shows that
\[
h_0^s\tau U,
\qquad
h_0^sd_0U
\]
admit cycle representatives in $B(N)$.

Moreover, these classes lie in the
lowest cohomological filtration $n=0$,
so no incoming differential can hit them.
So, they survive to $E_\infty$.
\end{proof}

\subsubsection*{The higher differentials}

We now complete the calculation in
the range $11\le m+n\le14$.

Inspection of the bidegrees in
Table~\ref{E2smn} shows that
$d_3=0$
in the range under consideration.
Hence
$
E_4^{s,m,n}=E_3^{s,m,n}.$

The differential $d_2$ is determined by the
$\Sq^2$-action on $N$.
After making the computation, we obtain the
nonzero $E_4$-terms listed in
Table~\ref{E4smn}.

For the families already proved to
survive to $E_\infty$, we simply
record their generators.

\begin{table}[htbp]
\centering
\small
\renewcommand{\arraystretch}{1.2}
\caption{Nonzero $E_4^{s,m,n}$-terms for $12\le m+n\le14$}
\label{E4smn}

\begin{tabular}{|c|l|}
\hline
$(s,m,n)$ & Generators\\
\hline

$(0,0,13)$
&
$vx^2U,\quad x\Sq^4uU=\Sq^6uU$
\\
\hline

$(1,1,12)$ & $h_1x^2\Sq^1uU$\\
\hline
$(1,1,13)$ & $h_1vx^2U$\\
\hline

$(1,3,9)$ & $h_2uxU,\quad h_2vU$\\
\hline

$(1,3,11)$
&
$h_2ux^2U,\quad h_2vxU,\quad h_2\Sq^4uU$
\\
\hline

$(2,2,11)$
&
$h_1^2ux^2U,\quad h_1^2vxU$
\\
\hline

$(2,6,7)$ & $h_2^2uU$\\
\hline
$(2,6,8)$ & $h_2^2\Sq^1uU$\\
\hline

$(3,3,10)$
&
$h_0^2h_2(\Sq^2\Sq^1u+x\Sq^1u)U$
\\
\hline

$(3,8,4)$ & $c_0x^2U$\\
\hline

$(s+4,8,4)$
&
$h_0^s\omega x^2U$
\\
\hline

$(4,9,4)$ & $h_1c_0x^2U$\\
\hline
$(5,9,4)$ & $h_1\omega x^2U$\\
\hline

$(6,10,4)$ & $h_1^2\omega x^2U$\\
\hline

$(5,11,0)$
&
$h_2\omega U$
\\
\hline

$(6,11,0)$
&
$h_0h_2\omega U$
\\
\hline

$(s+5,11,2)$
&
$h_0^sh_2\omega xU$
\\
\hline

$(s+3,12,0)$ & $h_0^s\tau U$\\
\hline
$(s+3,12,2)$ & $h_0^s\tau xU$\\
\hline

$(s+4,14,0)$ & $h_0^sd_0U$\\
\hline
\end{tabular}

\medskip
\parbox{0.94\textwidth}{\small
The ranges of $s$ are as in Table~\ref{E1smn}.
}
\end{table}

The next lemma determines the remaining
differentials on the classes involving
$h_2$ and $h_2^2$.

\begin{lemma}\label{infty}
The following equalities hold:
\[
E_\infty^{0,0,13}
=
E_5^{0,0,13}
=
\langle vx^2U\rangle,
\]
\[
E_\infty^{1,3,9}
=
E_5^{1,3,9}
=
\langle h_2(ux+v)U\rangle,
\]
\[
E_\infty^{1,3,11}
=
E_5^{1,3,11}
=
\langle h_2vxU\rangle,
\]
$$
E_5^{2,6,7}=0.$$
Moreover, the class
$
h_2^2\Sq^1uU
\in E_5^{2,6,8}$
survives to $E_6$, but
$
E_7^{2,6,8}=0.$
\end{lemma}

\begin{proof}
We first calculate the $d_4$-differentials.
By \eqref{eq:AAHSS-differential},
\[
d_4:
E_4^{s,m,n}
\longrightarrow
E_4^{s-1,m-3,n+4}.
\]

The relevant differentials are
\[
d_4:
E_4^{1,3,9}
\longrightarrow
E_4^{0,0,13},
\]
\[
d_4:
E_4^{2,6,7}
\longrightarrow
E_4^{1,3,11},
\]
\[
d_4:
E_4^{1,3,11}
\longrightarrow
E_4^{0,0,15}.
\]

The differentials above are determined by
the $\Sq^4$-action on $N$.

In particular, the calculation of
\[
d_4:E_4^{1,3,9}\longrightarrow E_4^{0,0,13}
\]
shows that its image is generated by
$
x\Sq^4uU=\Sq^6uU,$
while its kernel is generated by
$
h_2(ux+v)U.$
Consequently,
\[
E_5^{1,3,9}
=
\langle h_2(ux+v)U\rangle,
\quad
E_5^{0,0,13}
=
\langle vx^2U\rangle.
\]

Next, consider the differential sequence
\[
E_4^{2,6,7}
\xrightarrow{d_4}
E_4^{1,3,11}
\xrightarrow{d_4}
E_4^{0,0,15}.
\]
By the $\Sq^1$- and $\Sq^2$-actions,
\[
E_4^{0,0,15}
=
\langle
x^2\Sq^4uU,\,
u\Sq^1uU
\rangle.
\]
The $\Sq^4$-calculation shows that
the incoming and outgoing $d_4$-differentials
leave precisely the class
$
h_2vxU
\in E_5^{1,3,11}$.
It also shows that
\[
d_4:
E_4^{2,6,7}
\longrightarrow E_4^{1,3,11}
\]
is nonzero. Hence
\[
E_5^{2,6,7}=0,
\qquad
E_5^{1,3,11}
=
\langle h_2vxU\rangle.
\]

By the bidegrees of the remaining
differentials, the classes in
\[
E_5^{0,0,13},\qquad
E_5^{1,3,9},\qquad
E_5^{1,3,11}
\]
support no further differential and
cannot be hit by one. Therefore they
survive to $E_\infty$.

We also note that
$
E_4^{1,3,12}=0,$
so it is zero that
\[
d_4:
E_4^{2,6,8}
\longrightarrow
E_4^{1,3,12}.
\]

It remains to examine
$
h_2^2\Sq^1uU
\in E_5^{2,6,8}.$
We now compute the next possible
nontrivial differential, namely
\[
d_6:
E_6^{2,6,8}
\longrightarrow
E_6^{1,1,14}.
\]

Recall that an AAHSS differential records
only the leading filtration term of a
boundary in the bar complex.
Thus, to identify $d_6$, we must examine
the higher-filtration terms of the
corresponding bar boundary.

Consider the following element of $B(N)$:
\[
\begin{aligned}
\mathcal C
={}&
\bigl(
\Sq^4\otimes\Sq^4
+\Sq^7\otimes\Sq^1
+\Sq^6\otimes\Sq^2
\bigr)
\otimes\Sq^1uU
\\
&+
\bigl(
\Sq^2\otimes\Sq^4
+\Sq^5\otimes\Sq^1
\bigr)
\otimes
(\Sq^2\Sq^1u+x\Sq^1u)U
\\
&+
\bigl(
\Sq^4\otimes\Sq^1
+\Sq^2\otimes\Sq^3
+\Sq^1\otimes\Sq^4
\bigr)
\otimes
(\Sq^4u+ux^2)U.
\end{aligned}
\]
A direct calculation of the bar differential gives
\[
\partial\mathcal C
=
\Sq^2\otimes
(u^2+x^2\Sq^2\Sq^1u)U
\pmod{N_{15}}.
\]
The leading term of $\mathcal C$ represents
$
h_2^2\Sq^1uU
\in E_6^{2,6,8}.$

On the other hand, by the $\Sq^1$ and $\Sq^2$ actions, we have
$
E_6^{1,1,14}
=
\langle h_1u^2U\rangle.$
The boundary formula above therefore
shows that
\[
d_6(h_2^2\Sq^1uU)
=
h_1u^2U\ne0.
\]
Consequently,
$
E_7^{2,6,8}=0$.
\end{proof}

We finally record the remaining classes
needed in the subsequent Adams calculation.

\begin{lemma}
The classes represented by
\[
E_4^{3,3,10},\qquad
E_4^{1,1,13},\qquad
E_4^{1,1,12},\qquad
E_4^{2,2,11}
\]
survive to $E_\infty$.
\end{lemma}

\begin{proof}
The assertion follows from the bidegree
of the AAHSS differential.
Indeed, inspection of
Table~\ref{E4smn} shows that
all possible incoming and outgoing
differentials involving the indicated
terms have zero source or target
for $r\ge4$.
Thus they survive to
$E_\infty$.
\end{proof}

\section{The class $b_0xU$ at the prime $3$}
\label{3primary}

Let
$
 \M_{\mathscr B}
=
\mm{MO}\langle 8\rangle\wedge\mm{M}\eta.$ Recall the Atiyah--Hirzebruch spectral sequence
\begin{equation}\label{AHSSforMO8Meta7}
E^2_{p,q}
=
H_p\bigl(
\mm{M}\eta;
\pi_q(\mm{MO}\langle 8\rangle)
\bigr)
\Longrightarrow
\pi_{p+q}(\M_{\mathscr B}).
\end{equation}

In this section we mainly show that 
$$\Z_3\subset E^2_{2,10}$$
survives to $E^\infty$.
Indeed we complete the proof by considering the mod $3$ Adams spectral sequence
\[
E_2^{s,t}
=
\mm{Ext}_{\mathcal A}^{s,t}
\bigl(H^*(\M_{\mathscr B};\mb{Z}_3),\mb{Z}_3\bigr)
\Longrightarrow
\pi_{t-s}(\M_{\mathscr B})_{(3)},
\]
where $\mathcal A$ denotes the mod $3$ Steenrod algebra. 

For the Adams $E_2$-terms, we have the AAHSS
$$E_1^{s,m,n}=\mm{Tor}_{\mathcal A}^{s,s+m}
\bigl(H^*(\MO \langle 8\rangle ;\mb{Z}_3),\mb{Z}_3\bigr)\otimes H^n(\M\eta;\Z_3)\Longrightarrow E_2^{s,s+m+n}.$$ 
 We first determine $\mm{Tor}_{\mathcal A}^{s,s+m}
\bigl(H^*(\MO \langle 8\rangle ;\mb{Z}_3),\mb{Z}_3\bigr)$.
  
From \cite{Aikawa} and the Vanishing Theorem \cite{Ravenel1986}, the generators of the group $\mm{Tor}_{\mathcal{A}}^{s,t}(\mb{Z}_3,\mathbb{Z}_3)$ are in the following table for $0<t-s<14$.
\begin{center}
\resizebox{0.9\textwidth}{!}{
\begin{tabular}{|c|c|c|c|c|c|c|c|c|c|c|c|c|c|c|c|}
\hline
$s$ & $s>0$ & $1$ & $1$ &$2$ &$2$&$2$&$3$&$3$&$3$\\
\hline
$t-s$ & 0 & $3$ & $11$ &$7$ &$10$&$11$&$10$&$11$&$13$\\
\hline
$\mm{Tor}$ &$a^s_0$& $h_{1,0}$ &$h_{1,1}$&$\Pi_0h_0$&$b_{0}$ &$a_0h_{1,1}$&$a_0b_{0}$&$a^2_0h_{1,1}$&$b_0h_{1,0}$\\
\hline
\end{tabular}
}
\end{center}
where $a_0$ is the dual of Bockstein $\beta$, $h_{1,0}$ is the dual of $\mathcal{P}^1\in \mathcal{A}$, $h_{1,1}$ is the dual of $\mathcal{P}^3$.
From \cite{Giambalvo}, $H^l(\mm{MO} \langle  8 \rangle; \mathbb{Z}_3)$ ($0< l \le 13$) is nontrivial if $l=8$ and $12$, with generators $m_8$ and $m_{12}$.
From \cite[Figure 1]{HoRa1995} and the AAHSS
$$\widetilde E_1^{s,m,n}=\mm{Tor}^{\mathcal{A}}_{s,s+m}(\mathbb{Z}_3, \mb{Z}_3)\otimes  H^n(\mm{MO} \langle  8 \rangle;\mathbb{Z}_3)$$
$$\Longrightarrow \mm{Tor}^{\mathcal{A}}_{s,s+m+n}(\mb
{Z}_3, H^\ast(\mm{MO} \langle  8 \rangle;\mathbb{Z}_3)),$$
 the nontrivial $\mm{Tor}_{\mathcal{A}}^{s,s+l}(H^\ast(\mm{MO} \langle  8 \rangle;\mathbb{Z}_3),\mb{Z}_3)$ are as follows for $0\le l\le 14$
\begin{center}
\resizebox{1.0\textwidth}{!}{
\begin{tabular}{|c|c|c|c|c|c|c|c|c|c|c|c|c|c|c|c|}
\hline
$s$ &$s\ge 0$ & $1$ & $1$ &$2$ &$2$& $3$&$s\ge3$&$3$\\
\hline
$t-s$ &$8$ & $3$ & $11$ &$7$ &$10$&$10$&$12$ &$13$\\
\hline
Tor &$a_0^sm_8$ & $h_{1,0}$ &$h_{1,0}m_8$&$\Pi_0h_0$&$b_{0}$ & $a_0b_{0}=h_{1,0}\Pi_0h_0$&$a_0^sm_{12}$&$b_0h_{1,0}$\\
\hline
\end{tabular}
}
\end{center}

Besides, $d_2(m_8)=\Pi_0h_0U$ \cite{HoRa1995} in the ASS for $\pi_*(\mm{MO} \langle  8 \rangle )$.

\begin{lemma}\label{Mod3_ASS_013}
	The Adams $E_2^{0,13}$-term is $\Z_3$ and survives to $E_\infty$ in the ASS.
\end{lemma}
\begin{proof}
	By Theorem \ref{image} and Lemma \ref{detectbyorientableclass}(3), the edge homomorphism of Adams spectral sequence gives $E_\infty^{0,13}\ne 0$. This implies 
	\begin{equation}\label{Ext_013_mod3}
		\mm{Ext}_{\mathcal A}^{0,13}
\bigl(H^*(\M_{\mathscr B};\mb{Z}_3),\mb{Z}_3\bigr)\ne 0.
	\end{equation}
	
	By Lemma \ref{3cohomology} and $\mathcal P^1U=p_1(\eta)\pmod 3= kx^2\pmod 3$, we have
	$$\mathcal P^1(z_9U)=z_{13}+kz_9x^2.$$
	Then the AAHSS differential 
	$$d_4:E_4^{1,3,9}\to E_4^{0,0,13}=E_2^{0,0,13}=\Z_3\oplus \Z_3$$
	kills a $\Z_3$-summand. We chose $z_9x^2U$ as the generator of $E_5^{0,0,13}$. Therefore, $E_\infty^{0,0,13}=E_5^{0,0,13}$ follows by \eqref{Ext_013_mod3}. And the $\mm{Ext}_{\mathcal A}^{0,13}$-group is $\Z_3$.  
\end{proof}

\begin{lemma}
	\label{b0xASS}
The class
$
b_0xU$ denote a nonzero class in the Adams $E_2^{2,14}$-term. It
survives to $E_\infty$ in the mod $3$ Adams spectral
sequence for $\M_{\mathscr B}$.
\end{lemma}

\begin{proof}
We first prove that $b_0xU$ is a nonzero class in the Adams $E_2^{2,14}$-term.
In the AAHSS for the Adams $E_2$-terms, $b_0xU\in E_1^{2,10,2}$. 
By the degree reason, all incoming differentials for $E_r^{2,10,2}$ are zero. The targets of possible outgoing  differentials are $E_8^{1,3,10}$ and $E_{11}^{1,0,13}$.

By Lemma \ref{P1betaz9}, the differential $d_4:E_4^{1,3,10}\to E_4^{0,0,14}$ is nontrivial. Thus $$E_8^{1,3,10}=E_5^{1,3,10}=0.$$
By Lemma \ref{3cohomology}, the Bockstein $\beta:H^{13}(\M\eta;\Z_3)\to H^{14}(\M\eta;\Z_3)$
is an isomorphism. So, the differential $d_1:E_1^{1,0,13}\to E_1^{0,0,14}$ is nontrivial. Thus $$E_{11}^{1,0,13}=E_{2}^{1,0,13}=0.$$

The above equations implies that all outgoing differentials for $E_r^{2,10,2}$ are zero. Therefore, $E_\infty^{2,10,2}=E_1^{2,10,2}$, i.e.,
$b_0xU\in \mm{Ext}_{\mathcal A}^{2,14}
\bigl(H^*(\M_{\mathscr B};\mb{Z}_3),\mb{Z}_3\bigr).
$

By checking the $E_1$-terms of
the AAHSS, we have 
\[
\mm{Ext}_{\mathcal A}^{s,s+11}
\bigl(H^*(\M_{\mathscr B};\mb{Z}_3),\mb{Z}_3\bigr)=0
\quad\text{for }s\ge4.
\]
So all outgoing Adams differentials of $b_0xU$ are zero. 

By the degree of
 Adams differential, the possible incoming Adams differential for $b_0xU$ is 
$d_2:E_2^{0,13}
\longrightarrow
E_2^{2,14}$. By Lemma \ref{Mod3_ASS_013}, this differential vanishes. 
Therefore $b_0xU$ survives in the ASS.
\end{proof}

\begin{lemma}\label{b0x2}
In the AHSS \eqref{AHSSforMO8Meta7}, $\Z_3\subset E^2_{2,10}$ survives to $E^\infty$.
\end{lemma}
\begin{proof}
Let
$
\M_{\mathscr B}|_1
=
\mm{MO}\langle8\rangle
\wedge
\bigl(\mm{M}\eta|_{\CP^1}\bigr),$
and let
\[
i:\M_{\mathscr B}|_1\longrightarrow\M_{\mathscr B}
\]
be the map induced by the inclusion
$\CP^1\hookrightarrow\CP^\infty$.

By the AAHSS computation, the Adams $E_2$-terms for $\M_{\mathscr B}|_1$
satisfies
\[
0\ne b_0xU\in E_2^{2,14}(\M_{\mathscr B}|_1),
\]
together with
$
E_2^{0,13}(\M_{\mathscr B}|_1)=0,$
$E_2^{s,s+11}(\M_{\mathscr B}|_1)=0$
for $s\ge4.$
The degree argument shows that $b_0xU$
survives in the ASS for $\M_{\mathscr B}|_1$.

Let
$
\gamma_1\in\pi_{12}(\M_{\mathscr B}|_1)_{(3)}$
be an element detected by this class.
By naturality,
$\gamma= i_*(\gamma_1)$ is detected by $b_0xU$ in
the ASS for $\M_{\mathscr B}$.

Since $\CP^1$ has cells only in dimensions
$0$ and $2$, we have $\mm{AH}(\gamma_1)\le 2$.
Naturality consequently gives $\mm{AH}(\gamma)\le 2$. Next we show $\mm{AH}(\gamma)= 2$.

It remains to exclude filtration $0$.
The mod $3$ Adams computation for
$\mm{MO}\langle8\rangle$ in
\cite{HoRa1995} shows that the classes in
stem $12$ occur only in Adams filtration
at least $3$; indeed, they are
represented by the classes
$
a_0^s m_{12},$ $s\ge3.$

Thus every element in the image of
\[
\pi_{12}(\mm{MO}\langle8\rangle)_{(3)}
\longrightarrow
\pi_{12}(\M_{\mathscr B})_{(3)}
\]
has Adams filtration at least $3$, by
naturality of the ASS. Moreover, by the edge homomorphism of AHSS, the image is isomorphic to the $E^\infty_{0,12}$-term of the AHSS \eqref{AHSSforMO8Meta7}. 
However, $\mm{Ad}(\gamma)=2$.
Consequently, $\mm{AH}(\gamma)\ne 0$. 

Therefore, $\mm{AH}(\gamma)= 2$. This implies that $(E^\infty_{2,10})_{(3)}=(E^2_{2,10})_{(3)}=\Z_3$.
\end{proof}

\section{The spectral sequences for certain spectra}
\begin{proposition}\label{MetaASS_h12vxU}
In the Adams spectral sequence
\[
\Ext_{\mathscr A_1}^{s,t}
\bigl(
H^*(\mm M\eta;\mathbb Z_2),\mathbb Z_2
\bigr)
\Longrightarrow
{}_2\pi_{t-s}(ko\wedge \mm M\eta),
\]
$
h_1^2vxU
\in
\Ext_{\mathscr A_1}^{2,15}
\bigl(
H^*(\mm M\eta;\mathbb Z_2),\mathbb Z_2
\bigr)$ survives to $E_\infty$.
\end{proposition}

\begin{proof}
We first use the AAHSS to calculate the Ext-groups. 
\[
E_1^{s,m,n}
=
\Tor^{\mathscr A_1}_{s,s+m}
(\mathbb Z_2,\mathbb Z_2)
\otimes
H^n(\mm M\eta;\mathbb Z_2)
\Longrightarrow
\Tor^{\mathscr A_1}_{s,t}
(\mathbb Z_2,H^*(\mm M\eta;\mathbb Z_2)),
\]
whose differentials have the form
$d_r:
E_r^{s,m,n}
\longrightarrow
E_r^{s-1,m+1-r,n+r}$. One can refer to \cite[Theorem 3.1.25, p.65]{Ravenel1986} for $\Tor^{\mathscr A_1}_{\ast,\ast}
(\mathbb Z_2,\mathbb Z_2)
$.

We first consider the class
$h_1^2vxU.$ The possible outgoing differential is
$$d_2:E_2^{2,2,11}\to E_2^{1,1,13}.$$
By $w_2(\eta)\ne 0$ and the $\Sq^2$ operation
(cf. Table~\eqref{HBZ2},\eqref{SqBtable}),
we have
\[
\Sq^2(vxU)=0.
\]
Hence $d_2(h_1^2vxU)=0$.
  The possible incoming differentials are
excluded by degree reasons together with the vanishing
\[
E_2^{3,4,8}=0,\qquad E_2^{s,0,14}=0\quad s>0.
\]
Therefore,
$
h_1^2vxU$
survives to
$
E_\infty^{2,2,11}.
$ 

By the $\Sq^1$ and $\Sq^2$ operations,
$E_3^{0,0,14}=0$ in the AAHSS.

Next we consider the Adams differentials.

The map $\alpha_8:\MO\langle 8\rangle \to ko$ induces a homomorphism
$$\alpha_{8\ast}: \Ext_{\mathscr A_2}^{s,t}
\bigl(
H^*(\mm M\eta;\mathbb Z_2),\mathbb Z_2
\bigr)\to \Ext_{\mathscr A_1}^{s,t}
\bigl(
H^*(\mm M\eta;\mathbb Z_2),\mathbb Z_2
\bigr)$$
that maps $h_1^2vxU$ to $h_1^2vxU$. In the ASS for $\bord_{13}(\eta)$, $h_1^2vxU$ is a permanent cycle. Hence, in the ASS for $ko_\ast(\M\eta)$, all outgoing Adams differentials of $h_1^2vxU$ are zero.

$E_3^{0,0,14}=0$ deduces $\Ext_{\mathscr A_1}^{0,14}
\bigl(
H^*(\mm M\eta;\mathbb Z_2),\mathbb Z_2
\bigr)=0$. This implies that the incoming Adams differential $d_2$ for $h_1^2vxU$ is zero, and finishes the proof. 
\end{proof}

\begin{corollary}\label{AHSS_Meta_112}
	The coset $\{h_1^2vxU\}\subset ko_{13}(\M\eta)$ has AH filtration $11$. 
\end{corollary}
\begin{proof}
	The induced homomorphism $\alpha_{8\ast}:\bord_{13}(\eta)\to ko_{13}(\M\eta)$ sends the coset $\{h_1^2vxU\}$ to $\{h_1^2vxU\}$. The former is provided by $E^\infty_{11,2}(\MO\langle 8\rangle ,\M\eta)$ in the AHSS for $\bord_{\ast}(\eta)$. On the other hand, the induced homomorphism 
	$$E^\infty_{11,2}(\MO\langle 8\rangle ,\M\eta)\to E^\infty_{11,2}(ko ,\M\eta)$$
	is injective. Hence, any element of the later $\{h_1^2vxU\}$ has AH filtration at most $11$. Then applying the argument in the proof for Proposition \ref{221} gives that any element of the later $\{h_1^2vxU\}$ has AH filtration at least $11$. Therefore, the desired result follows.
	
	The final key is to check $\Ext_{\mathscr A_1}^{s,s+13}
\bigl(H^\ast(\mm{M}\eta_{10};\Z_2),\Z_2\bigr)=0$ for $s\le 2$, where $\mm{M}\eta_{10}$ is the $10$-skeleton of $\mm{M}\eta$. Here we also use the AAHSS. 
\end{proof}

\begin{proposition}\label{Mvarphih12b9zeta}
The class
\[
h_1^2\beta_9x U
\in
\Ext_{\mathscr A_1}^{2,15}
\bigl(
H^*(\mm M\varphi;\mathbb Z_2),\mathbb Z_2
\bigr)
\]
survives to $E_\infty$ in the Adams spectral sequence
\[
\Ext_{\mathscr A_1}^{s,t}
\bigl(
H^*(\mm M\varphi;\mathbb Z_2),\mathbb Z_2
\bigr)
\Longrightarrow
{}_2\pi_{t-s}(ko\wedge\mm M\varphi).
\]
\end{proposition}

\begin{proof}
We use the AAHSS to compute the Adams $E_2$-term.
Its $E_1$-page is
\[
E_1^{s,m,n}
=
\Tor_{\mathscr A_1}^{s,s+m}
(\mathbb Z_2,\mathbb Z_2)
\otimes
H^n(\mm M\varphi;\mathbb Z_2),
\]
with differentials
$
d_r:
E_r^{s,m,n}
\longrightarrow
E_r^{s-1,m+1-r,n+r}.$

We first show that
$
h_1^2\beta_9x U
\in E_1^{2,2,11}
$
survives to $E_\infty$ in the AAHSS.

By $w_2(\varphi)\ne 0$ and \eqref{SqBq}, we  obtain
$
\Sq^2(\beta_9x U)=0.$

Consequently, the outgoing differential
$
d_2:
E_2^{2,2,11}
\longrightarrow
E_2^{1,1,13}$
vanishes. All other outgoing
differentials have zero targets by degree
considerations.

On the other hand, every possible incoming
differential originates from a term of the form
$
E_r^{3,m,12-m}.$
By degree
considerations, we have
\[
E_1^{3,m,12-m}=0
\qquad(m\ge0).
\]
Thus all incoming differentials vanish as well.
We conclude that
\[
E_\infty^{2,2,11}
=
E_1^{2,2,11}
=
\langle h_1^2\beta_9x U\rangle.
\]

We next determine the Adams
differentials for this class. 

Recall that an Adams differential has bidegree
\[
d_r:E_r^{s,t}
\longrightarrow
E_r^{s+r,t+r-1}.
\]

Since
$
H^{14}(\mm M\varphi;\mathbb Z_2)=0,$
we have
$
E_1^{0,0,14}=0,$
and consequently
\[
\Ext_{\mathscr A_1}^{0,14}
\bigl(
H^*(\mm M\varphi;\mathbb Z_2),
\mathbb Z_2
\bigr)=0.
\]
This excludes the only possible incoming
Adams differential to
$E_2^{2,15}$.

On the other hand, the target of a differential originating from
$ E_r^{2,15}$ is
$
E_r^{2+r,14+r}$ for $r\ge2$, which has 
Adams filtration at least $4$.

We therefore examine
\[
\Ext_{\mathscr A_1}^{s,s+12}
\bigl(
H^*(\mm M\varphi;\mathbb Z_2),
\mathbb Z_2
\bigr),
\qquad s\ge4.
\]

By the degrees and relations of
$\Ext_{\mathscr A_1}^{*,*}
(\mathbb Z_2,\mathbb Z_2)$,
the only possibly nonzero terms
$
E_1^{s,m,12-m},$ $s\ge4,$
occur when
$
m=8, 10, 12.$

More precisely, the relevant generators are
\[
h_0^j\omega x^2U,\qquad
h_1^2\omega x U,\qquad
h_0^j\nu\omega U,
\]
where $j\ge0$.

The Steenrod action in these low degrees satisfies
\[
\Sq^1U=0,
\qquad
\Sq^2(x U)=0,
\qquad
\Sq^2(x^2U)=x^3U=0.
\]
Here the last equality uses the fact that
$x^3$ has odd order.

Together with the degree restrictions,
the above formulas show that the relevant
AAHSS differentials vanish. Thus
\begin{equation}\label{E_12_varphi}
\begin{aligned}
E_\infty^{j+7,12,0}
&=
E_1^{j+7,12,0}
=
\langle h_0^j\nu\omega U\rangle,
\\
E_\infty^{j+4,8,4}
&=
E_1^{j+4,8,4}
=
\langle h_0^j\omega\zeta^2U\rangle,
\\
E_\infty^{6,10,2}
&=
E_1^{6,10,2}
=
\langle h_1^2\omega\zeta U\rangle,
\end{aligned}
\qquad j\ge0.
\end{equation}

In particular, \eqref{E_12_varphi}
describes all possible targets of outgoing
Adams differentials from
$h_1^2\beta_9x U$.
Thus we want to show that none of
these classes can be hit by an Adams
differential.

We now use the Thom spectrum map
\[
c:\mm M\varphi\longrightarrow\mm M\xi,
\]
where $\xi=-kL$ is the
bundle over $\CP^\infty$, and $c$ is introduced in Section \ref{oriented}.

The induced morphism of Adams spectral
sequences allows us to compare the classes
in \eqref{E_12_varphi} with the corresponding
classes for $\mm M\xi$.

We first compute the relevant Adams
$E_2$-terms for $\mm M\xi$.
The same AAHSS calculation gives
\begin{equation}\label{E_12_xi}
\begin{aligned}
E_\infty^{j+7,12,0}
&=
E_1^{j+7,12,0}
=
\langle h_0^j\nu\omega U\rangle,
\\
E_\infty^{j+4,8,4}
&=
E_1^{j+4,8,4}
=
\langle h_0^j\omega x^2U\rangle,
\\
E_\infty^{6,10,2}
&=
E_1^{6,10,2}
=
\langle h_1^2\omega xU\rangle,
\end{aligned}
\qquad j\ge0.
\end{equation}

We claim that  none of
these classes in
\eqref{E_12_xi} can be hit by an Adams
differential
in the ASS for
$ko\wedge\mm M\xi$.

By the degrees of the
$\mathscr A_1$-Ext generators,
the only nonzero $E_1$-terms
with $m+n=13$ are
\[
E_1^{1,1,12}
=
\mathbb Z_2\langle h_1x^6U\rangle,
\quad E_1^{5,9,4}
=
\mathbb Z_2\langle h_1\omega x^2U\rangle.
\]

Since $w_2(\xi)=x$, we have
$
\Sq^2(x^6U)=x^7U$
and
$
\Sq^2(x^2U)=x^3U.$
Consequently, the AAHSS differentials
\[
d_2(h_1x^6U)=x^7U,
\quad
d_2(h_1\omega x^2U)
=
\omega x^3U
\]
are nonzero.
Thus both classes disappear at the
$E_3$-page.
Therefore
\[
\Ext_{\mathscr A_1}^{s,s+13}
\bigl(
H^*(\mm M\xi;\mathbb Z_2),
\mathbb Z_2
\bigr)
=0
\qquad(s\ge0).
\]
This  implies the claim.

We now use the Thom map $c$ to induce a morphism
of ASSs
\[
\Ext_{\mathscr A_1}^{s,t}
\bigl(
H^*(\mm M\varphi;\mathbb Z_2),
\mathbb Z_2
\bigr)
\longrightarrow
\Ext_{\mathscr A_1}^{s,t}
\bigl(
H^*(\mm M\xi;\mathbb Z_2),
\mathbb Z_2
\bigr).
\]
 Under the induced map on the
Adams $E_2$-terms, the generators in
\eqref{E_12_varphi} map to the
corresponding generators in
\eqref{E_12_xi}.
Combining the claim with
naturality gives that no
class in \eqref{E_12_varphi} can be
the target of an Adams differential.
So, $d_r(h_1^2\beta_9x U)=0$ for $r\ge 2$. This finishes the proof.
\end{proof}

\begin{proposition}\label{E2112h12beta9zeta}
The coset
$
\{h_1^2\beta_9x U\}
\subset ko_{13}(\mm M\varphi)
$
has AH filtration $11$.
\end{proposition}

\begin{proof}
Let
$\mm M\varphi_{11}$
be the \(11\)-skeleton of
\(\mm M\varphi\).

The class
$h_1^2\beta_9x U$
only involves cohomology classes of degree at most \(11\).
Hence the AAHSS calculation for
\(\mm M\varphi_{11}\) gives a nonzero class
\[
h_1^2\beta_9x U
\in \Ext_{\mathscr A_1}^{2,15}
\bigl(
H^*(\mm M\varphi_{11};\mathbb Z_2),
\mathbb Z_2
\bigr).
\]
Furthermore, the same argument gives that $h_1^2\beta_9x U$ survives to $E_\infty$ in the ASS for $ko \wedge \M\varphi_{11}$.

On the other hand,
$
E_2^{p,13-p}
(ko,\mm M\varphi_{11})
=
0$
$(p=12,13)$,
because the \(11\)-skeleton has no homology in these dimensions.
Therefore, any element in the coset $\{h_1^2\beta_9x U\}\subset ko_{13}( \M\varphi_{11})$ has
 AH filtration at most $11$.
Thus naturality implies that any element in the coset $\{h_1^2\beta_9x U\}\subset ko_{13}( \M\varphi)$ has AH filtration at most $11$.

On the other hand,
the Thom map $\M\varphi \to \M\eta$ maps the coset $\{h_1^2\beta_9x U\}$ to $\{h^2_1 vxU\}$. By Corollary \ref{AHSS_Meta_112}, the later coset has AH filtration $11$.
 Thus, naturality gives that any element in $\{h_1^2\beta_9x U\}$ has AH filtration at least $11$. 
 
 This finishes the proof.  
\end{proof}

\section*{Declaration on the Use of AI Tools}

During the preparation of this manuscript, the authors used
ChatGPT (OpenAI) as an auxiliary tool for improving the English
language, including grammar correction, refinement of mathematical
exposition, and enhancement of clarity and readability.

The research problems, central ideas, and main conceptual framework
of this work were conceived by the authors. The authors have
critically examined and verified all mathematical arguments and
results presented in the manuscript and take full responsibility
for the correctness and integrity of the work.

\end{document}